\documentclass[a4paper,12pt]{article}
\usepackage{amsmath,latexsym,amssymb,amsfonts,natbib}
\usepackage[a4paper,hmargin=1.87 cm]{geometry}
\usepackage[dvips]{graphicx}
\usepackage{pifont}
\usepackage{pstcol,pst-fill,pst-grad}
\usepackage{rotating}
\usepackage{subfigure}
\usepackage{graphics}
\usepackage{hyperref}
\usepackage{multirow}
\usepackage{multicol}
\usepackage{courier}
\usepackage{pstricks}
\usepackage[latin1]{inputenc}
\usepackage{indentfirst}
\usepackage{ccaption}
\usepackage{verbatim}
\usepackage{makeidx}
\usepackage{listings}
\usepackage{dsfont}
\usepackage{mathrsfs}

\newcommand{\proof}     {\paragraph{Proof}}
\newcommand{\carre}     {\hfill$\Box$}
\numberwithin{equation}{section}
\newtheorem{defi}{Definition}
\newtheorem{lem}{Lemma}
\newtheorem{theo}{Theorem}

\newtheorem{prop}{Proposition}
\newtheorem{rem}{Remark}
\newtheorem{hyp}{Hypothesis}
\newtheorem{ex}{Example}
\newcommand{\N}{\mathbb N}

\newcommand{\R}{\mathbb R}
\newcommand{\E}{\mathbb E}

\title{Weak convergence of a measure-valued process for social networks}
\author{Khader Khadraoui and Ahmed Sid-Ali\\
{\small Laval University, Department of Mathematics and Statistics}\\
{\small Qu\'ebec, Qc, G1V 0A6, Canada}\\
}
\date{}
\begin{document}
\maketitle
\begin{abstract}
This article formalizes the problem of modeling social networks into an interacting particle system on random geometric graphs. Each vertex of the graph is associated with a geometric position in a latent space describing the unobserved affinity between the network members and characterizes the strength of the interaction between them.  We endow the system with two recruitment mechanisms that depend on the positions of the particles in the latent space and one departure mechanism independent of the latent position. We characterize each spatial position by a Dirac measure and the system state by a measure-valued process which is the sum of all Dirac masses. Therefore, we investigate the large-scale behavior of the system. In particular, using a renormalization technique, we study the system's behavior when the initial number of particles goes to infinity. We thus show the weak convergence of the rescaled measure-valued process towards the solution of a deterministic integrodifferential equation. Finally, we present some Monte Carlo simulations with different parameter sets to illustrate the merit of our approach in modeling some phenomena encountered in real-world networks. 
\end{abstract}

\section{Introduction}
\label{sec0}

Let consider the closure $\bar{\mathscr{D}}_G$ of an open connected subset $\mathscr{D}_G$ of $\mathbb{R}^d$, for 
some $d\geq 1$, and by $\mathcal{B}(\bar{\mathscr{D}}_G)$ we denote the $\sigma$-algebra of Borel subsets in $\bar{\mathscr{D}}_G$. Let $(\bar{\mathscr{D}}_G,\mathcal{B}(\bar{\mathscr{D}}_G))$ be a locally compact and separable metric space, endowed with a Borel $\sigma$-field, state space. Denote by $\mathscr{P}(\bar{\mathscr{D}}_G)$ the set of probability measures on $\bar{\mathscr{D}}_G$ with weak topology. In the study of Markov models of interacting particles, one looks at a $N\in \mathbb{N}$ particle system $(x^1,\ldots, x^N)$ satisfying a stochastic dynamical equation evolving according to the infinitesimal generator that characterizes the Markov transition probability kernel on a product space $\bar{\mathscr{D}}_G^N$. The chief purpose of this work is the design of a stochastic particle system approach on finite and connected graphs for the modeling of a general continuous-time and measure-valued dynamical social network given by the empirical measure: 
\begin{align*}
   \mathscr{N}_\tau = \sum_{i=1}^{N_\tau} \delta_{x_\tau^i},   \qquad \forall \; N_\tau \in \mathbb{N}, \;\; \forall \; \tau\in[0,\infty), 
\end{align*}
where $\mathscr{N}_0\in \mathscr{P}(\bar{\mathscr{D}}_G)$. Such empirical measures and particle systems have arisen in such diverse scientific disciplines. In fact, interacting particle methods have been developed in advanced signal processing applications and particularly provide powerful tools for solving a large class of nonlinear filtering problems (see, e.g., \cite{Carv+al97, Cris+Lyons97, DelM96}), in  physics, mainly for problems of fluid mechanics (see, e.g., \cite{McKean67, Sznit84, Mel+Roell87}) and in statistical mechanics (see, e.g., \cite{Ligg85, Spit69, Dobr71}.) In recent years their application area has grown, establishing unexpected connections with some other fields such as, for instance, ecology and biology (see, e.g., \cite{Four+Mel04, Fink+al09, Tran08, Decr+al12}) and the numerical approximations of Feynman-Kac measures (see, e.g., \cite{Del+Douc14}.) For further information or an account on the general theory of interacting particle systems and measure-valued processes with a full treatment of this new class of stochastic processes, the reader can find a thorough introduction in \cite{Daw17} and the references given there.

Our primary motivation stems from the revolutionary changes in complex networks science over the past two decades, creating many new challenges and open problems, particularly in social network modeling and analysis. Guided by the well-known basic events that govern a social network, we will see that the dynamic structure of such a system can be viewed, under mild assumptions, as a measure-valued process. The aim is to take advantage of advances in stochastic analysis that make possible the modeling and analysis of social network populations having complex structures and dynamics.  This article will focus on these developments. In our formulation, the most important measure of complexity of the problem is now reduced to the time evolution of sizes $(N_\tau)_{\tau\geq 0}$ and of non-static interactions between vertices (particles) of the system. The main advantage of dealing with a measure-valued process rather than other tools used in networks literature is that
the model of $(\mathscr{N}_\tau)_{\tau\geq 0}$ is Markovian and assesses the complete characteristics of the population distributed in latent space over continuous-time as well as their dynamic interactions will be considered. Our claim that this formulation is the natural framework
for formulating a social network model will be amply
justified by the following results.

In the present paper, a new measure-valued process framework is introduced for modeling a simple social network with a latent space. By latent space, we mean here, in the spirit of \cite{Hof+all02}, that the probability of a  relation (affinity) between vertices depends on the positions of particles in an unobserved social space. In particular, we introduce a stochastic particle system approach to design our modeling. The particle system described in this work will consist of finitely many particles, and the system as a whole will be a Markov process with space given by a subset of all finite point measures on $\bar{\mathscr{D}}_G$. The model to be constructed here is obtained by imposing various types of events on the motions of the particles. To this end,  we give a rigorous pathwise representation of the model, which is summed up into three Poisson point measures. It is noteworthy that these Poisson measures will manage the incoming of new particles by receiving an invitation from an existing population member or by affiliating with existing population members on the one hand and the withdrawal of particles on the other hand. Our objective is to make the construction flexible enough to encompass linear and nonlinear
interaction functions. In some applications and real-world problems, it is desirable to use the Rayleigh fading connection function or a class of connection functions that decay exponentially in some fixed positive power of the distance between particles. In this context, we emphasize that the connection function plays a crucial role in the study of the connectivity of any network in the geometric graphs literature. For not only scratch the surface and for a detailed discussion with a complete treatment of the above quite profound question, but the reader is also referred to \cite{Pen03}  and the references given there.  After fixing the model within which we work, we solve the crucial question related to the asymptotic convergence of the empirical measure of the particle system on the path space. We prove that the rescaled random empirical measure of the
particle system weakly converges to the solution of a deterministic integrodifferential equation as the initial number of vertices grows. We mean here by rescaled the renormalization of the empirical measure, which must have the effect that the density of the vertices must grow to infinity. The proof of convergence involves essentially three steps: first, the uniqueness of the solution of the deterministic equation limit, second, the tightness of the sequence of distribution of the rescaled empirical measure, and third the convergence in distribution of this sequence. This latter equation limit establishes a parsimonious deterministic approximation to describe the system dynamics (in other words, it encapsulates the evolution of the process) when the number of vertices is large.

In our development, knowledge of complex networks theory is not a requirement. For a detailed discussion of the subject of the network, the reader is referred to the books of  \cite{Rek+Bar02} and \cite{New10} where we find a general introduction with examples that help appreciate the relevance of our assumptions (see also \cite{Durr07} for some dynamic models on graphs). The use of interacting particle systems applied to model social networks opens the door to other complex problems such as studying the connectivity (see \cite{Ahm+Kha18}) and establishing the random graph configuration of the system, with the degree distribution of the vertices will be given and with the edges will be randomly matched dynamically. The coherence of the measure-valued processes paradigm with configuration using random graphs has been recently emphasized as attested by the recent paper of  \cite{Decr+al12} about the SIR epidemic model. However, studying these graph properties (e.g., the degree distribution, the size of connected components, the density of triangles and other moments, etc.), which are of crucial interest, is beyond the scope of this paper and will be the focus of future work. We restrict our efforts here to the description of the model and the study of its asymptotic properties. To conclude this enlightenment, we emphasize that the contact process (studied on integer lattices or homogeneous trees) which is one of the most studied interacting particle systems was introduced recently in the graph theory, and probabilists started investigating this process on some families of random networks like configuration models, or preferential attachment graphs (see, e.g., \cite{Chatt+Dur09, Mout+al13, Shap+Val17}.) Our work can be seen as a new step in this direction, whereas using a measure-valued process instead of contact processes (the space here is a continuous subset of $\mathbb{R}^d$ instead of $\mathbb{Z}^d$) from the class of interacting particle systems.

The article is organized as follows. In Section \ref{sec2},  after fixing the context within which we work,  the design of our particle system approach and its key assumptions are carefully described. Moreover, the exact Monte Carlo algorithm useful for numerical computation is given for completeness. In Section \ref{sec3}, a rigorous pathwise representation of the dynamic in terms of a stochastic differential equation driven by Poisson point measures is given, and the infinitesimal generator is derived. 
In Section \ref{sec5}, we deal with the weak convergence, and we prove that the process converges in law to a deterministic equation. To validate the computational performance of the approach, we present some numerical tests in Section \ref{sec6}.  We finish the paper with a general discussion of results and the outlook in Section \ref{sec7}.

\bigskip

\section{Description of the model}
\label{sec2}

  We address a stylized interacting particle system description problem where the task is to model rigorously  the dynamic of a simple social network. We describe the graphical construction of the measure-valued process by giving in detail the model  and present its key assumptions.

\subsection{General assumptions and notations}
\label{sec2.1}

First, we define some general assumptions and notations although more notations will be introduced as needed. A geometric graph $G=(V,E)$ will be understood as a set $V$ of vertices and a set $E\subseteq\{ \{x,y\}\subseteq V : x\neq y \} $ of edges. Thus, for convenience we will not explicitly treat graphs with loops (edges that start and end at the same vertex) and parallel edges between vertices, though one can define the measure-valued process on those graphs as well and our results could then be readily adapted.  The graphs we consider will always be connected where we connect two vertices by  undirected  link if and only if their distance is smaller than a certain neighborhood radius $a_\mathsf{f}>0$ to be introduced afterwards. We denote by $N=|G|$ the number of vertices of $G$. We will often abuse notation by omitting the graph from some notations; for example, we may write $N$ in place of $N^G$ and $\mathscr{N}$ in place of $\mathscr{N}^G$. Vertices of the graph are interpreted as individuals in a social network; each individual is represented at a latent state $x\in\bar{\mathscr{D}}_G$ by the Dirac measure $\delta_x$. The measure-valued process on finite and connected graph is characterized  by the point measure:
\begin{align}
 \mathscr{N}_\tau(dx)=\sum_{i=1}^{N_\tau}\delta_{x_{\tau}^i}(dx),
\label{net}
\end{align}
where $N_{\tau}\in\N$ stands for the size of the system at time $\tau$ and $\{x_{\tau}^i,\; i=1,\ldots,N_{\tau}\}$ describes the latent states of 
vertices in $\bar{\mathscr{D}}_G$. The indexes $i$ are ordered here from an arbitrary order point of view. Note that 
representing networks by point processes has created recent advances in statistical analysis of networks, such as graphexes and
stretched graphons \citep{Veit+Roy15,Brog+al18,Car+Fox17}.  Denote by  $\mathcal{S}_F(\bar{\mathscr{D}}_G)$ the 
set of finite nonnegative measures on $\bar{\mathscr{D}}_G$. Moreover, $\mathcal{S}_G\subset\mathcal{S}_F(\bar{\mathscr{D}}_G)$ consists of the subset of all finite point measures on $\bar{\mathscr{D}}_G$:
\begin{align*}
 \mathcal{S}_G=\Big\{\sum_{i=1}^N \delta_{x^i},\;N\geq0,\;x^i\in\bar{\mathscr{D}}_G    \Big\},
\end{align*}
where by convention $\sum_{i=1}^0 \delta_{x^i}$ is the null measure. Notation (\ref{net}) designating the random system seems somewhat abstract  but it will 
be more clear in the rest of the article that $(\mathscr{N}_{\tau})_{\tau\geq0}$ is a stochastic process (precisely a measure-valued process), taking its values in $\mathcal{S}_G$ and  describing the dynamic of the interacting particle system at each time $\tau\in[0,\infty)$. Furthermore, for any measure $\mu(dx)$ defined on $\bar{\mathscr{D}}_G$ and any function  $\psi:\bar{\mathscr{D}}_G\mapsto \mathbb{R}$, we use the angle brackets $\langle \mu,\psi\rangle$ to denote the function-measure duality, i.e.,  $\langle \mu,\psi\rangle=\int_{\bar{\mathscr{D}}_G}\psi(x)\mu (dx)$. The last notation is valid for continuous measures as well as for point measure $\mathscr{N}_{\tau}(dx)$ given by (\ref{net}), in the latter case $\langle \mathscr{N}_{\tau},\psi\rangle=\sum_{i=1}^{N_{\tau}} \psi(x_{\tau}^i)$. The same notation allows us to write the size of the system at time $\tau$ as  $N_{\tau}=\langle \mathscr{N}_{\tau},1\rangle$.

To make the following description of the model clear, we present now the simple three events that will govern the system in continuous time; first the event of recruitment by invitation, second the event of recruitment by affinity and third the event of withdrawal. These late three events occur in asynchronous time. Of course, the simplicity in the modeling  makes possible the development of efficient tools for fitting the model to real network data and the improvement of its probabilistic and statistical analysis. For $\mathscr{N}\in\mathcal{S}_G$, some quantities that we 
will use in the sequel are as follows:  
\begin{itemize}
 \item Let $\alpha> 0$ denotes the invitation rate of each vertex at some state $x\in\bar{\mathscr{D}}_G$.
 \item Let $K(x,dz)$ denotes the dispersion law  of each new vertex after receiving an invitation from an existing vertex located at $x$. It is assumed to satisfy, for each $x\in\bar{\mathscr{D}}_G$,
\begin{align}
\begin{split}
 \int_{\mathscr{E}} K(x,dz)=
 \left\{\begin{array}{ll}
  1
& \textrm{ if } \mathscr{E}= \{z\in\mathbb{R}^d: x+z\in \bar{\mathscr{D}}_G\} ,
\\
 0
 & \textrm{ while if } \mathscr{E}= \{z\in\mathbb{R}^d: x+z\notin \bar{\mathscr{D}}_G\}.
   \end{array}\right.
   \end{split}
   \label{HypKern1}
\end{align}
 \item Let $w^{\mathsf{af}}(y,\mathscr{N})\in[0,\infty)$ denotes the affinity rate which describes the strength of affinity between a new vertex recruited by affinity at state $y\in \bar{\mathscr{D}}_G $ with the existing system $\mathscr{N}$. 
 \item Let $K^{\mathsf{af}}(dy)$ denotes the dispersion law of each new vertex after recruitment by affinity at some state $y\in\bar{\mathscr{D}}_G$. It
is assumed to satisfy,
 \begin{align}
  \begin{split}
 \int_{\mathscr{E}'} K^{\mathsf{af}}(dy)=
 \left\{\begin{array}{ll}
  1
& \textrm{ if } \mathscr{E}'= \bar{\mathscr{D}}_G  ,
\\
 0
 & \textrm{ while if } \mathscr{E}'\cap\bar{\mathscr{D}}_G= \emptyset.
   \end{array}\right.
    \end{split}
  \label{HypKern2}  
\end{align}
 \item For $x,y\in\bar{\mathscr{D}}_G$, let consider by $\mathsf{aff}(x,y)\in[0,\infty)$ a local affinity kernel which describes 
 the contribution of a vertex in state $x$ to the affinity affecting another vertex in state $y$.
 \item Let $\beta  > 0$ denotes the withdrawal rate of each vertex at some  $x\in\bar{\mathscr{D}}_G$.
\end{itemize}

Let us briefly discuss our first notations, it is assumed that the rates $\alpha$ and $\beta$ are space independent and are the same for all vertices just for ease of exposition. A generalization by considering functions  
$\alpha(x)$ and $\beta(x)$ with $x\in \bar{\mathscr{D}}_G$ might allow us to take into account external effects such as attractive vertices (hubs), unattractive  vertices and so forth. However, the rate $w^{\mathsf{af}}(y,\mathscr{N})$ of the affinity component of the system  depends on its state $y\in \bar{\mathscr{D}}_G$ and on the state  of the whole system $\mathscr{N}$. This affinity component will play an essential role in any connectivity study of the system under construction and its modeling should not be designed lightly.

\subsection{Interacting particle system and graphical construction}
\label{sec2.2}

Considering the state $\mathscr{N}=\sum_{i=1}^{N}\delta_{x^i}$ of the random system at a given time, to each particle located at some  $x\in\bar{\mathscr{D}}_G$ is associated an independent invitation exponential clock with parameter $\alpha>0$ and an independent withdrawal exponential clock with parameter $\beta>0$. To the whole system $\mathscr{N}$ is associated an independent affinity exponential clock with 
parameter $w^{\mathsf{af}}(y,\mathscr{N})\in[0,\infty)$ for $y\in\bar{\mathscr{D}}_G$. Specifically, the random system will be subject to 3 types of point events occurring at specific clocks where the first of all these clocks that rings determines the next event as follows: 

\begin{enumerate}
\item According to the  recruitment by invitation event, a particle in state $x$ can give an invitation to another particle to join an empty state $y\in\bar{\mathscr{D}}_G$ (in the neighborhood of $x$) and if the invitation is accepted this particle immediately becomes a member of the network. Its state is given by $y=x+z$ where $z$ is chosen randomly according to an invitation dispersion kernel $K(x,dz)$ and then the system earns new member: $\mathscr{N} \mapsto \mathscr{N}+\delta_{x+z}$. From now on, we assume that the dispersion kernel induces a density w.r.t. the Lebesgue measure on $\mathbb{R}^d$ such that this density is given by $K(x,dz)=k(x,z)dz$. 
\item When a recruitment by affinity event occurs, a new particle chooses a region to join in $\bar{\mathscr{D}}_G$ and affiliates with an existing vertices. Particularly, its state $y\in\bar{\mathscr{D}}_G$  is chosen randomly according to an affinity dispersion kernel $K^{\mathsf{af}}(dy)$ and we assume again that this kernel admits a density on $\bar{\mathscr{D}}_G$ such that $K^{\mathsf{af}}(dy)=k^{\mathsf{af}}(y)dy$. Moreover, for all $\mathscr{N}\in \mathcal{S}_G$, we assume that 
\begin{align}
 w^{\mathsf{af}}(y,\mathscr{N})=\sum_{x\in V: x\sim y}\mathsf{aff}(x,y)=\int_{\bar{\mathscr{D}}_G}\mathsf{aff}(x,y) \mathscr{N}(dx),
 \label{model-aff}
\end{align}
where $x \sim y$ means that $x$ and $y$ are neighbors. To complete the model (\ref{model-aff}), we specify the local affinity function as, among others,  a triangular truncated (linear) function given by (see Figure \ref{Fig.aff}):
\begin{align}
\forall x,y\in\bar{\mathscr{D}}_G,\qquad  \mathsf{aff}(x,y)
  =
  \left\{\begin{array}{ll}
     A_\mathsf{f}\,
       \big(1-\frac{1}{a_\mathsf{f}}
       \,\| x-y\|  \big)^+
     & \textrm{ if }x\neq y\,,
   \\
     0
     & \textrm{ while if } x=y,
  \end{array}\right.
\label{interaction}
\end{align}
where $(\cdot)^+=\max(\cdot,0)$, the parameter $A_{\mathsf{f}} > 0$ determines the high affinity level and 
$a_{\mathsf{f}} > 0$ is a second parameter specifying the radius of the affinity zone around each particle of the system.  
\begin{figure}
\center
\setlength{\unitlength}{5cm}
\medskip
\medskip
\setlength{\unitlength}{5cm}
\begin{picture}(1,1)
\put(-0.4,-0.25){\vector(0,1){1}}
\put(-0.4,-0.25){\vector(1,0){1.3}}
\put(-0.4,0.25){\line(2,-1){1}}
\put(0.3,0.2){\vector(-1,-1){0.3}}
\put(0.31,0.2){affinity area}
\put(-0.5,-0.25){$0$}
\put(0.9,-0.25){$\|x-y\|$}
\put(0.55,-0.32){$a_\mathsf{f}$}
\put(-0.5,0.8){$\mathsf{aff}(x,y)$}
\put(-0.55,0.25){$A_\mathsf{f}$}
\end{picture}
\medskip
\medskip
\medskip
\medskip
\medskip
\medskip
\caption{The triangular local affinity function given by (\ref{interaction}).}
\label{Fig.aff}
\end{figure}
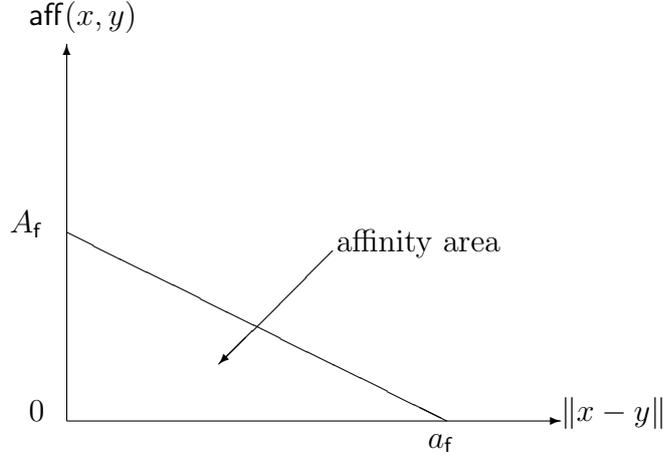
This affinity interaction is inspired from some real networks where member highly connected (with high degree) is logically highly attractive. Furthermore, this choice allows us to manage the  interaction by affinity in the system, one could imagine more complex affinity models when more complex scenarios are investigated. Then, after each affinity recruitment the system earns new vertex: $\mathscr{N} \mapsto \mathscr{N}+\delta_{y}$. 
\item  When a withdrawal event occurs, a vertex in state $x\in\bar{\mathscr{D}}_G$ disappears at rate $\beta$. This withdrawal is intrinsic as it does not depend on the state of the system. After this event the size of the system jumps instantaneously from $N$ to $N-1$ and then the system loses one vertex: $\mathscr{N} \mapsto \mathscr{N}-\delta_{x}$. 
\end{enumerate}

\begin{figure}
\center
\subfigure[\it $\mathscr{N}_\tau$ with $N_\tau=198$.]
    {\includegraphics[width=8cm,height=8cm]
                 {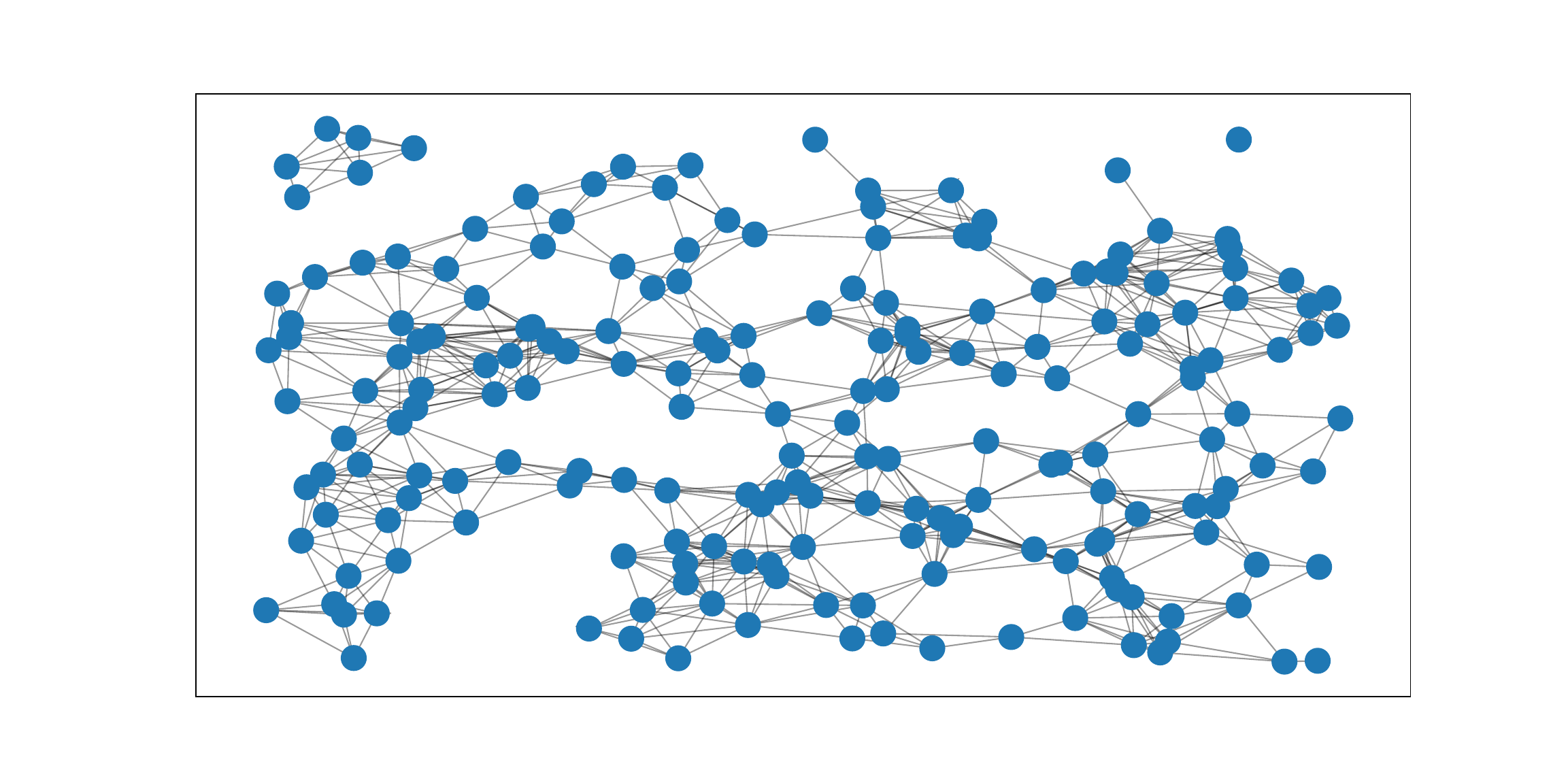}}
\subfigure[\it $\mathscr{N}_{\tau+\Delta\tau}$ with $N_{\tau+\Delta\tau}=197$.]
    {\includegraphics[width=8cm,height=8cm]
                 {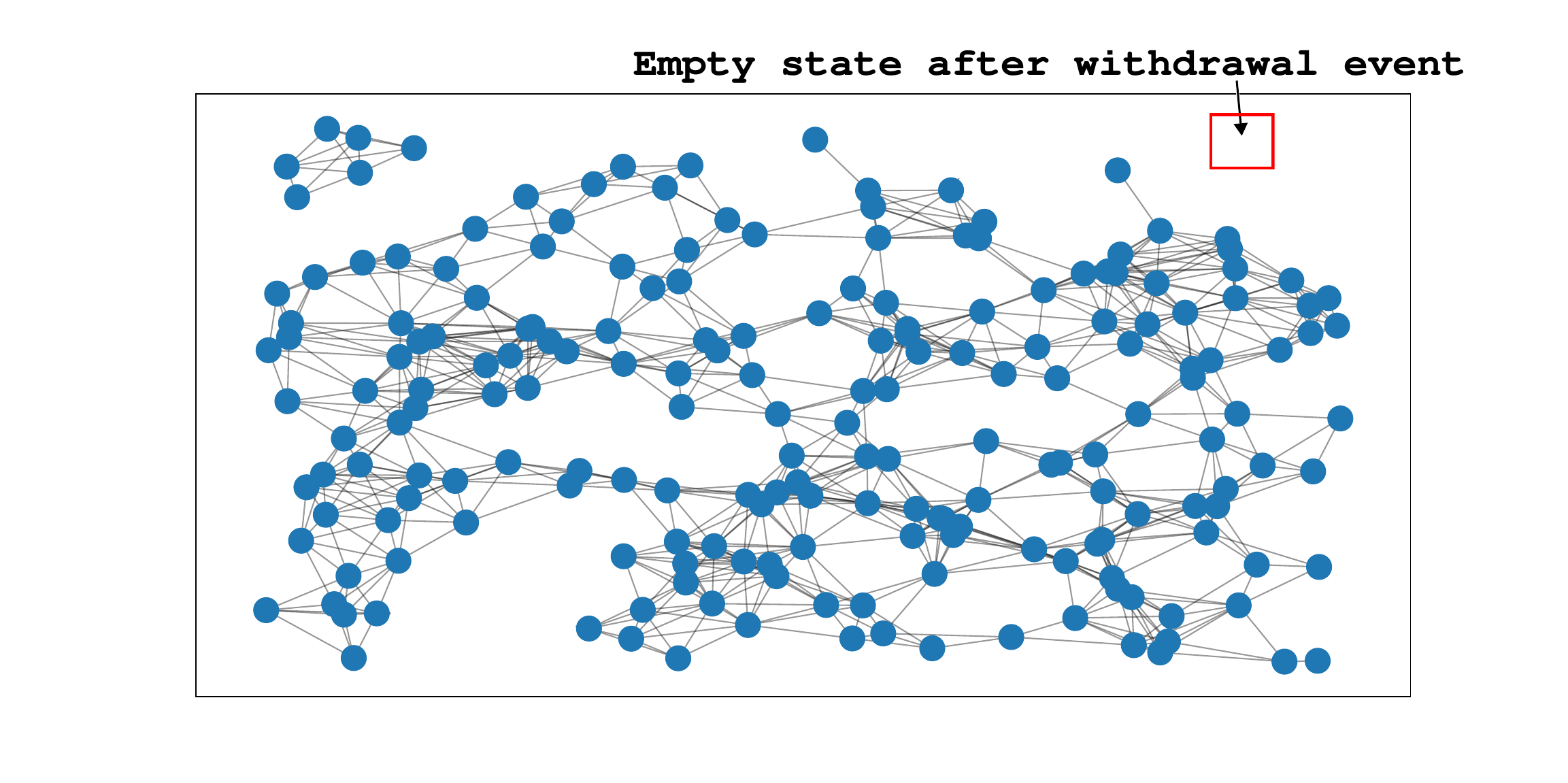}}
\subfigure[\it $\mathscr{N}_{\tau'+\Delta\tau'}$ with $N_{\tau'+\Delta\tau'}=198$.]
    {\includegraphics[width=8cm,height=8cm]
                 {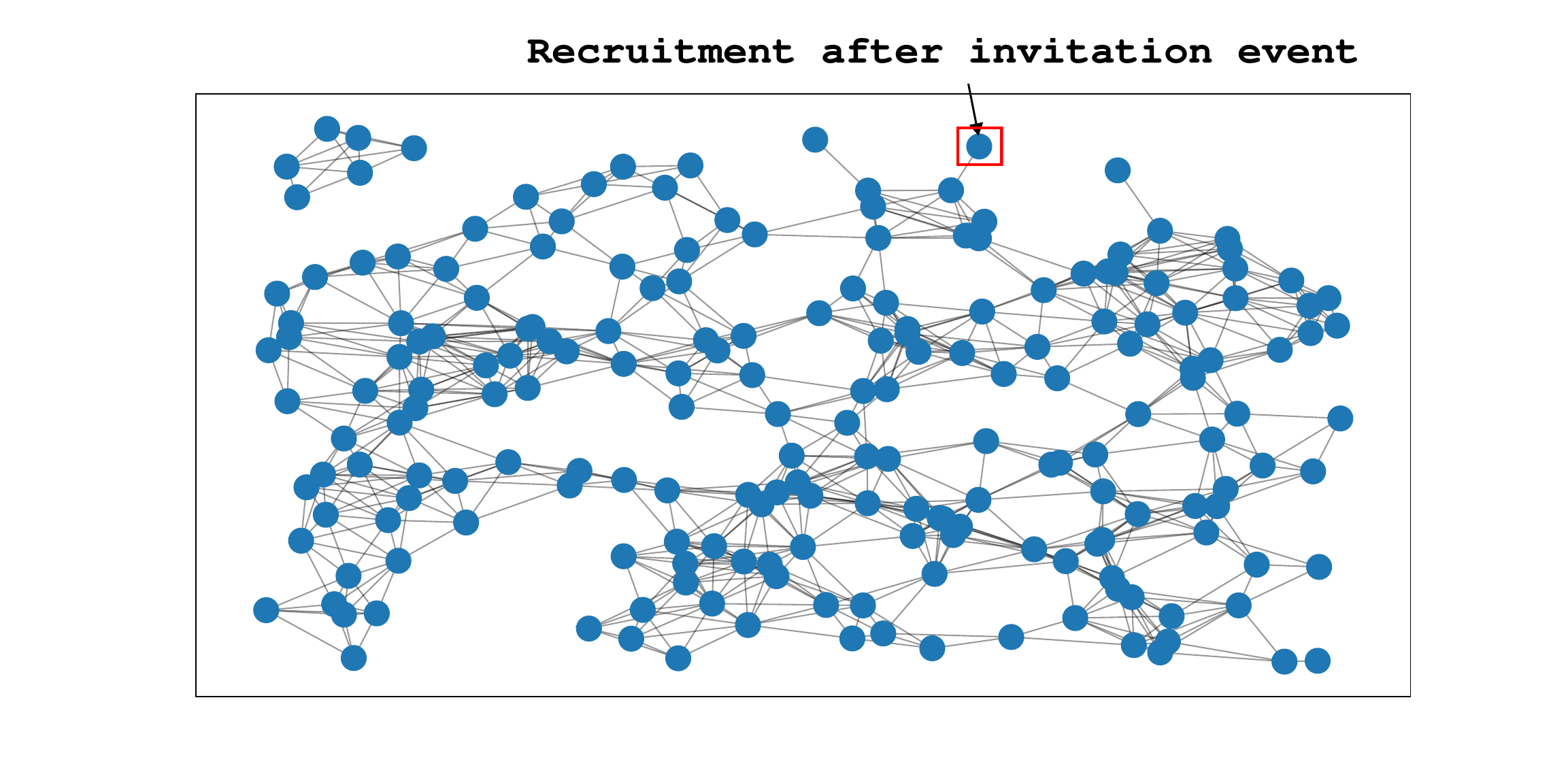}}  
 \subfigure[\it $\mathscr{N}_{\tau"+\Delta\tau"}$ with $N_{\tau"+\Delta\tau"}=199$.]
    {\includegraphics[width=8cm,height=8cm]
                 {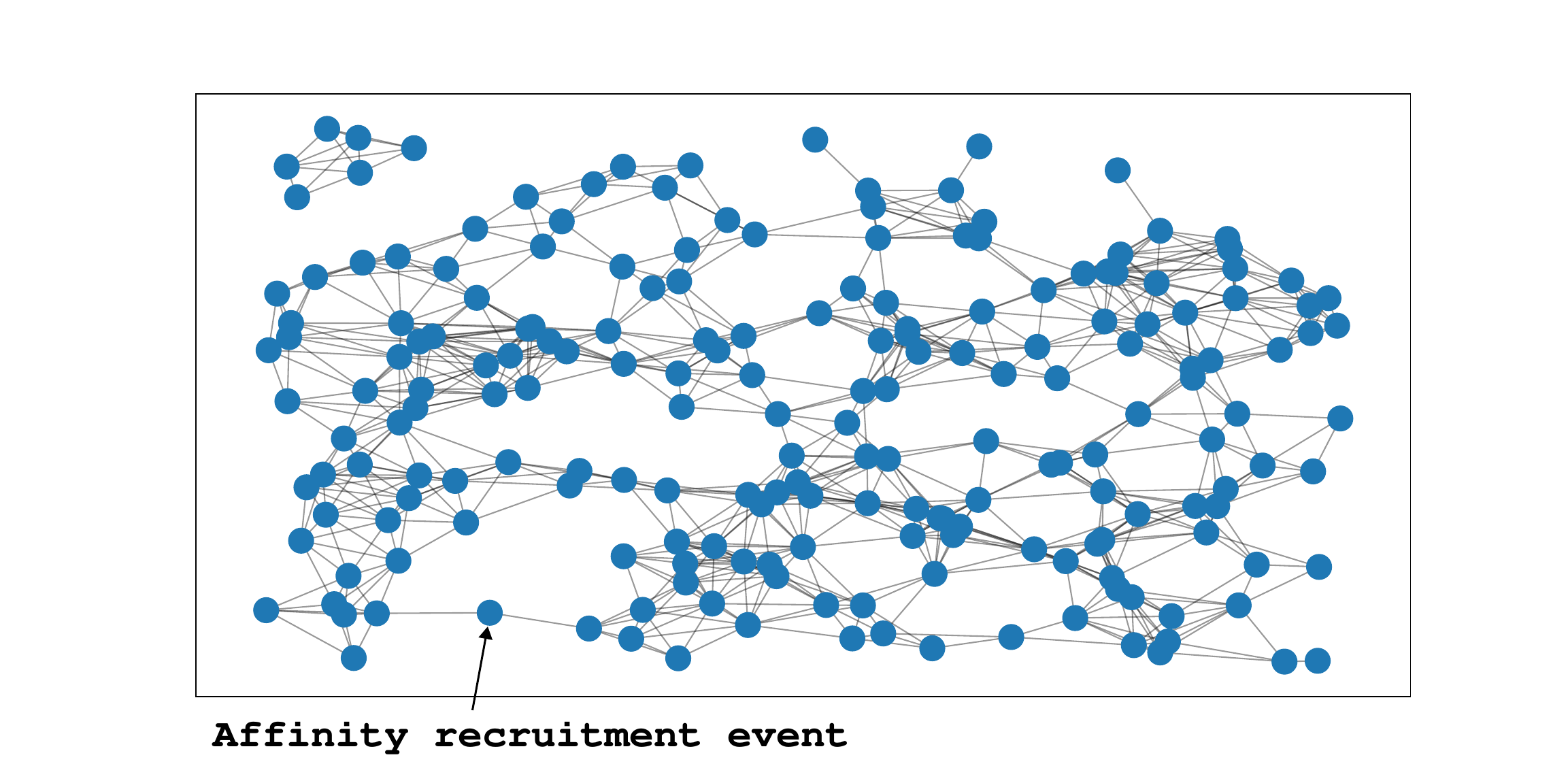}} 
\caption{{\small Geometric graphs shown interacting process in $\bar{\mathscr{D}}_G=[0,1]^2$, here $a_\mathsf{f}=0.125$. (a) Particles indicated by 198 small circles (vertices) connected by linear lines (edges) if the interpoint distance is small than 0.125. (b) A particle is chosen, and incurs a withdrawal event at time $\tau'=\tau+\Delta\tau$ where $\Delta\tau$ is a random exponential time. (c) A particle is chosen, and invites another particle to join the system at its neighborhood. (d) An empty state is chosen, and a new particle joins the system by affinity.}}
\label{Fig-graph}
\end{figure}

A natural hypothesis will be considered from now on is to assume that these mechanisms of recruitment with dispersion and withdrawal are mutually independent. Nevertheless, considering individual exponential clocks is relatively cumbersome and a better efficient Monte Carlo procedure will rely on the existence of one global exponential clock that dominates all  point phenomena. That existence holds true when all the different individual clocks are uniformly bounded. For simplicity and from now on, we assume  that the spatial dependence of the introduced kernels and rates is bounded in some sense by assuming that there exist some positive reals $\gamma_1,\gamma_2$ and some probability densities $\tilde{k}$ on $\mathbb{R}^d$ and $\tilde{k}^{\mathsf{af}}$ on $\bar{\mathscr{D}}_G$ such that
\begin{align}
 k(x,z)\leq \gamma_1 \tilde{k}(z) \qquad \textrm{ and } \qquad  k^{\mathsf{af}}(y)\leq \gamma_2 \tilde{k}^{\mathsf{af}}(y), 
 \label{Ass1}
\end{align}
for all $x,y\in\bar{\mathscr{D}}_G$ and $z\in \mathbb{R}^d$. As well, for all $\mathscr{N}\in\mathcal{S}_G$, we assume that there exists a constant $A_\mathsf{f}$, introduced firstly 
in (\ref{interaction}), such that
\begin{align}
\mathsf{aff}(x,y)\leq A_\mathsf{f} \qquad \textrm{ which implies } \qquad  w^{\mathsf{af}}(y,\mathscr{N})\leq A_\mathsf{f} N.
  \label{Ass2}
\end{align}
For completeness and to write more rigorously the Monte Carlo algorithm for the system simulation, we describe now how to simulate the different events. Given a vertex chosen at random in the system,  the type of point phenomenon to be considered is determined by a sampling technique, and it is decided whether the chosen phenomenon is actually applied or not by an acceptance/rejection sampling technique.  The existence of a uniform bound avoids explosion phenomena due to the accumulation of infinitely many events at a given time. The global clock can be computed easily thanks to the properties of the exponential distribution by,
\begin{align}
\begin{split}
  \mathsf{H}_\tau= \mathsf{h}_\tau^{\alpha} +\mathsf{h}_\tau^{\mathsf{af}} +\mathsf{h}_\tau^{\beta},\qquad \textrm{ where   } \quad  
  \left\{\begin{array}{ll}
     \mathsf{h}_\tau^{\alpha}= \alpha N_\tau,
     & 
   \\
     \mathsf{h}_\tau^{\mathsf{af}}= A_{\mathsf{f}}N_\tau,
     & 
   \\
    \mathsf{h}_\tau^{\beta} = \beta N_\tau.
    &
  \end{array}\right.
\end{split}
\label{GlobH}
\end{align}
Let $T_0= 0$ and start with a randomly chosen initial state $\mathscr{N}_0$. For $k=1,2,3,\ldots$, suppose the time of the last event $T_{k-1}$ and the corresponding state of the system $\mathscr{N}_{T_{k-1}}$ given, we describe how to simulate $\mathscr{N}_{T_{k}}$ starting from $\mathscr{N}_{T_{k-1}}$.  In order to determine the instant $T_k$ when the next event could take place, we draw one realization from the exponential distribution with parameter $\mathsf{H}_{T_{k-1}}$ using (\ref{GlobH}). Then, from the instant $T_{k-1}$ to the instant $T_{k}$ of the next event, i.e.,  along the time interval $[T_{k-1},T_{k})$ nothing is happening and one iteration of the scheme is 
as follows:

\paragraph{Iteration $\mathscr{N}_{T_{k-1}} \to \mathscr{N}_{T_{k}}$:}
\begin{enumerate}
 \item Computation of the global rate $\mathsf{H}_{T_{k-1}}$ given by (\ref{GlobH}).
 \item Simulation of the next event instant:
 \begin{align*}
  T_k=T_{k-1}+\Delta T_k \qquad \textrm{ with } \qquad  \Delta T_k\sim\textrm{Exp}(\mathsf{H}_{T_{k-1}}).
 \end{align*}
 \item Computation of the system evolution between the two instants:
 \begin{align*}
   \mathscr{N}_\tau=\mathscr{N}_{T_{k-1}}  \qquad \textrm{ for  } \qquad \tau \in[T_{k-1},T_k).
 \end{align*}
\item  Computation of the probabilities:
\begin{align*} 
  \rho_k^{\alpha}
  &= \frac{ \mathsf{h}_{T_{k-1}}^{\alpha}}{\mathsf{H}_{T_{k-1}}}\,,
  &
  \rho_k^{\mathsf{af}}
  &= \frac{ \mathsf{h}_{T_{k-1}}^{\mathsf{af}}}{ \mathsf{H}_{T_{k-1}}}\,,
  &
  \rho_k^{\beta}
  &= \frac{ \mathsf{h}_{T_{k-1}}^{\beta}}{  \mathsf{H}_{T_{k-1}}}.
\end{align*}
\item  We choose at random the nature of the next event according to the probability values $\rho_k^{\alpha}$, $\rho_k^{\beta}$ and $\rho_k^{\mathsf{af}}$:
\begin{itemize}
\item With probability $\rho_k^{\alpha}$ we determine if there will be an invitation recruitment  by acceptance/rejection. If there is acceptance of the 
event, we draw uniformly a vertex $x_{T_{k-1}}^i$ where its index   
$i\sim U\{1,\dots, N_{T_{k-1}}\}$, we draw $z\in\mathbb{R}^d$ using the dispersal kernel  $K(x_{T_{k-1}}^i,dz)$ and then
\begin{align*}
 \mathscr{N}_{T_k} = \left\{\begin{array}{ll}
     \mathscr{N}_{T_{k-1}}  +\delta_{\big\{x_{T_{k-1}}^i+z\big\}}  ,
     &  \textrm{ with probability } \frac{k(x_{T_{k-1}}^i,z)}{\gamma_1 \tilde{k}(z)},
   \\
     \mathscr{N}_{T_{k-1}} ,
     &  \textrm{ with probability } 1- \frac{k(x_{T_{k-1}}^i,z)}{\gamma_1 \tilde{k}(z)}.
  \end{array}\right.
\end{align*} 
\item With probability $\rho_k^{\beta}$ we determine if a withdrawal event will occur. If acceptance, we draw uniformaly a vertex $x_{T_{k-1}}^i$ where its index $i\sim U\{1,\dots, N_{T_{k-1}}\}$ and we set $\mathscr{N}_{T_k}=\mathscr{N}_{T_{k-1}}-\delta_{\{x_{T_{k-1}}^i\}}$. 
\item With probability $\rho_k^{\mathsf{af}}$ a recruitment by affinity event occurs. If acceptance, we draw an empty state $y$ using the affinity kernel 
$K^{\mathsf{af}}(dy)$, we draw uniformly a vertex $x_{T_{k-1}}^i$ from the current system vertices   
(i.e., $i\sim U\{1,\dots, N_{T_{k-1}}\}$)  and let 
\begin{align*}
 \mathscr{N}_{T_k} = \left\{\begin{array}{ll}
     \mathscr{N}_{T_{k-1}}  +\delta_{y}  ,
     &  \textrm{ with probability } \frac{\mathsf{aff}(x_{T_{k-1}}^i,y)k^{\mathsf{af}}(y)}{{A}_{\mathsf{f}}\gamma_2\tilde{k}^{\mathsf{af}}(y)},
   \\
     \mathscr{N}_{T_{k-1}} ,
     &  \textrm{ with probability } 1-\frac{\mathsf{aff}(x_{T_{k-1}}^i,y)k^{\mathsf{af}}(y)}{{A}_{\mathsf{f}}\gamma_2\tilde{k}^{\mathsf{af}}(y)}.
  \end{array}\right.
\end{align*} 
\end{itemize}
\end{enumerate}
This simulation procedure enables us to provide some numerical tests  in Section \ref{sec6}.  
\begin{ex}
To illustrate the construction of the exact computational scheme, let us consider an example (see  Figure \ref{Fig-graph}) in which we briefly illustrate a dynamic shown respectively a withdrawal, a recruitment by invitation and a recruitment by affinity events starting from a system size of $N_\tau=198$ particles in the unit square ($\bar{\mathscr{D}}_G=[0,1]^2$) with affinity threshold $a_\mathsf{f}= 0.125$. As shown, the particles $x_\tau^1,\ldots,x_\tau^{198}$ are related through a random network and are represented by the vertices of an undirected graph. Between two neighbors, we place an edge if $\|x_\tau^i-x_\tau^j\|\leq a_\mathsf{f}$, for any $i\neq j$, to highlight the affinity interaction mechanism. 
\end{ex}

We have left open so far the precise law of the process $(\mathscr{N}_{\tau})_{\tau\geq 0}$ other than subsequently specifying a Monte Carlo
algorithm that simulates a trajectory of the system that provides an empirical representation of this law. More precisely, our previous scheme 
generates a trajectory of a process which has the same law as the considered process. As the law of the Markov process under study is characterized by its infinitesimal generator,  we shall introduce this generator further on in the sequel. The system is described by the evolution in time of the empirical measure (\ref{net}), and this evolution must follow the dynamic defined above. Therefore, we are looking for a $\mathcal{S}_G$-valued stochastic process  $(\mathscr{N}_{\tau})_{\tau\geq 0}$ with infinitesimal generator $\mathcal{L}$, defined for a large class of test functions $\Phi:\mathcal{S}_F(\bar{\mathscr{D}}_G)\mapsto\R$, of the form (unless otherwise stated)
\begin{align}
 \Phi(\mathscr{N})=F(\langle \mathscr{N},f\rangle),
 \label{Testfn}
\end{align}
where $F\in C^1(\R)$ and $f\in C(\bar{\mathscr{D}}_G)$, by:
\begin{align}
\begin{split}
\mathcal{L}\Phi(\mathscr{N})&=\alpha\int_{\bar{\mathscr{D}}_G}\mathscr{N}(dx) \int_{\mathbb{R}^d}\Big\{\Phi(\mathscr{N}+\delta_{x+z})-\Phi(\mathscr{N})\Big\}k(x,z)dz \\
&\qquad + \beta \int_{\bar{\mathscr{D}}_G} \Big\{\Phi(\mathscr{N}-\delta_{x})-\Phi(\mathscr{N})\Big\} \mathscr{N}(dx)   \\
&\qquad + \int_{\bar{\mathscr{D}}_G}\Big\{ \int_{\bar{\mathscr{D}}_G}\big( \Phi(\mathscr{N}+\delta_{y})-
\Phi(\mathscr{N})\big) \mathsf{aff}(x,y) {k}^{\mathsf{af}}(y)dy \Big\}\mathscr{N}(dx).
\end{split}
\label{generat}
\end{align}
The above generator prescribes that vertices disappear with rate $\beta$ and join the system by invitation with rate $\alpha$. The third term in (\ref{generat}) describes the recruitment by affinity component.  Note also that the three terms of the generator are linear in $\mathscr{N}$. The idea in the following is to construct a $\mathcal{S}_G$-valued stochastic process with infinitesimal generator $\mathcal{L}$. 

\begin{rem}
The class of functions $\Phi$ defined by $(\ref{Testfn})$ is separating and convergence determining  \citep[Theorem 3.2.6]{Daw93}. Therefore  it is 
sufficient to restrict ourselves to this class of functions where the expression of the generator is explicit. Moreover, the expression of the generator on 
this class of functions determines the law of the process. In particular, we refer the reader to \cite[p. 111-112]{Eth+Kur86} for the definitions of 
separating and convergence determining sets. 
\end{rem}

In the next Section \ref{sec3}, we shall give a rigorous pathwise representation of the model in terms of some Poisson point measures.

\section{Poisson measures of the particle system}
\label{sec3}

Let $\mathscr{N}_0$ denotes the initial condition of the process, it is a random variable with values in $\mathcal{S}_G $. The model  has a vector of unknown parameters, assumed to lie in $\theta=(\alpha,\beta,A_\mathsf{f},a_\mathsf{f},\theta^\alpha,\theta^{\mathsf{af}})\in (\mathbb{R}_+^*)^4 \times \mathbb{R}^{ |\theta^\alpha|+|\theta^{\mathsf{af}}|}$ where $\theta^\alpha$ denotes the parameters of the invitation dispersion kernel and $\theta^{\mathsf{af}}$ denotes the parameters of the affinity dispersion kernel. We also denote by $\mathbb{P}_\theta$ a probability measure under which the measure-valued process with parameters $\theta$ is defined on the graph $G$ (which will be clear from the context, as will the initial configuration of the process); later to avoid overly burdensome notation we omit $\theta$ from the notation as well. We denote by $\mathbb{E}_\theta$, or simply $\mathbb{E}$, the associated expectation. Let $(\Omega,\mathcal{A},\mathbb{P})$ be a 
sufficiently large probability space. On this space, we consider the objects $w^{\mathsf{af}}(y,\mathscr{N})$, $K(x,dz)$ and $K^{\mathsf{af}}(dy)$ which are assumed mutually independent and  space dependent. To manage the incoming of new vertices by invitation or affinity, and the withdrawal of vertices, we propose a stochastic differential equation (SDE) driven by Poisson point measures (PPMs) to describe the evolution of the point measure (\ref{net}). First, we need the following definition.
\begin{defi} 
We consider, on the probability space $(\Omega,\mathcal{A},\mathbb{P})$, the three independent PPMs: 
\begin{enumerate}
\item Let $P^{\alpha}$ be a PPM on $[0,\infty)\times\N^*\times\mathbb{R}^d\times[0,1]$ 
of intensity measure:
\begin{align*}
 I^{\alpha}(dt,di,dz,d\rho)=\alpha \gamma_1\tilde{k}(z)dt \Big(\sum_{k\geq 1}\delta_k(di)\Big)dzd\rho,
\end{align*}
where $dt,dz,d\rho$ are the Lebesgue measures on, respectively, $[0,\infty),\mathbb{R}^d, [0,1]$,  and $\sum_{k\geq 1}\delta_k(di)$ is the counting measure 
on $\N^*$. 
\item Let $P^{\beta}$ be a PPM on $[0,\infty)\times\N^*\times[0,1]$ 
of intensity measure:
\begin{align*}
 I^{\beta}(dt,di,d\rho)=\beta dt \Big(\sum_{k\geq 1}\delta_k(di)\Big)d\rho.
\end{align*}
\item Let $P^{\mathsf{af}}$ be a PPM on $[0,\infty)\times\N^*\times\bar{\mathscr{D}}_G\times[0,1]$ 
of intensity measure:
\begin{align*}
 I^{\mathsf{af}}(dt,di,dy,d\rho)= A_{\mathsf{f}}\gamma_2 \tilde{k}^{\mathsf{af}}(y) dt \Big(\sum_{k\geq 1}\delta_k(di)\Big)dyd\rho.
\end{align*}
\end{enumerate}
\label{defpois}
\end{defi}

The point measure $P^{\alpha}$ provides possible times at which invitation recruitment may occur. Each of its atoms is associated with a possible invitation event time $t$, a real $z$ which gives the dispersion of the vertex being possibly recruited and an integer $i$ which gives the vertex that has to produce the invitation. The mark $\rho$ is as in the previous section an auxiliary variable used for the construction of the acceptance/rejection sampling. Similar interpretation holds for $P^{\beta}$ and $P^{\mathsf{af}}$. We consider the canonical filtration $(\mathcal{F}_\tau)_{\tau \geq0}$ generated by the previous three PPMs. We aim now to write the random process in terms of these stochastic objects.  We 
shall describe the system by the evolution in time of the empirical measure $\mathscr{N}_\tau$ through its infinitesimal generator. To this end, let first from now on assume that $\bar{\mathscr{D}}_G$ is a compact and  
consider the path space $\mathbb{T}_G\subset\mathbb{D}([0,\infty),\mathcal{S}_F(\bar{\mathscr{D}}_G))$  defined by
\begin{align}
 \mathbb{T}_G=\Big\{(\mathscr{N}_\tau)_{\tau\geq0}\Big/
\begin{tabular}{l}
 $\forall \tau\geq0,\mathscr{N}_\tau \in\mathcal{S}_G, \textrm{ and } \exists  0=\tau_0<\tau_1<\tau_2<\cdots,$\\
  $\lim_{n\to\infty} \tau_n=\infty \textrm{ and } \mathscr{N}_\tau=\mathscr{N}_{\tau_i}\; \forall \tau\in[\tau_i,\tau_{i+1}) $\\
\end{tabular}
\Big\},
\end{align}
where $\mathbb{D}([0,\infty),\mathcal{S}_F(\bar{\mathscr{D}}_G))$ is the Skorokhod space of $\mathcal{S}_F(\bar{\mathscr{D}}_G)$-valued c\`adl\`ag functions on $[0,\infty)$ (more details will be given in Section \ref{sec5}). Therefore, $\mathbb{T}_G$ is the subspace of $\mathcal{S}_G$-valued pure jump processes on $[0,\infty)$. 
Note that for $(\mathscr{N}_\tau)_{\tau\geq0}\in\mathbb{T}_G$ and $\tau>0$, we can define $\mathscr{N}_{\tau-}$ by
\begin{align}
 \mathscr{N}_{\tau-}=
\left\{\begin{array}{ll}
  \mathscr{N}_\tau
& \textrm{ if } \tau \notin \cup_i \{\tau_i\} ,
\\
 \mathscr{N}_{\tau_{i-1}}
 & \textrm{ if } \tau=\tau_i \textrm{ for some } i\geq1.
   \end{array}\right.
\label{Xi}
\end{align}
In other words, for $\tau=\tau_i$, $\mathscr{N}_{\tau-}$ is the value of $\mathscr{N}$ just before the jump time $\tau_i$. The process $(\mathscr{N}_\tau)_{\tau\geq 0}$ features a jump dynamics (recruitments and withdrawals) and we therefore recall a well-known formula for the pure jump processes:
\begin{align}
 \Phi(\mathscr{N}_\tau)=\Phi(\mathscr{N}_0)+\sum_{t\leq \tau}[\Phi(\mathscr{N}_{t-}+\{\mathscr{N}_t-\mathscr{N}_{t-}\})-\Phi(\mathscr{N}_{t-})], \textrm{ a.s. for }\tau\geq 0, 
 \label{jump}
\end{align}
for any function $\Phi:\mathcal{S}_G\mapsto\mathbb{R}$. Note that the sum $\sum_{t\leq \tau}$ in (\ref{jump}) contains only a finite number of terms as the system $(\mathscr{N}_\tau)_{\tau\geq0}$ admits only a finite number of jumps over any finite time interval. In fact, the formula (\ref{jump}) includes information on the process, i.e. the dynamics between the jumps,
but no information on the jumps themselves. Hence, for this reason we have introduced the previous three PPMs to overcome the derivation
of the following SDE describes the evolution of the system for all $\tau\geq 0$
as stated in Definition \ref{defnet}.
\begin{defi}
Assume that conditions $(\ref{Ass1})$ and $(\ref{Ass2})$ are satisfied. We say that a $(\mathcal{F}_{\tau})_{\tau}$-adapted stochastic process that belongs 
a.s. to $\mathbb{T}_G$ describes our model if a.s., for all $\tau\geq 0$,
\begin{align}
\begin{split}
 \mathscr{N}_\tau=&\mathscr{N}_0+\int_{0}^\tau \int_{\N^*}\int_{\mathbb{R}^d}\int_0^1 \mathds{1}_{\{i\leq N_{t-}\}}
\mathds{1}_{\big\{\rho\leq ( k(x_{t-}^i,z))\big/(\gamma_1\tilde{k}(z)) \big\}}
  \delta_{(x_{t-}^i+z)} P^{\alpha}(dt,di,dz,d\rho)\\
&-\int_{0}^\tau \int_{\N^{\ast}}\int_0^1 \mathds{1}_{\{i\leq N_{t-}\}}\mathds{1}_{\{\rho\leq(\beta/\beta)\}}
\delta_{(x_{t-}^i)}P^{\beta}(dt,di,d\rho)\\
&+\int_{0}^\tau\int_{\N^{\ast}}\int_{\bar{\mathscr{D}}_G}\int_0^1 \mathds{1}_{\{i\leq N_{t-}\}}
\mathds{1}_{\big\{\rho\leq (\mathsf{aff}(x_{\tau-}^i,y)k^{\mathsf{af}}(y))\big/(A_\mathsf{f}\gamma_2 \tilde{k}^{\mathsf{af}}(y) )\big\}}\delta_{(y)}
P^{\mathsf{af}}(dt,di,dy,d\rho),
\end{split}
\label{netdef}
\end{align}
where $\mathds{1}_{A}$ denotes the indicator function of the set $A$ and the three terms of integrals 
are associated to the three basic independent mechanisms.
\label{defnet}
\end{defi}
From $(\ref{Xi})$, it is straightforward that $N_{\tau-}=\langle \mathscr{N}_{\tau-},1\rangle$ is the size of the system before the jump time.  Although the SDE $(\ref{netdef})$ looks somewhat complicated, the principle is easy to interpret. The indicator functions that involve $\rho$ are related to the rates and appear to make use of the acceptance/rejection sampling technique described previously in Section $\ref{sec2}$. 
It is interesting to remark that in the second term (associated with the withdrawal component) the integral and the indicator function that involve $\rho$ may be canceled because the rate here does not depend on the space variable $x$. 

Now, we state in Theorem \ref{theo1} that if the stochastic process $\mathscr{N}$ solves $(\ref{netdef})$, then it follows the dynamic described by the infinitesimal generator $\mathcal{L}$ given  by $(\ref{generat})$.

\begin{theo}
Assume that conditions $(\ref{Ass1})$ and $(\ref{Ass2})$ are satisfied. Consider a $(\mathcal{F}_\tau)_{\tau\geq0}$-adapted stochastic process 
$(\mathscr{N}_\tau)_{\tau\geq0}$ 
that belongs a.s. to $\mathbb{T}_G$ and solves $(\ref{netdef})$. Then  $(\mathscr{N}_\tau)_{\tau\geq0}$ is Markovian  and its infinitesimal 
generator $\mathcal{L}$ is defined in particular for the class of test functions $\Phi$ (given by (\ref{Testfn})) by formula $(\ref{generat})$. 
\label{theo1} 
\end{theo} 
\proof The proof for the fact that $(\mathscr{N}_\tau)_{\tau\geq0}$ is  Markovian is classical and it will be omit to save place. 
Using the formula (\ref{jump}) and Definition \ref{defnet}, for any $\Phi$ given by (\ref{Testfn}) and all $\tau\geq0$, $\Phi(\mathscr{N}_\tau)$ is given a.s. by
\begin{align}
\begin{split}
 \Phi(\mathscr{N}_\tau)=\Phi(\mathscr{N}_0)&+\int_{0}^\tau\int_{\N^*}\int_{\mathbb{R}^d}\int_0^1  \mathds{1}_{\{i\leq N_{t-}\}}
 \mathds{1}_{\{\rho\leq ( k(x_{t-}^i,z))/(\gamma_1\tilde{k}(z)) \}} \\
&\qquad\qquad \times[\Phi(\mathscr{N}_{t-}+\{\delta_{(x_{t-}^i+z)}\})-\Phi(\mathscr{N}_{t-})]  P^{\alpha}(dt,di,dz,d\rho) \\
&+\int_{0}^\tau\int_{\N^{\ast}}\int_0^1 \mathds{1}_{\{i\leq N_{t-}\}}
 \mathds{1}_{\{\rho\leq(\beta/\beta)\}}  \\
&\qquad\qquad \times[\Phi(\mathscr{N}_{t-}-\{\delta_{(x_{t-}^i)}\})-\Phi(\mathscr{N}_{t-})]  P^{\beta}(dt,di,d\rho) \\
&+\int_{0}^\tau \int_{\N^{\ast}}\int_{\bar{\mathscr{D}}_G}\int_0^1 \mathds{1}_{\{i\leq N_{t-}\}}
 \mathds{1}_{\{\rho\leq (\mathsf{aff}(x_{t-}^i,y)k^{\mathsf{af}}(y))/(A_\mathsf{f}\gamma_2 \tilde{k}^{\mathsf{af}}(y) )\}} \\
&\qquad\qquad \times[\Phi(\mathscr{N}_{t-}+\{\delta_{(y)}\})-\Phi(\mathscr{N}_{t-})]
P^{\mathsf{af}}(dt,di,dy,d\rho).
\end{split}
\label{Phifst}
\end{align}
  
At first sight, let us  consider the Definition \ref{defpois} and taking expectation of (\ref{Phifst}) which give
\begin{align}
\nonumber
 \mathbb{E}[\Phi(\mathscr{N}_\tau)]&=\mathbb{E}[\Phi(\mathscr{N}_0)]+\int_0^\tau \mathbb{E}\Big[\alpha\gamma_1\tilde{k}(z)
 \sum_{i=1}^{N_{t-}}\frac{1}{\gamma_1\tilde{k}(z)}\\
\nonumber
&\qquad\qquad\quad \times \int_{\mathbb{R}^d}\Big\{\Phi(\mathscr{N}_{t-}+\{\delta_{(x_{t-}^i+z)}\})-
\Phi(\mathscr{N}_{t-})\Big\}k(x_{t-}^i,z)dz             \Big]dt\\
\nonumber
&\quad+\int_0^\tau \E\Big[\beta
 \sum_{i=1}^{N_{t-}}\frac{\beta}{\beta}
\Big\{\Phi(\mathscr{N}_{t-}-\{\delta_{(x_{t-}^i)}\})-
\Phi(\mathscr{N}_{t-})\Big\}  \Big]dt\\
\nonumber
&\quad+\int_0^\tau \E\Big[A_\mathsf{f}\gamma_2 \tilde{k}^{\mathsf{af}}(y) \sum_{i=1}^{N_{t-}}
\frac{1}{ A_\mathsf{f}\gamma_2 \tilde{k}^{\mathsf{af}}(y) } \\
\nonumber
&\qquad\qquad\quad \times \int_{\bar{\mathscr{D}}_G} \Big\{\Phi(\mathscr{N}_{t-}+\{\delta_{(y)}\})-\Phi(\mathscr{N}_{t-})\Big\}\mathsf{aff}(x_{t-}^i,y)
{k}^{\mathsf{af}}(y)dy\Big]dt\\
\begin{split}
&=\E[\Phi(\mathscr{N}_0)]+\beta\int_0^\tau \E\Big[ \int_{\bar{\mathscr{D}}_G}\mathscr{N}_{t}(dx) \Big\{\Phi(\mathscr{N}_{t}-\{\delta_{(x)}\})-
\Phi(\mathscr{N}_{t})\Big\}   \Big] dt \\
& \qquad+\alpha\int_0^\tau \E\Big[\int_{\bar{\mathscr{D}}_G} \mathscr{N}_{t}(dx) \int_{\mathbb{R}^d}\Big\{\Phi(\mathscr{N}_{t}+\{\delta_{(x+z)}\})-
\Phi(\mathscr{N}_{t})\Big\}k(x,z)dz \Big]dt\\
&\qquad+\int_0^\tau \E\Big[\int_{\bar{\mathscr{D}}_G}\Big\{ \int_{\bar{\mathscr{D}}_G}\big( \Phi(\mathscr{N}_{t}+\delta_{(y)})-
\Phi(\mathscr{N}_{t})\big) \mathsf{aff}(x,y) {k}^{\mathsf{af}}(y)dy \Big\}\mathscr{N}_{t}(dx)
\Big]dt.
\end{split}
\label{espgenet}
\end{align}
Now, differentiating the expression (\ref{espgenet}) at $\tau=0$ and using $\mathcal{L}\Phi(\mathscr{N}_0)=\partial_\tau\E[\Phi(\mathscr{N}_\tau)]_{\tau=0}$ lead immediately to (\ref{generat}). This completes the proof. \carre

In the sequel, the following remark will be important for further analysis. More importantly,  Remark \ref{theo2} ensures that the size of the vertices will grow with bound which guarantees that the system does not explode.
\begin{rem}
Admit assumptions (\ref{Ass1}) and (\ref{Ass2}). Consider that $\mathbb{E}[\langle \mathscr{N}_0,1\rangle ^m]<\infty$ for some $m\geq 1$, 
then for any  $0<T<\infty$:
 \begin{align}
  \mathbb{E}\Big[\sup_{\tau\in[0,T]}\langle \mathscr{N}_\tau,1\rangle^m\Big]<\infty.
  \label{noexplo}
 \end{align}
\label{theo2}
\end{rem}

In the following section, we shall study the infinite particle limit in order to establish the weak convergence of the measure-valued process.  

\section{Weak convergence}
\label{sec5}

We now turn to study the large graph limit by starting from a finite graph and taking the limit when the vertices size of the graph tends to infinity. Remark that the graph in the limit is infinite. Formally, we will prove that the process $(\mathscr{N}_\tau)_{\tau\geq0}$ converges in distribution towards a deterministic process in the space $\mathbb{D}([0,T],\mathcal{S}_{F}(\bar{\mathscr{D}}_G))$ equipped 
with the Skorokhod metric. We refer the reader to \citet[p.117]{Eth+Kur86} for a complete definition of the Skorokhod metric. Clearly, the vertices size $N_\tau\in\mathbb{N}$ is a random variable here and cannot be used to observe the asymptotic regime of the measure-valued process. Instead, one may fix $N_0=n\in \mathbb{N}^*$ and study the weak convergence when the initial size of vertices goes to infinity ($n\to \infty$). Globally, the weak convergence that we shall study typically goes as follows: one fixes $\theta\in (\mathbb{R}_+^*)^4 \times \mathbb{R}^{ |\theta^\alpha|+|\theta^{\mathsf{af}}|}$ and some sequence of graphs $(G_n)_{n\geq 1}$ (usually converging or increasing, in some sense, to an infinite graph, or belonging to some class of random graphs), and then studies the asymptotic behavior of the random point measure $\mathscr{N}_\tau(dx)$ given by (\ref{net}).

First, we gather a brief treatment  of some classical properties useful for weak convergence to keep forthcoming asymptotic analysis more clear by following the discussion in \cite{Bil99} to which 
we also refer for further details. We endow  the space of finite measures $\mathcal{S}_{F}(\bar{\mathscr{D}}_G)$ on $\bar{\mathscr{D}}_G$  with the topology of the weak convergence of measures, that is the smallest topology for which the applications $\xi \mapsto\langle\xi,f\rangle=\int_{\bar{\mathscr{D}}_G} f(x)\xi(dx)$ are  continuous for any $f\in C(\bar{\mathscr{D}}_G)$ (remark that $f$ is continuous and bounded where bounded since defined on a compact set). This topology is metrized by the Prokhorov metric  defined in \citet{Prok56} as the analogue of the L\'evy metric for more general spaces than $\mathbb{R}$ and is given, by
\begin{align}
 \pi_{\textrm{P}}(\xi,\xi')= \inf \Big\{\epsilon\geq 0:   \xi(E_0)\leq \xi'(E_0^\epsilon)+
 \epsilon \textrm{ and } \xi'(E_0)\leq \xi(E_0^\epsilon)+\epsilon\mbox{ for all closed $E_0 \subseteq \bar{\mathscr{D}}_G$} \Big \},
 \label{pro.dist}
\end{align}
where the $\epsilon$-inflation of a set is given by $E_0^\epsilon = \{x\in \bar{\mathscr{D}}_G; \inf_{y\in E_0}\|x-y\|<\epsilon \}$. This Prokhorov metric is bounded by the total variation metric \citep[p.34]{Hub81} 
denoted by  $\pi_{\textrm{TV}}(\xi,\xi')=\|\xi-\xi'\|_{\textrm{TV}}$ and associated with the norm, for any finite and signed measure $\xi$,
\begin{align}
 \|\xi\|_{\textrm{TV}}=\sup_{A\in\mathcal{B}(\bar{\mathscr{D}}_G)}|\xi(A)+\xi(A^c)|= \xi^+(\bar{\mathscr{D}}_G) + \xi^-(\bar{\mathscr{D}}_G)=
 \sup_{f \in C(\bar{\mathscr{D}}_G),\, \|f \|_\infty \leq 1} |\langle\xi,f\rangle|,
 \label{TVNorm}
\end{align}
where we recall that $\mathcal{B}(\bar{\mathscr{D}}_G)$ is the Borel $\sigma$-algebra on $\bar{\mathscr{D}}_G$ and $\xi=\xi^+-\xi^-$ is the Jordan-Hahn decomposition of $\xi$. 

Let also briefly give a  characterization of the  convergence for the Skorohod metric: A sequence $(\xi^n)_{n\in\mathbb{N}}$ converges to $\xi$ in 
 $\mathbb{D}([0,T],\mathcal{S}_{F}(\bar{\mathscr{D}}_G))$ (which means  $\pi_{\textrm{P}}(\xi^n,\xi)\to 0$) if and only if there exists a 
 sequence $\lambda_n(\tau)$ of time  change functions (i.e. strictly increasing bijective functions on $[0,T]$, with $\lambda_n(0)=0$ and 
 $\lambda_n(T)=T$) satisfying: 
\begin{align}
 \sup_{\tau\in[0,T]} \pi_{\textrm{P}}(\xi_\tau^n,\xi_{\lambda_n(\tau)}) {\stackrel{\hbox{{\scriptsize $n\to\infty$}}}{\longrightarrow}} 0 
 \qquad \textrm{ and } \qquad
 \sup_{\tau\in[0,T]} |\lambda_n(\tau)-\tau| \to 0.
 \label{conv.skoro}
\end{align}
When the sequence $(\xi^n)_{n\in\mathbb{N}}$ converges to $\xi$ in $\mathbb{D}([0,T],\mathcal{S}_{F}(\bar{\mathscr{D}}_G))$ and if 
$\xi\in C([0,T],\mathcal{S}_{F}(\bar{\mathscr{D}}_G))$, then 
\begin{align}
 \sup_{\tau\in[0,T]} \pi_{\textrm{P}}(\xi_\tau^n,\xi_\tau) \leq \sup_{\tau\in[0,T]} \pi_{\textrm{P}}(\xi_\tau^n,\xi_{\lambda_n(\tau)}) + 
 \sup_{\tau\in[0,T]} \pi_{\textrm{P}}(\xi_{\lambda_n(\tau)},\xi_\tau) \to 0
 \label{conv.unif}
\end{align}
because of (\ref{conv.skoro}) and the uniform continuity of $\xi$ in $[0,T]$. Therefore $\xi^n$ converges to $\xi$ in 
$\mathbb{D}([0,T],\mathcal{S}_{F}(\bar{\mathscr{D}}_G))$ and at the same time for the uniform metric. In other words, the Prokhorov  metric coincides with the uniform metric on ${C}([0,T],\mathcal{S}_{F}(\bar{\mathscr{D}}_G))$, the space of $\mathcal{S}_F$-valued continuous functions on $[0,T]$ (see \cite{Eth+Kur86}, Chapter 3 for more details).

Let start the study of the convergence in law of  $(\mathscr{N}_\tau)_{\tau\in[0,T]}$ on the space $\mathbb{D}([0,T],\mathcal{S}_{F}(\bar{\mathscr{D}}_G))$ of  c\`adl\`ag functions on $[0,T]$ with values in $\mathcal{S}_F(\bar{\mathscr{D}}_G)$ by establishing some martingale properties.

\subsection{Law of large numbers scaling}

As previously emphasized, we give now martingale properties that result from standard stochastic calculus for jump processes and SDE driven by PPMs (see, e.g., \cite{Jaco+Shir87}). 
\begin{prop}
Assume that $(\ref{Ass1})$ and $(\ref{Ass2})$ are satisfied. Also assume that $\mathbb{E}[\langle \mathscr{N}_0,1\rangle ^m]<\infty$ for some $m\geq 1$. Consider the process $(\mathscr{N}_\tau)_{(\tau\geq 0)}$ given by Definition $\ref{defnet}$ and consider $\mathcal{L}$ its infinitesimal generator defined by 
$(\ref{generat})$. Then,
\begin{enumerate}
 \item For all functions $\Phi$ given by (\ref{Testfn}), with  $F\in C^1(\R)$ and $f\in C(\bar{\mathscr{D}}_G)$, such that for some constant $C$ and for all $\mathscr{N}\in\mathcal{S}_G$, 
$|\Phi(\mathscr{N})|+|\mathcal{L}\Phi(\mathscr{N})|\leq C(1+\langle \mathscr{N},1\rangle ^m),$ the process
\begin{align}
\Phi(\mathscr{N}_\tau)-\Phi(\mathscr{N}_0)-\int_0^\tau  \mathcal{L}\Phi(\mathscr{N}_t) dt
\label{marting}
\end{align}
is a c\`adl\`ag $L^1-(\mathcal{F}_\tau)_{\tau\geq0}$-martingale starting from $0$.
\item  For any  function $f\in C(\bar{\mathscr{D}}_G)$, the process 
\begin{align}
\begin{split}
M_\tau^f &=\langle \mathscr{N}_{\tau},f\rangle-\langle \mathscr{N}_{0},f\rangle-\alpha \int_0^{\tau} dt \int_{\bar{\mathscr{D}}_G}\mathscr{N}_t(dx) \int_{\mathbb{R}^d}f(x+z)k(x,z)dz \\ 
& +\beta \int_0^{\tau}dt  \int_{\bar{\mathscr{D}}_G} f(x) \mathscr{N}_t(dx)  -\int_0^{\tau}dt\int_{\bar{\mathscr{D}}_G}\Big\{ \int_{\bar{\mathscr{D}}_G}f(y) \mathsf{aff}(x,y) {k}^{\mathsf{af}}(y)dy \Big\}\mathscr{N}_t(dx)
\end{split}
\end{align}
is a c\`adl\`ag $L^2$-martingale starting from $0$ with quadratic variation:

\begin{align}
\begin{split}
[ M^f]_\tau&=\alpha \int_0^{\tau} dt \int_{\bar{\mathscr{D}}_G}\mathscr{N}_t(dx)\int_{\mathbb{R}^d}f^2(x+z)k(x,z)dz
+\beta \int_0^{\tau}dt \int_{\bar{\mathscr{D}}_G} f^2(x) \mathscr{N}_t(dx)\\
& \qquad + \int_0^{\tau}dt\int_{\bar{\mathscr{D}}_G}\Big\{ \int_{\bar{\mathscr{D}}_G}f^2(y) \mathsf{aff}(x,y) {k}^{\mathsf{af}}(y)dy \Big\}\mathscr{N}_t(dx).
\end{split}
\label{brack}
\end{align}
\end{enumerate}
\label{Propmarting}
\end{prop}

\proof The point $(i)$ is trivial.  In fact, $\forall \tau\geq 0$, $\mathscr{N}_\tau$  is a Markov Process. Therefore, the Dynkin's formula and 
$(\ref{noexplo})$ give immediately that $(\ref{marting})$ is a c\`adl\`ag $L^1-(\mathcal{F}_{\tau})_{\tau\geq 0}$-martingale. To prove the point $(ii)$, 
assume that $\mathbb{E}[\langle \mathscr{N}_{0},1\rangle^3]<\infty$. Then, applying the point $(i)$ with $\Phi(\mathscr{N})=\langle \mathscr{N},f\rangle$ yields that $M^f$ is a martingale. To compute its brackets, we first apply $(i)$ with $\Phi(\mathscr{N})=\langle \mathscr{N},f\rangle^2$ and deduce that:
\begin{align}
\begin{split}
\langle \mathscr{N}_{\tau},f\rangle^2-&\langle \mathscr{N}_{0},f\rangle^2-\int_0^\tau  \mathcal{L}\langle \mathscr{N}_t,f\rangle^2 dt=
\langle \mathscr{N}_{\tau},f\rangle^2-\langle \mathscr{N}_{0},f\rangle^2 \\
&-\alpha \int_0^{\tau} dt\int_{\bar{\mathscr{D}}_G}\mathscr{N}_t(dx) \int_{\mathbb{R}^d}\Big\{2\,f(x+z)\langle \mathscr{N}_t,f\rangle+f^2(x+z)\Big\}\,k(x,z)dz \\ 
&+\beta \int_0^{\tau}dt\, \int_{\bar{\mathscr{D}}_G}(f^2(x)-2f(x)\langle \mathscr{N}_t,f\rangle )\,\mathscr{N}_t(dx)   \\
&-\int_0^{\tau}dt\int_{\bar{\mathscr{D}}_G}\Big\{ \int_{\bar{\mathscr{D}}_G}\Big\{2\,f(y)\langle \mathscr{N}_t,f\rangle+f^2(y)\Big\} \mathsf{aff}(x,y) {k}^{\mathsf{af}}(y)dy 
\Big\}\mathscr{N}_t(dx)
\end{split}
\label{brack1}
\end{align}
is a martingale. Hence we apply the It\^o's formula  in order to compute $\langle \mathscr{N}_{\tau},f\rangle^2$ from $M_\tau^f$ to find that 
\begin{align}
\begin{split}
\langle \mathscr{N}_{\tau},f\rangle^2-\langle \mathscr{N}_{0},f\rangle^2&-2\alpha  \int_0^{\tau}dt \int_{\bar{\mathscr{D}}_G}\mathscr{N}_t(dx) \int_{\mathbb{R}^d}f(x+z)\langle \mathscr{N}_{t},f\rangle\,k(x,z)\,dz \\
&+2 \beta  \int_0^{\tau}dt\, \int_{\bar{\mathscr{D}}_G}f(x)\langle \mathscr{N}_{t},f\rangle \mathscr{N}_t(dx)\\
&-2\int_0^{\tau}dt\,\int_{\bar{\mathscr{D}}_G}\Big\{ \int_{\bar{\mathscr{D}}_G}f(y)\langle \mathscr{N}_{t},f\rangle\,\mathsf{aff}(x,y) {k}^{\mathsf{af}}(y)dy \Big\}\mathscr{N}_t(dx)-\langle M^f\rangle_\tau
\end{split}
\label{brack2}
\end{align}
is a martingale using the Doob-Meyer martingale representation theorem. Finally, comparing $(\ref{brack1})$ and $(\ref{brack2})$, leads to $(\ref{brack})$. \carre

For each $n\in\mathbb{N}^{\star}$, we put $n$ the initial size of the system, i.e. $\langle \mathscr{N}^n_0,1 \rangle=n$, and define a set of parameters and kernels  
$(\alpha^n,K^n,w^{\mathsf{af}}_n,\mathsf{aff}^n, K^{\mathsf{af}}_n, \beta^n,a_\mathsf{f}^n,A_\mathsf{f}^n)$ as in section $(\ref{sec2.2})$ satisfying $(\ref{Ass1})$ and $(\ref{Ass2})$. As well, we consider sequence of measures $(\mathscr{N}^n)_{n\in \mathbb{N}}$ such that, for any $n\in \mathbb{N}^*$, $\mathscr{N}^n$ satisfies SDE (\ref{netdef}) stated in Definition \ref{defnet}
with initial condition $\mathscr{N}_0^n$.  Consider the subset $\mathcal{S}_G^n$  of $\mathcal{S}_F(\bar{\mathscr{D}}_G)$ embedded with the weak convergence topology and defined by,
\begin{align*}
\mathcal{S}_G^n=\Big\{\frac{1}{n}\mathscr{N}, \; \mathscr{N}\in\mathcal{S}_G\Big\}.
\end{align*}
Let now scale the measures in the following way, for all $n\in\mathbb{N}^*$
\begin{align}
\bar{\mathscr{N}}_\tau^n=\frac{1}{n} \mathscr{N}_\tau^n, \qquad \tau\geq 0.
\label{normlr}
\end{align}
The measures $(\bar{\mathscr{N}}_\tau^n)$ represent the mean-field approximation when the initial size $n$ grows to infinity.  From now 
on, we assume that:
\begin{align}
 \bar{\mathscr{N}}_0^n =  \frac{1}{n}\mathscr{N}_0^n {\stackrel{\hbox{{\scriptsize $n\to\infty$}}}{\longrightarrow}} \xi_0  \textrm{ in distribution in } \mathcal{S}_F(\bar{\mathscr{D}}_G),
\label{asymnorm}
 \end{align}
where $\xi_0$ is the limit measure after renormalization at the initial time. We suppose that $\xi_0$ is deterministic together with 
$\langle \xi_0,1\rangle>0$. We obtain rescaled SDE which is the same as the SDE (\ref{netdef}) parameterized by $n$. Hence, the fact that $\bar{\mathscr{N}}^n$ is a Markov process is straightforward given that $\mathscr{N}^n$ is Markovian. Moreover, we may state that its infinitesimal generator can be easily deduced from $(\ref{generat})$.

\begin{prop}
Let $n\in\mathbb{N}^{\star}$. Denote $\mathcal{L}^n:\mathcal{S}_G^n\mapsto\mathbb{R}$ the infinitesimal generator of $(\bar{\mathscr{N}}_\tau^n)_{\tau\geq 0}$.  Then, it is defined for all functions $\Phi$ given by $(\ref{Testfn})$ and for all $\bar{\mathscr{N}}^n\in\mathcal{S}_G^n$ by 
\begin{align}
\begin{split}
\mathcal{L}^n\Phi(\bar{\mathscr{N}}^n)&=n\,\alpha^n\int_{\bar{\mathscr{D}}_G}\bar{\mathscr{N}}^n(dx) \int_{\mathbb{R}^d}\Big\{\Phi(\bar{\mathscr{N}}^n+\frac{\delta_{x+z}}{n})-\Phi(\bar{\mathscr{N}}^n)\Big\}k^n(x,z)dz \\
          &\qquad +n\,\beta^n \int_{\bar{\mathscr{D}}_G} \Big\{\Phi(\bar{\mathscr{N}}^n-\frac{\delta_{x}}{n})-\Phi(\bar{\mathscr{N}}^n)\Big\} \bar{\mathscr{N}}^n(dx)   \\
          &\qquad+n\int_{\bar{\mathscr{D}}_G}\Big\{ \int_{\bar{\mathscr{D}}_G}\big( \Phi(\bar{\mathscr{N}}^n+\frac{\delta_{y}}{n})-\Phi(\bar{\mathscr{N}}^n)\big) \mathsf{aff}^n(x,y)
          {k}_n^{\mathsf{af}}(y)dy \Big\}\bar{\mathscr{N}}^n(dx).
\end{split}
\label{genenor}
\end{align}
\label{theogenenor}
\end{prop}

\proof Let $n\in\mathbb{N}^{\star}$ and  $\Phi$ a function defined by $(\ref{Testfn})$. Furthermore, define the function $\Phi^n$ by 
$\Phi^n(\mathscr{N}_\tau^n)=\Phi(\frac{\mathscr{N}_\tau^n}{n})$. Let $ \bar{\mathscr{N}}^n\in\mathcal{S}_G^n$, and $\tilde{\mathcal{L}}$  the generator of $(\mathscr{N}_\tau^n)_{\tau\geq 0}$ given by 
$(\ref{generat})$. Note that $\tilde{\mathcal{L}}$ is simply obtained in a similar manner as $\mathcal{L}$. Then,
\begin{align*}
\mathcal{L}^n\Phi(\bar{\mathscr{N}}_\tau^n)&=\partial_\tau\mathbb{E}[\Phi(\bar{\mathscr{N}}_\tau^n)]_{\tau=0}=\partial_\tau\mathbb{E}[\Phi(\frac{\mathscr{N}_\tau^n}{n})]_{\tau=0}=
\partial_\tau\mathbb{E}[\Phi^n(\mathscr{N}_\tau^n)]_{\tau=0}\\
              &=\alpha^n\int_{\bar{\mathscr{D}}_G}\mathscr{N}^n(dx) \int_{\mathbb{R}^d}\Big\{\Phi^n(\mathscr{N}^n+\delta_{x+z})-\Phi^n(\mathscr{N}^n)\Big\}k^n(x,z)dz \\
          &\qquad +\beta^n \int_{\bar{\mathscr{D}}_G} \Big\{\Phi^n(\mathscr{N}^n-\delta_{x})-\Phi^n(\mathscr{N}^n)\Big\} \mathscr{N}^n(dx)   \\
          &\qquad+\int_{\bar{\mathscr{D}}_G}\Big\{ \int_{\bar{\mathscr{D}}_G}\big( \Phi^n(\mathscr{N}^n+\delta_{y})-\Phi^n(\mathscr{N}^n)\big) \mathsf{aff}^n(x,y) {k}_n^{\mathsf{af}}(y)dy \Big\}\mathscr{N}^n(dx).
\end{align*}
We complete the proof by replacing $\Phi^n$ by $\Phi$ and $\mathscr{N}^n$ by $n\,\bar{\mathscr{N}}^n$ to get $(\ref{genenor})$. \carre

We are  now in position to reformulate the martingale property for the renormalized measure $\bar{\mathscr{N}}^n$, for which the proof follows the same spirit as  Proposition \ref{Propmarting}. 
\begin{prop}
Admit assumptions (\ref{Ass1}) and (\ref{Ass2}) and assume that $\mathbb{E}[\langle \bar{\mathscr{N}}_0^n,1\rangle ^m]<\infty$ for some $m\geq 1$. Define the process $(\bar{\mathscr{N}}_\tau^n)_{\tau\geq 0}$ by (\ref{normlr}) with the infinitesimal generator  $\mathcal{L}^n$, then
\begin{enumerate}
\item  For all  functions $\Phi$ given by (\ref{Testfn}), with  $F\in C^1(\R)$ and $f\in C(\bar{\mathscr{D}}_G)$, such that for some constant $C$ and for all $\bar{\mathscr{N}}^n\in\mathcal{S}_G^n$, 
$|\Phi(\bar{\mathscr{N}}^n)|+|\mathcal{L}^n\Phi(\bar{\mathscr{N}}^n)|\leq C(1+\langle \bar{\mathscr{N}}^n,1\rangle ^m),$ the process
\[
\Phi(\bar{\mathscr{N}}_\tau^n)-\Phi(\bar{\mathscr{N}}_0^n)-\int_0^\tau dt\mathcal{L}^n\Phi(\bar{\mathscr{N}}_t^n)
\]
is a c\`adl\`ag $L^1-(\mathcal{F}_\tau)_{\tau\geq0}$-martingale starting from $0$.
\item  For any  function $f\in C(\bar{\mathscr{D}}_G)$, the process 
\begin{align}
\begin{split}
M_\tau^{n,f} =\langle \bar{\mathscr{N}}_\tau^n,f\rangle-\langle \bar{\mathscr{N}}_0^n,f\rangle&-\alpha^n\int_0^{\tau} dt \int_{\bar{\mathscr{D}}_G}\bar{\mathscr{N}}_t^n(dx) 
\int_{\mathbb{R}^d}f(x+z)k^n(x,z)dz \\ 
&+\beta^n\int_0^{\tau}dt  \int_{\bar{\mathscr{D}}_G} f(x) \bar{\mathscr{N}}_t^n(dx)   \\
&-\int_0^{\tau}dt\int_{\bar{\mathscr{D}}_G}\Big\{ \int_{\bar{\mathscr{D}}_G}f(y) \mathsf{aff}^n(x,y) {k}_n^{\mathsf{af}}(y)dy \Big\}\bar{\mathscr{N}}_t^n(dx)
\end{split}
\label{martingn}
\end{align}
is a c\`adl\`ag $L^2$-martingale starting from $0$ with quadratic variation:
\begin{align}
\begin{split}
[ M^{n,f}]_\tau=&\frac{\alpha^n}{n}\int_0^{\tau} dt \int_{\bar{\mathscr{D}}_G}\bar{\mathscr{N}}_t^n(dx)\int_{\mathbb{R}^d}f^2(x+z)k^n(x,z)dz\\
&+\frac{\beta^n}{n}\int_0^{\tau}dt  \int_{\bar{\mathscr{D}}_G} f^2(x) \bar{\mathscr{N}}_t^n(dx)\\
&+ \frac{1}{n}\int_0^{\tau}dt\int_{\bar{\mathscr{D}}_G}\Big\{ \int_{\bar{\mathscr{D}}_G}f^2(y) \mathsf{aff}^n(x,y) {k}_n^{\mathsf{af}}(y)dy \Big\}\bar{\mathscr{N}}_t^n(dx).
\end{split}
\label{it12}
\end{align}
\end{enumerate}
\label{Proposition}
\end{prop}

\proof The proof proceeds from (\ref{netdef}) and generator (\ref{genenor}) together with standard stochastic calculus for jump processes  as in the proof of  proposition \ref{Propmarting}.  \carre

\subsection{Large graph limit}

Now, we focus on the study of the system limit when $n\to \infty$. Particularly,  we establish  the convergence of the normalized measure and show that can be approximated by a deterministic equation that might be another  way  to describe the proposed model here. Again, we assume that the local affinity between two vertices of the rescaled process  is uniformly bounded in such a way that the Hypothesis $\ref{H}$ is as follows:  
\begin{hyp}
Let us assume that:
\begin{enumerate}
\item There exists a constant 
$A_\mathsf{f}^n$ such that, for all $x,y\in\bar{\mathscr{D}}_G$ and for all $\bar{\mathscr{N}}^n\in\mathcal{S}_G^n$,
\begin{align*}
\mathsf{aff}^n(x,y)\leq A_\mathsf{f}^n \qquad \textrm{ which implies } \qquad  w_n^{\mathsf{af}}(y,\bar{\mathscr{N}}^n)\leq A_\mathsf{f}^n 
\langle \bar{\mathscr{N}}^n,1 \rangle.
\end{align*}
As well, there exists a bounded nonnegative function $\mathsf{aff}$ on $\bar{\mathscr{D}}_G\times \bar{\mathscr{D}}_G $  such that $\mathsf{aff}^n(x,y)=\frac{\mathsf{aff}(x,y)}{n}$ and there exist some constants $\alpha,\beta, A_\mathsf{f},a_\mathsf{f}$ such that   $A_\mathsf{f}^n=A_\mathsf{f} $, $\alpha^n=\alpha$, $a_\mathsf{f}^n=a_\mathsf{f}$ and $\beta^n=\beta$.
\item  There exist some reals $\gamma_3,\gamma_4>0$ and  probability densities $\tilde{k}^n$ on $\mathbb{R}^d$ and $\tilde{k}_n^{\mathsf{af}}$ on $\bar{\mathscr{D}}_G$ such that, for all $x\in\bar{\mathscr{D}}_G$,
\begin{align*}
 k^n(x,z)\leq \gamma_3 \tilde{k}^n(z) \qquad \textrm{ and } \qquad  k_n^{\mathsf{af}}(y)\leq \gamma_4 \tilde{k}_n^{\mathsf{af}}(y).
\end{align*}
Further, there exist a continuous nonnegative functions $k$ on $\bar{\mathscr{D}}_G\times \mathbb{R}^d$ that satisfies (\ref{HypKern1}) and $k^{\mathsf{af}}$ on $\bar{\mathscr{D}}_G$ that satisfies (\ref{HypKern2}) such that  $k^n(x,z)=k(x,z)$ and $k_n^{\mathsf{af}}(y)=k^{\mathsf{af}}(y)$.  
 \item For $m\geq 1$, $\sup_n\mathbb{E}[\langle \bar{\mathscr{N}}_0^n,1\rangle^m]<\infty$.
\end{enumerate}
\label{H}
\end{hyp}
In the Hypothesis $\ref{H}$, the first assertion ($i$) is convenient since it avoids explosion phenomena. The second assertion ($ii$) allows to bound  the dispersion of new vertices after recruitments by invitation or affinity. The third assertion ($iii$) is easily satisfied through an appropriate choice of the initial distribution $\xi_0$.             

Our main result is the following theorem.
\begin{theo}[Weak convergence]
Assume that the three assertions of Hypothesis $\ref{H}$ are fulfilled, then the rescaled  process $(\bar{\mathscr{N}}_\tau^n)_{\tau\geq 0}$ converges in law  in the space 
$\mathbb{D}([0,T],\mathcal{S}_{F}(\bar{\mathscr{D}}_G))$ towards the unique solution $(\xi_\tau)_{\tau\geq 0}$, satisfying 
$\sup_{\tau\in[0,T]}\langle\xi_\tau,1\rangle
<\infty$, of the deterministic equation:
\begin{align}
\begin{split}
\langle\xi_\tau,f\rangle &=\langle\xi_0,f\rangle+\alpha \int_0^\tau dt\,\int_{\bar{\mathscr{D}}_G}\xi_t(dx) \int_{\mathbb{R}^d}f(x+z)k(x,z)dz -
\beta \int_0^\tau dt \int_{\bar{\mathscr{D}}_G} f(x)\xi_t(dx)   \\
          &\qquad +\int_0^\tau dt \int_{\bar{\mathscr{D}}_G}\Big\{ \int_{\bar{\mathscr{D}}_G} f(y)\mathsf{aff}(x,y) {k}^{\mathsf{af}}(y)dy \Big\}\xi_t(dx),
\end{split}
\label{integrodiff}
\end{align}
where (\ref{integrodiff}) is defined for all  functions $f\in C(\bar{\mathscr{D}}_G)$. 
\label{theo5}
\end{theo}

In some words, the last theorem states that under some assumptions, the sequence of measures $(\bar{\mathscr{N}}^n)_{n\in\mathbb{N}}$ tends in law to a deterministic continuous measure-valued function $(\xi_\tau)_{\tau\geq 0}$  which is the unique solution of the last integrodifferential equation (\ref{integrodiff}). In practice, it seems that the deterministic equation (\ref{integrodiff}) will be valid in large system size and can be taken as a limit approximation for the random model in large network sizes.

\begin{rem}
The previous results shown with the class of measurable and bounded functions $f$ from $\bar{\mathscr{D}}_G$ into $\mathbb{R}$ can be 
straightforwardly generalized to the class of functions that depend on time. In particular, the solution $(\xi_\tau)_{\tau\geq 0}$ of the equation 
$(\ref{integrodiff})$ is solution of
 \begin{align*}
\langle\xi_\tau,f_\tau\rangle &=\langle\xi_0,f_0\rangle+\alpha \int_0^\tau dt\,\int_{\bar{\mathscr{D}}_G}\xi_t(dx) \int_{\mathbb{R}^d}f_t(x+z)k(x,z)dz -
\beta \int_0^\tau dt \int_{\bar{\mathscr{D}}_G} f_t(x)\xi_t(dx)   \\
          &\qquad +\int_0^\tau dt \int_{\bar{\mathscr{D}}_G}\Big\{ \int_{\bar{\mathscr{D}}_G} f_t(y)\mathsf{aff}(x,y) {k}^{\mathsf{af}}(y)dy \Big\}\xi_t(dx),
\end{align*}
where we consider here all time-dependent functions $f\in C([0,\infty)\times \bar{\mathscr{D}}_G): (\tau,x)\mapsto f_\tau(x)$.
\end{rem}

We split the proof of Theorem \ref{theo5} into several technical lemmas. Before presenting the following statements, we fix $T>0$. The aim at first in the following lemma is to show the uniqueness of the solution for the  integrodifferential equation (\ref{integrodiff}). 
\begin{lem}
Under the assumptions of Hypothesis $\ref{H}$, the solution of the deterministic integrodifferential equation (\ref{integrodiff}) is unique.
\label{unik}
\end{lem}
\proof We are  interested in proving  the uniqueness of the solution. For this,  let $(\xi_\tau)_{\tau\geq 0}$ and $(\xi_\tau')_{\tau\geq 0}$ two 
solutions of $(\ref{integrodiff})$ satisfying
\[
\sup_{\tau\in[0,T]}\langle\xi_{\tau}+\xi_\tau',1\rangle=B_T<+\infty.
\] 
There may be some doubt as to why $B_T$ is finite. Let us verify this for one solution, say $(\xi_\tau)_{\tau\geq 0}$,  we have
\begin{align*}
\langle\xi_\tau,1\rangle &=\langle\xi_0,1\rangle+\alpha \int_0^\tau dt\,\int_{\bar{\mathscr{D}}_G}\xi_t(dx) \int_{\mathbb{R}^d}k(x,z)dz -
\beta\int_0^\tau dt \int_{\bar{\mathscr{D}}_G} \xi_t(dx)   \\
          &\qquad +\int_0^\tau dt \int_{\bar{\mathscr{D}}_G}\Big\{ \int_{\bar{\mathscr{D}}_G} \mathsf{aff}(x,y) {k}^{\mathsf{af}}(y)dy \Big\}\xi_t(dx)\\
&\leq \langle\xi_0,1\rangle + \alpha \gamma_3 \int_0^\tau \langle\xi_t,1\rangle dt + \beta\int_0^\tau \langle\xi_t,1\rangle dt + 
\gamma_4 A_\mathsf{f} \int_0^\tau \langle\xi_t,1\rangle dt  \\
&\leq  \langle\xi_0,1\rangle  \exp\{(\alpha \gamma_3 + \beta + \gamma_4 A_\mathsf{f})\tau\}<+\infty
\end{align*}
where the last line is established by the Gronwall's lemma. We consider now the variation norm $\|\nu_1-\nu_2\|_{\textrm{TV}}$ defined 
in $(\ref{TVNorm})$ for all $\nu_1,\nu_2\in\mathcal{S}_F(\bar{\mathscr{D}}_G)$. Also, we consider some function $f\in C(\bar{\mathscr{D}}_G)$ such that 
$\|f\|_\infty\leq 1$. Therefore, by a straightforward calculation we obtain
\begin{align*}
\left|\langle\xi_\tau-\xi_\tau',f\rangle\right|&=\Big| \langle\xi_0-\xi_0',f\rangle  +\alpha \int_0^\tau dt\int_{\bar{\mathscr{D}}_G}(\xi_t(dx)-\xi_t'(dx)) \int_{\mathbb{R}^d}f(x+z)k(x,z)dz\\
 &\qquad\qquad\qquad \; - \beta \int_0^\tau dt \int_{\bar{\mathscr{D}}_G} f(x)(\xi_t(dx)-\xi_t'(dx))\\
 &\qquad\qquad\qquad \;  + \int_0^\tau  dt \int_{\bar{\mathscr{D}}_G}\Big\{ \int_{\bar{\mathscr{D}}_G} f(x)\mathsf{aff}(x,y)  {k}^{\mathsf{af}}(y)dy \Big\}(\xi_t(dx)-\xi_t'(dx))    \Big|\\        
 &\leq  \int_0^\tau dt\Big|\int_{\bar{\mathscr{D}}_G}(\xi_t(dx)-\xi_t'(dx)) \Big\{ \alpha\int_{\mathbb{R}^d}f(x+z)k(x,z)dz - \beta f(x)\Big\} \Big|\\
 &\qquad\qquad\qquad \; + \int_0^\tau dt \Big|\int_{\bar{\mathscr{D}}_G}\Big\{ \int_{\bar{\mathscr{D}}_G} f(x) \mathsf{aff}(x,y) {k}^{\mathsf{af}}(y)dy \Big\}(\xi_t(dx)-\xi_t'(dx)) \Big|.                                           \end{align*}
Now, using the assertions of Hypothesis $\ref{H}$ together with $\|f\|_\infty\leq 1$, we can establish easily the following bounds:
\[
\Big|\alpha \int_{\mathbb{R}^d}f(x+z)k(x,z)dz-\beta f(x)\Big|\leq \alpha \gamma_3 +\beta
\]
and
\[
\Big|\int_{\bar{\mathscr{D}}_G} f(x)\mathsf{aff}(x,y) {k}^{\mathsf{af}}(y)dy\Big| \leq {A}_{\mathsf{f}}\gamma_4
\]
which immediately yields that

\begin{align*}
\left|\langle\xi_\tau-\xi_\tau',f\rangle\right|& \leq (\alpha \gamma_3 +\beta+{A}_{\mathsf{f}}\gamma_4)\int_0^\tau \|\xi_t-\xi_t'\|_{\textrm{TV}}dt.
\end{align*}

Finally, we take the supremum over all functions $f\in C(\bar{\mathscr{D}}_G)$ such that $\|f\|_{\infty}\leq 1$ to establish
\begin{align*}
\sup_{f\in C(\bar{\mathscr{D}}_G),\|f\|_\infty\leq 1}\left|\langle\xi_\tau-\xi_\tau',f\rangle\right|= \|\xi_\tau-\xi_\tau'\|_{\textrm{TV}} 
&\leq  (\alpha \gamma_3 +\beta+{A}_{\mathsf{f}}\gamma_4)\int_0^\tau \|\xi_t-\xi_t'\|_{\textrm{TV}} dt\\
&\leq 0\times \exp\{ (\alpha \gamma_3 +\beta+{A}_{\mathsf{f}}\gamma_4) \tau \},
\end{align*}
where the last line is obtained from the Gronwall's inequality.  We conclude that, for all $\tau\geq 0$, $\xi_\tau = \xi_\tau'$. Hence,  (\ref{integrodiff}) has a unique solution. \carre

Now, let us prove some moment estimates. 
\begin{lem}
Under the assertions of the Hypothesis $\ref{H}$ and for some $m\geq 1$, we have 
\begin{align}
\sup_{n \in \mathbb{N}^*}\mathbb{E}\Big[\sup_{\tau\in[0,T]}\langle \bar{\mathscr{N}}_\tau^n,1\rangle^m\Big]< \infty.
\label{boundnorm}
\end{align}
\label{lemboundnorm}
\end{lem}
\proof Let $n\in\mathbb{N}^{\star}$ and consider the process $(\mathscr{N}_\tau^n)_{\tau\geq 0}$. Then, by Remark \ref{theo2}, 
we find that
\[
\mathbb{E}\Big[\sup_{\tau\in[0,T]}\langle \mathscr{N}_\tau^n,1\rangle^m\Big]<C_{m,T}\mathbb{E}[\langle \mathscr{N}_0^n,1\rangle^m].
\]
We recall that the constant $C_{m,T}$ is not dependent on $n$ and that the rescaled process is defined by $\bar{\mathscr{N}}_\tau^n=\frac{1}{n}\mathscr{N}_\tau^n$. Hence, by 
assertion $(iii)$ of Hypothesis $\ref{H}$ we conclude that,   
\begin{align*}
\sup_n\mathbb{E}\Big[\sup_{\tau\in[0,T]} \langle \frac{\mathscr{N}^n_\tau}{n},1\rangle^m\Big]\leq \sup_n\mathbb{E}[\langle \frac{\mathscr{N}_0^n}{n},1\rangle ^m] C_{m,T}<\infty,
\end{align*}
which completes the proof. \carre

We shall now prove the tightness of the sequence $(\bar{\mathscr{N}}^n)_{n\in \mathbb{N}^*}$ in $\mathbb{D}([0,T],\mathcal{S}_{F}(\bar{\mathscr{D}}_G))$.
We first endow $\mathcal{S}_F(\bar{\mathscr{D}}_G)$ with the weak topology  and $\mathbb{D}([0,T],\mathcal{S}_{F}(\bar{\mathscr{D}}_G))$ with the Skorokhod  topology. Note that the space $\mathbb{D}([0,T],\mathcal{S}_{F}(\bar{\mathscr{D}}_G))$ equipped with the Skorokhod topology is complete  \citep[Theorem 12.2]{Bil68}. Therefore, 
according to Prokhorov's theorem  \citep{Prok56}, the sequence $(\bar{\mathscr{N}}^n)_{n\in \mathbb{N}^*}$ is tight in 
$\mathbb{D}([0,T],\mathcal{S}_{F}(\bar{\mathscr{D}}_G))$ if and only if it is relatively compact (i.e. its closure is compact). Hence,  the tightness of 
$(\bar{\mathscr{N}}^n)_{n\in \mathbb{N}^*}$ is equivalent to the fact that from any subsequence one can extract a subsequence that converges in distribution in the space $\mathbb{D}([0,T],\mathcal{S}_{F}(\bar{\mathscr{D}}_G))$.

\begin{lem}
Endow the space of finite measures $\mathcal{S}_F(\bar{\mathscr{D}}_G)$ on $\bar{\mathscr{D}}_G$ with the topology of the weak convergence of measures (metrized by the Prokhorov metric), and suppose that the assertions of Hypothesis $\ref{H}$ are fulfilled. Then, the sequence of probability laws $\mathscr{L}^n=\mathscr{L}(\bar{\mathscr{N}}^n)$ of the family $(\bar{\mathscr{N}}^n)_{n\in\mathbb{N}^*}$ is tight in 
$\mathscr{P}(\mathbb{D}([0,T],\mathcal{S}_{F}(\bar{\mathscr{D}}_G))$, the space of probability measures on $\mathbb{D}([0,T],\mathcal{S}_{F}(\bar{\mathscr{D}}_G))$.  
 \label{tightness}
 \end{lem}
\proof First, note that according to \cite[Theorem 2.1]{Copp+86}, proving the tightness of the sequence $(\bar{\mathscr{N}}^n)_{n\in \mathbb{N}^*}$  is amounts to proving the tightness of its projections $\langle\bar{\mathscr{N}}^n,f\rangle$ in  $\mathbb{D}([0,T],\mathbb{R})$ for all functions $f$ in a dense sequence of  $C_b(\bar{\mathscr{D}}_G)$, the set of bounded  continuous functions on $\bar{\mathscr{D}}_G$. Moreover, to show the latter result, it is sufficient, 
following \cite{Jof+Met86},  to verify the following Aldous-Rebolledo criteria  \citep{Aldous78}:
\begin{itemize}
 \item [(A)]: For any $0<\tau\leq T$, the sequence $(\langle\bar{\mathscr{N}}_\tau^n,f\rangle)_{n\in\mathbb{N}^*}$ is tight in $\mathbb{R}$ for any bounded continuous 
 function $f$ on $\bar{\mathscr{D}}_G$.
 \item [(A')]: Consider from (\ref{martingn}) the following semimartingale decomposition,
 \begin{align*}
  \langle\bar{\mathscr{N}}_\tau^n,f\rangle = \langle\bar{\mathscr{N}}_0^n,f\rangle +A_\tau^n +M_\tau^{n,f}
 \end{align*}
where $A_\tau^n$ is the adapted finite variation process and $M_\tau^{n,f}$ is the martingale defined in (\ref{martingn}). Then, for every 
$T > 0$, $\epsilon > 0$, $\eta > 0$, there exists $\Delta>0$ and $n_0\in\mathbb{N}^{\star}$ such that for any sequence $(\tau_n)_{n\in\mathbb{N}}$ of stopping 
times with $\tau_n \leq T$ we have:
\begin{align}
 \label{mark1}
 &\sup_{n\geq n_0} \sup_{\delta\in[0,\Delta]}\mathbb{P}\big( |A_{\tau_n+\delta}^n - A_{\tau_n}^n| \geq \eta  \big)\leq \epsilon, \\
 \label{mark2}
 &\sup_{n\geq n_0} \sup_{\delta\in[0,\Delta]}\mathbb{P}\big( | M^{n,f}_{\tau_n+\delta} -  M^{n,f}_{\tau_n} | \geq \eta  \big)\leq \epsilon. 
\end{align}
\end{itemize}
Let us then check the conditions (A) and (A'). Let $n\in\mathbb{N}^{\star}$. Fix $0<\tau\leq T$ and $\epsilon >0$. Then, for any bounded continuous function $f$ on $\bar{\mathscr{D}}_G$ and any $c_{\epsilon}>0$, we find using the Markov inequality that,
\begin{align*}
 \mathbb{P}\big(|\langle\bar{\mathscr{N}}_\tau^n,f\rangle| \geq  c_{\epsilon}  \big) &\leq \frac{1}{c_{\epsilon}} \|f\|_\infty  \sup_{n\in\mathbb{N}^*} \mathbb{E}\big[\langle\bar{\mathscr{N}}_\tau^n,1\rangle \big]\\
 &\leq  \frac{1}{c_{\epsilon}} \|f\|_\infty  \sup_{n\in\mathbb{N}^*} \mathbb{E}\big[\sup_{\tau\in [0,T]}\langle\bar{\mathscr{N}}_\tau^n,1\rangle \big]
\end{align*}
which clearly, thanks to Lemma \ref{lemboundnorm} with $m=1$ and that $f$ is bounded, completes the verification of (A). Indeed, for all $\epsilon>0$, we can find a compact $B_{\epsilon}$ such that $[-c_{\epsilon}, c_{\epsilon}]\subset B_{\epsilon}$ and 
$\mathbb{P}(\langle\bar{\mathscr{N}}_\tau^n,f\rangle\in B_{\epsilon}^c)<\epsilon$.  We move now to the verification of (A'). We start by $(\ref{mark1})$. First, we note 
from $(\ref{martingn})$ that the finite variation part of the semimartigale  $\langle\bar{\mathscr{N}}_\tau^n,f\rangle$ is given by,
\begin{align*}
A_\tau^n =&\alpha\int_0^{\tau} dt \int_{\bar{\mathscr{D}}_G}\bar{\mathscr{N}}_t^n(dx) \int_{\mathbb{R}^d}f(x+z)k(x,z)dz  
-\beta\int_0^{\tau}dt  \int_{\bar{\mathscr{D}}_G} f(x) \bar{\mathscr{N}}_t^n(dx)   \\
&+\int_0^{\tau}dt\int_{\bar{\mathscr{D}}_G}\Big\{ \int_{\bar{\mathscr{D}}_G}f(y) \mathsf{aff}(x,y) {k}^{\mathsf{af}}(y)dy \Big\}\bar{\mathscr{N}}_t^n(dx).
\end{align*}
Let us consider a stopping time $\tau_n$ and $\delta>0$ such that: $0\leq\tau_n<\tau_n+\delta\leq T$ where $T>0$. Then, straightforward computations using conditions $(i)$ and $(ii)$ of Hypothesis $\ref{H}$, the fact that $f$ is bounded  and then Lemma $\ref{lemboundnorm}$ with $m=1$ give,  
\begin{align*}
 \mathbb{E}\big[|A_{\tau_n+\delta}^n - A_{\tau_n}^n| \big] &\leq (\alpha \gamma_3  + \beta  + {A}_{\mathsf{f}}\gamma_4 )\|f\|_\infty \mathbb{E}\Big[\Big|\int_{\tau_n}^{\tau_n+\delta} dt\int_{\bar{\mathscr{D}}_G}\bar{\mathscr{N}}_t^n(dx)\Big|\Big]\\
 &\leq (\alpha \gamma_3  + \beta  + {A}_{\mathsf{f}}\gamma_4 )\|f\|_\infty \mathbb{E}\Big[\Big|\int_{\tau_n}^{\tau_n+\delta} dt\sup_{t\in [0,T]}\langle\bar{\mathscr{N}}_t^n,1\rangle\Big|\Big]\\
 &\leq \delta(\alpha \gamma_3  + \beta  + {A}_{\mathsf{f}}\gamma_4 )\|f\|_\infty C_{T}.
\end{align*}
Therefore, applying the Markov inequality shows that  $(\ref{mark1})$ holds. It remains now to show $(\ref{mark2})$. Recall that $M^{n,f}$ is the martingale  
given by $(\ref{martingn})$. Then, by Cauchy-Schwarz inequality we find,

\begin{align*}
 \mathbb{E}\big[| M^{n,f}_{\tau_n+\delta} -  M^{n,f}_{\tau_n}| \big]&\leq \mathbb{E}\big[| M^{n,f}_{\tau_n+\delta} -  M^{n,f}_{\tau_n}|^2 \big]^{\frac{1}{2}}\leq 
 \mathbb{E}\big[|[ M^{n,f}]_{\tau_n+\delta} - [ M^{n,f}]_{\tau_n}| \big]^{\frac{1}{2}}.
\end{align*}

Now, using the expression of the quadratic variation  $[ M^{n,f}]$ given by $(\ref{it12})$, the assumptions $(i)$ and $(ii)$ of Hypothesis 
($\mathbb{H}$) together with the fact that $f$ is bounded and Lemma \ref{lemboundnorm} with $m=1$, we establish, following again a straightforward 
computations,  the bound,

\begin{align*}
 \mathbb{E}\big[| M^{n,f}_{\tau_n+\delta} -  M^{n,f}_{\tau_n}| \big] &\leq 
 \mathbb{E}\big[|[ M^{n,f}]_{\tau_n+\delta} - [ M^{n,f}]_{\tau_n}| \big]^{\frac{1}{2}}\\
 &\leq \sqrt{ \frac{\delta}{n}(\alpha \gamma_3  + \beta  + {A}_{\mathsf{f}}\gamma_4 
 ) \|f\|_\infty^2 C_{T}}.
\end{align*}

Applying again the Markov inequality we obtain $(\ref{mark2})$. Thus, condition $(A')$ is also satisfied and we conclude that the sequence  
$(\bar{\mathscr{N}}^n)_{n\in\mathbb{N}^*}$ is tight. \carre

As an immediate consequence of Lemma $\ref{tightness}$, we can extract from $(\bar{\mathscr{N}}^n)_{n\in \mathbb{N}^*}$ a convergent subsequence. Let 
us consider a subsequence $\bar{\mathscr{N}}^{\tilde{n}}$ of $\bar{\mathscr{N}}^{n}$ that converges in distribution to its limit, say $\tilde{\mathscr{N}}$, in the Skorokhod space 
$\mathbb{D}([0,T],\mathcal{S}_{F}(\bar{\mathscr{D}}_G))$. The aim of the following lemma is to overcome the continuity of $\tilde{\mathscr{N}}$.
\begin{lem}
 Let $(\bar{\mathscr{N}}^{\tilde{n}})$ denotes a convergent subsequence of $(\bar{\mathscr{N}}^n)$  and let $\tilde{\mathscr{N}}$ be its limit. Then, the limiting process $\tilde{\mathscr{N}}$ is a.s. continuous which says that,
 \begin{align*}
 \textrm{for almost all}\quad\omega\in\Omega,\quad \tilde{\mathscr{N}}(\omega) \in C([0,T],\mathcal{S}_{F}(\bar{\mathscr{D}}_G)),
 \end{align*} 
  except possibly for some $\omega$ in a set of $\mathbb{P}$-null measure. 
 \label{lemlim}
\end{lem}

\proof  By recalling the definition (\ref{Xi}) of $\mathscr{N}_{\tau-}$, we find that,
\begin{align*}
 | \langle \mathscr{N}_\tau^{\tilde{n}}, 1\rangle - \langle \mathscr{N}_{\tau-}^{\tilde{n}}, 1\rangle |=
\left\{\begin{array}{ll}
 |N_{\tau}^{\tilde{n}} - N_{\tau}^{\tilde{n}}|=0
& \textrm{ if } \tau \notin \cup_i \{\tau_i\} ,
\\
  |N_{\tau}^{\tilde{n}} - N_{\tau_{i-1}}^{\tilde{n}}|\leq 1
 & \textrm{ if } \tau=\tau_i \textrm{ for some } i\geq1.
   \end{array}\right.
\end{align*}
Indeed, the difference between the number of vertices between times $\tau$ and $\tau-$ is at most $1$. Therefore, it is straightforward that, for all 
$f\in C(\bar{\mathscr{D}}_G)$ such that $\|f\|_\infty \leq 1$, 
\begin{align*}
 | \langle \bar{\mathscr{N}}_\tau^{\tilde{n}}, f\rangle - \langle \bar{\mathscr{N}}_{\tau-}^{\tilde{n}}, f\rangle | \leq | \langle \frac{\mathscr{N}_\tau^{\tilde{n}}}{\tilde{n}}, 1\rangle - 
 \langle \frac{\mathscr{N}_{\tau-}^{\tilde{n}}}{\tilde{n}}, 1\rangle | \leq \frac{1}{\tilde{n}}.
\end{align*}
Since the Prokhorov metric is bounded by the total variation metric, we take the supremum over all $\tau\in[0,T]$ to establish
\begin{align*}
 \sup_{\tau\in[0,T]} \pi_{\textrm{P}}( \bar{\mathscr{N}}_\tau^{\tilde{n}},\bar{\mathscr{N}}_{\tau-}^{\tilde{n}}) \leq \sup_{\tau\in[0,T]} \| \bar{\mathscr{N}}_\tau^{\tilde{n}} - 
 \bar{\mathscr{N}}_{\tau-}^{\tilde{n}} \|_{\textrm{TV}} = \sup_{\tau\in[0,T]} \sup_{f\in C(\bar{\mathscr{D}}_G),\|f\|_\infty \leq 1 }
 | \langle \bar{\mathscr{N}}_\tau^{\tilde{n}}, f\rangle - \langle \bar{\mathscr{N}}_{\tau-}^{\tilde{n}}, f\rangle |\leq \frac{1}{\tilde{n}},
\end{align*}
which completes the proof. \carre 

Now, in order to characterize the solution $\xi$ of $(\ref{integrodiff})$, we consider the following function defined, for any 
$\vartheta\in \mathbb{D}([0,T],\mathcal{S}_{F}(\bar{\mathscr{D}}_G))$, by
\begin{align}
\begin{split}
 \Upsilon_\tau(\vartheta)&=\langle\vartheta_\tau,f\rangle -\langle\vartheta_0,f\rangle - \alpha \int_0^\tau dt\,\int_{\bar{\mathscr{D}}_G}\vartheta_t(dx) 
 \int_{\mathbb{R}^d}f(x+z)k(x,z)dz \\ 
&+ \beta \int_0^\tau dt \int_{\bar{\mathscr{D}}_G} f(x)\vartheta_t(dx)   
           - \int_0^\tau dt \int_{\bar{\mathscr{D}}_G}\Big\{ \int_{\bar{\mathscr{D}}_G} f(y)\mathsf{aff}(x,y) {k}^{\mathsf{af}}(y)dy \Big\}\vartheta_t(dx).
\end{split}
          \label{Upsi}
          \end{align}
Thus, if $\Upsilon_\tau(\vartheta)=0$ for all $\tau\geq 0$ and all measurable bounded functions $f$ on $\bar{\mathscr{D}}_G$, then $\vartheta=\xi$, the unique solution of $(\ref{integrodiff})$.  Let us first prove that the function $\vartheta \to \Upsilon_\tau(\vartheta)$ is a.s. continuous at any 
$\vartheta\in C([0,T],\mathcal{S}_{F}(\bar{\mathscr{D}}_G))$.

\begin{lem}
 Suppose the assertions of Hypothesis $\ref{H}$ hold true. Then, for any $\tau\in [0,T]$  and any $f\in C_b(\bar{\mathscr{D}}_G)$, the function $\Upsilon_\tau$ given by (\ref{Upsi}) is continuous from $\mathbb{D}([0,T],\mathcal{S}_{F}(\bar{\mathscr{D}}_G))$ to $\mathbb{R}$ in any $\vartheta\in C([0,T],\mathcal{S}_{F}(\bar{\mathscr{D}}_G))$.
 \label{Upsicont}
\end{lem}

\proof Consider a sequence $\vartheta^n$  that converges to $\vartheta$ in $\mathbb{D}([0,T],\mathcal{S}_{F}(\bar{\mathscr{D}}_G))$ 
w.r.t the Skorokhod topology. Since the limit $\vartheta$ is continuous being in $C([0,T],\mathcal{S}_{F}(\bar{\mathscr{D}}_G))$, we deduce from $(\ref{conv.unif})$  that the sequence $\vartheta^n$ also converges to  $\vartheta$  w.r.t the uniform topology, namely 
\begin{align}
    \sup_{t\in [0,T]}\pi_p(\vartheta^n,\vartheta)\underset{n\rightarrow\infty}{\rightarrow}0.
\label{unif-conv}
\end{align}
Let $\tau\in [0,T]$. From (\ref{Upsi}), using assertions $(i)$ and $(ii)$ of Hypothesis $\ref{H}$ and given that $f$ is bounded, we find that, 
\begin{equation}
    \begin{split}
 |\Upsilon_\tau(\vartheta^n)-\Upsilon_\tau(\vartheta) | &\leq |\langle \vartheta_\tau^n-\vartheta_\tau,f \rangle|+ 
  |\langle \vartheta_0^n-\vartheta_0,f \rangle|\\
  &\qquad\qquad\qquad\quad+ \alpha \int_0^\tau dt\,\bigg|\int_{\bar{\mathscr{D}}_G}[\vartheta_t(dx)-\vartheta^n_t(dx)] 
 \int_{\mathbb{R}^d}f(x+z)k(x,z)dz \bigg|\\ 
&\qquad\qquad\qquad\quad+ \beta \int_0^\tau dt \bigg|\int_{\bar{\mathscr{D}}_G} f(x)[\vartheta_t(dx)-\vartheta^n_t(dx)]\bigg|   \\
          & \qquad\qquad\qquad\quad+ \int_0^\tau dt \bigg|\int_{\bar{\mathscr{D}}_G}\Big\{ \int_{\bar{\mathscr{D}}_G} f(y)\mathsf{aff}(x,y) {k}^{\mathsf{af}}(y)dy \Big\}[\vartheta_t(dx)-\vartheta^n_t(dx)]  \bigg|\\
 &\leq |\langle \vartheta_\tau^n-\vartheta_\tau,f \rangle|+ 
  |\langle \vartheta_0^n-\vartheta_0,f \rangle|\\
  &\qquad\qquad\qquad\quad+ \alpha \tau\sup_{t\in [0,\tau]}\,\bigg|\int_{\bar{\mathscr{D}}_G}[\vartheta_t(dx)-\vartheta^n_t(dx)] 
 \int_{\mathbb{R}^d}f(x+z)k(x,z)dz \bigg|\\ 
&\qquad\qquad\qquad\quad+ \beta\tau\sup_{t\in [0,\tau]}  \bigg|\int_{\bar{\mathscr{D}}_G} f(x)[\vartheta_t(dx)-\vartheta^n_t(dx)]\bigg|   \\
          & \qquad\qquad\qquad\quad+ \tau\sup_{t\in [0,\tau]}  \bigg|\int_{\bar{\mathscr{D}}_G}\Big\{ \int_{\bar{\mathscr{D}}_G} f(y)\mathsf{aff}(x,y) {k}^{\mathsf{af}}(y)dy \Big\}[\vartheta_t(dx)-\vartheta^n_t(dx)]  \bigg| , 
\label{bound.Upsi}  
    \end{split}
\end{equation}
Recall that the function $f$ is continuous and bounded on $\bar{\mathscr{D}}_G$. Thence, by the assertions of Hypothesis \ref{H}, it is easy  to see that the functions $\int_{\mathbb{R}^d}f(x+z)k(x,z)dz$ and $\int_{\bar{\mathscr{D}}_G} f(y)\mathsf{aff}(x,y) {k}^{\mathsf{af}}(y)dy$ are also continuous and bounded in any $x\in\bar{\mathscr{D}}_G$. Therefore, since the weak topology is metrized by the Prokhorov metric, we deduce from $(\ref{unif-conv})$ that the right side of the inequality $(\ref{bound.Upsi})$ goes to zero as $n\rightarrow\infty$. This concludes the proof. \carre

We prove in the next lemma the convergence in distribution of $\Upsilon_\tau(\bar{\mathscr{N}}^{\tilde{n}})$ to the limit $\Upsilon_\tau(\tilde{\mathscr{N}})$, where  
$(\bar{\mathscr{N}}^{\tilde{n}})_{\tilde{n}\in\mathbb{N}^*}$ is again a convergent subsequence extracted from $(\bar{\mathscr{N}}^{{n}})_{{n}\in\mathbb{N}^*}$ and 
$\tilde{\mathscr{N}}$ its limit. To do so, we make use of the  continuous mapping theorem which states that continuous functions are limit-preserving even if their arguments are sequences of random processes \citep[Theorem 2.7]{Bil99}.

\begin{lem}
 Let $\bar{\mathscr{N}}^{\tilde{n}}$ denotes a convergent subsequence of $(\bar{\mathscr{N}}^n)_{n\in \mathbb{N}^*}$ and let $\tilde{\mathscr{N}}$ be its limit. Denote by $\Upsilon_\tau$ 
 the function defined by $(\ref{Upsi})$. Then, the following convergence in law holds, 
 \begin{align}
  \Upsilon_\tau(\bar{\mathscr{N}}^{\tilde{n}}) {\stackrel{Law}{\longrightarrow}}  \Upsilon_\tau(\tilde{\mathscr{N}})  \qquad \textrm{ as  }\quad\tilde{n}\to\infty.
  \label{Uspsiconv}
 \end{align}
 \label{lempreser}
\end{lem}

\proof The subsequence $\bar{\mathscr{N}}^{\tilde{n}}$ is extracting from the sequence $(\bar{\mathscr{N}}^{{n}})_{{n}\in\mathbb{N}^*}$ which is tight (by Lemma \ref{tightness}). Then, this subsequence converges in distribution to its limit $\tilde{\mathscr{N}}$. Furthermore, this limit $\tilde{\mathscr{N}}$ is a.s. continuous (by Lemma \ref{lemlim}). Moreover, the function $\Upsilon_\tau$ is also a.s continuous (by Lemma \ref{Upsicont}), hence it remains just to apply the continuous mapping theorem to get 
the convergence in law given by $(\ref{Uspsiconv})$. This completes the proof. \carre

We proceed now to conclude the proof of Theorem $\ref{theo5}$ (our main result).

\proof First, from Lemma \ref{unik}, the equation $(\ref{integrodiff})$ has a unique solution. Furthermore, Lemma $\ref{tightness}$ states that the sequence 
$\bar{\mathscr{N}}^n$ is tight. Therefore, to conclude the proof of Theorem $\ref{theo5}$, it is suffice to show that $\bar{\mathscr{N}}^n$ has a unique accumulation point which coincides with the unique solution $\xi$  of $(\ref{integrodiff})$. Let again $\bar{\mathscr{N}}^{\tilde{n}}$ be a convergent subsequence and let $\tilde{\mathscr{N}}$ be its 
limit. Then, following $(\ref{Upsi})$, we have to show that a.s. $\Upsilon_\tau(\tilde{\mathscr{N}})=0$ for all $\tau\geq 0$ and any $f\in C_b(\bar{\mathscr{D}}_G)$. Let us prove this. Note from $(\ref{martingn})$ that $\Upsilon_\tau(\bar{\mathscr{N}}^{{n}})=M_\tau^{n,f}$. By taking the expectation and using the quadratic 
variation (\ref{it12}) and Lemma \ref{lemboundnorm}, we get,
\begin{align*}
 \mathbb{E}[|\Upsilon_\tau(\bar{\mathscr{N}}^{{n}})|]=\mathbb{E}[|M_\tau^{n,f}|] &\leq \mathbb{E}[|M_\tau^{n,f}|^2]^{\frac{1}{2}}=\mathbb{E}[\langle M^{n,f}\rangle_\tau]^{\frac{1}{2}}\\
 &\leq \Big(\frac{1}{n}(\alpha \gamma_3  + \beta  + {A}_{\mathsf{f}}\gamma_4 
 ) \|f\|_\infty^2 \mathbb{E} \int_0^\tau \langle \bar{\mathscr{N}}_t^{{n}},1\rangle dt \Big)^{\frac{1}{2}}\\
 &\leq\Big(\frac{1}{n} (\alpha \gamma_3  + \beta  + {A}_{\mathsf{f}}\gamma_4) \|f\|_\infty^2  C_{1,\tau} \Big)^{\frac{1}{2}},
\end{align*}
which tends to zero as $n$ goes to infinity. In addition, we easily check that, for any $\vartheta\in\mathbb{D}([0,T],\mathcal{S}_{F}(\bar{\mathscr{D}}_G))$,
\begin{align*}
 |\Upsilon_\tau(\vartheta)| &\leq \|f\|_\infty(2+\alpha \gamma_3 +\beta+{A}_{\mathsf{f}}\gamma_4 ) \sup_{t\in[0,T]}|\langle \vartheta_t,1 \rangle|.
\end{align*}
Hence using Lemma \ref{lemboundnorm}, we easily verify that the sequence $\Upsilon_\tau(\bar{\mathscr{N}}^{\tilde{n}})_{\tilde{n}}$ is uniformly integrable. Now, by applying the preserving convergence (\ref{Uspsiconv})  together with the dominated convergence theorem it follows that,
\begin{align*}
 \lim_{\tilde{n}\to\infty} \mathbb{E}[|\Upsilon_\tau(\bar{\mathscr{N}}^{\tilde{n}})|]=\mathbb{E}[|\Upsilon_\tau(\tilde{\mathscr{N}})|]=0,
\end{align*}
and consequently $\Upsilon_\tau(\tilde{\mathscr{N}})$ is a.s. equal to zero. Hence, the limit   $\tilde{\mathscr{N}}$ is a.s. equal to $\xi$ the unique solution of (\ref{integrodiff}), which completes the proof.  \carre

Thus,  we have demonstrated a weak convergence in large size ($n\to\infty$) of the stochastic model described in Section \ref{sec2}. Therefore, this convergence opens a door to other complex asymptotic problems where, for instance, it would be interesting to establish a central limit theorem associated with this convergence in order to obtain convergence rates. A thorough analysis of the convergence rates is beyond the scope of this paper and will be the focus of future investigation.

\section{Numerical simulations}
\label{sec6}

We explore now the numerical performance of the Monte Carlo algorithm described in Section \ref{sec2} to provide computational validation and better qualify the contribution. Particularly,  we are interested in understanding the influence of spatial interactions (invitation, affinity, dispersion, etc) on the formation of communities and other aspects of social network dynamics.  We consider the unit square $\bar{\mathscr{D}}_G=[0,1]^2$ as the virtual space. To avoid technicalities arising from irregularities around the 
borders, we consider $[0,1]^2$ as a torus. Therefore, each vertex of the network is characterized by its two spatial coordinates. Furthermore, the invitation kernel $k(x,z)$ is fixed to a  bivariate normal distribution $N_2(0,\sigma^2I_2)$ where $\sigma>0$ (a dispersion parameter) and $I_2$ denotes the $2\times 2$ identity matrix. The affinity kernel $k^{\mathsf{af}}(y)$ follows a bivariate uniform distribution on $[0,1]^2$. This choice for the affinity kernel seems somewhat noninformative since using this kernel gives rise to the fact that new recruited vertices by affinity tend to occupy in a equiprobably manner all the space of the unit square and their dispersion is uniformly in the vicinity of present vertices. Of course, other choices that take into account  the preference for some states in the system can be tested as well but this will be at the expense of additional model complexity. Finally, we consider the triangular linear function  (\ref{interaction}) for the  local affinity function $\mathsf{aff}(\cdot,\cdot)$. According to the proposed approach, we can identify the following parameters $\alpha$, $\beta$, $A_{\mathsf{f}} $, $a_\mathsf{f} $ and $\sigma$ that play a significant role in the tuning of the model.


We then test different affinity functions. Besides the triangular affinity function, the rectangular affinity function $\text{aff}(x,y)=A_f\textbf{1}_{\{||x-y||\leq r_{af}\}}$, we consider the logistic affinity function (Kisd's function) given by $$\text{aff}(x,y)=C_1(1-\frac{1}{1+C_2exp(-k||x-y||)}),$$ 
where $C_1$, $C_2$ and $k$ are some constants. Moreover we consider the Rayleigh fading connection function $$\text{aff}(x,y)=exp(-a(\frac{||x-y||}{b})^c),$$ 
where $a$, $b$ and $c$ are some constants. We shall test the impact of these 4 affinity functions on the dynamics. Examples of these functions are displayed in Fig.\ref{fig:affs}. The Monte Carlo scheme is run for $10^5$ updates with $N_0=1000$, $\alpha=1$, $\beta=1$, $p=0.8$, $A_f=1$, $r_{af}=0.1$ and $\sigma=0.01$ for each affinity function. The results are displayed in Fig.\ref{fig aff12} - Fig.\ref{fig dd4}. We have an increasing trend of size, and the trend is exponential for the logistic and the Rayleigh fading function, which is expected since the expression of both functions contains exponential function. The size increases more rapidly with rectangular function than with the triangular function because we accept more recruitment by affinity in the former case. In the group of exponential trend, the closer a function is to $\text{aff}=A_f=1$, the bigger is the final network size and the more dispersed are the particles. 

The initial particles being distributed uniformly on $\bar{\mathcal{D}}$, the initial degree distribution converges in distribution to a Gaussian distribution. With the triangular function, communities are formed, some of them are related to others, the degree distribution is no more normal but more widespread, particles whose degree is big appear. Compared to the triangular function, the rectangular function accepts more recruitments by affinity (uniform kernel), so we have more particles, communities are more connected, that's why the degree distribution has another small bump of 250-300 degrees. But at the same time, the particles are more dispersed, which shifts the degree distribution to the left. With exponential trend, the size explodes, each particle is connected to many others, that's why the final degree distribution has a shift to the right w.r.t the others. The shift is bigger if the function is closer to $\text{aff}=A_f=1$. With all the four functions, the distribution is multimodal (superposition of several normal distributions) and fat-tailed (there exist many hubs, whose degree far exceeds the mean degree). 


 We run three different scenarios using the Monte Carlo algorithm starting from three random initial states with different sizes. Each initial condition 
 is drawn at $\tau=T_0=0$ from a bivariate uniform distribution on $[0,1]^2$ as follows:
\begin{itemize}
\item A small initial network size starting from $N_0=100$ vertices.
\item A medium initial network size starting from $N_0=1000$ vertices.
\item A high initial network size starting from $N_0=10 000$ vertices.
\end{itemize} 

\begin{figure}[h]
\center
\subfigure[\it $\mathscr{N}_\tau$ at $\tau=0$ with $N_0=10^2$.]
    {\includegraphics[width=5.6cm,height=5cm]
                 {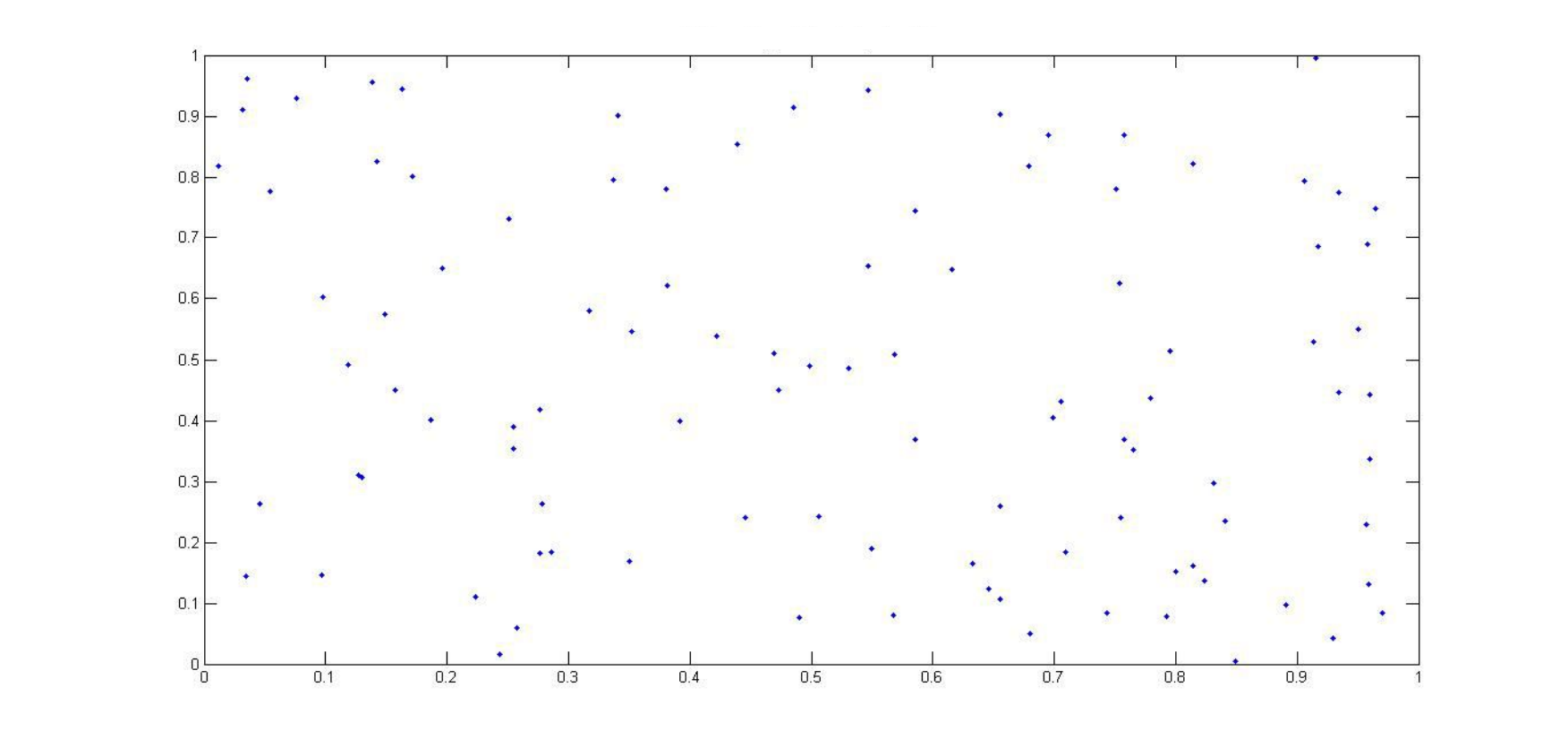}}
\subfigure[\it $\mathscr{N}_\tau$ at $\tau=T_{10^5}$.]
    {\includegraphics[width=5.6cm,height=5cm]
                 {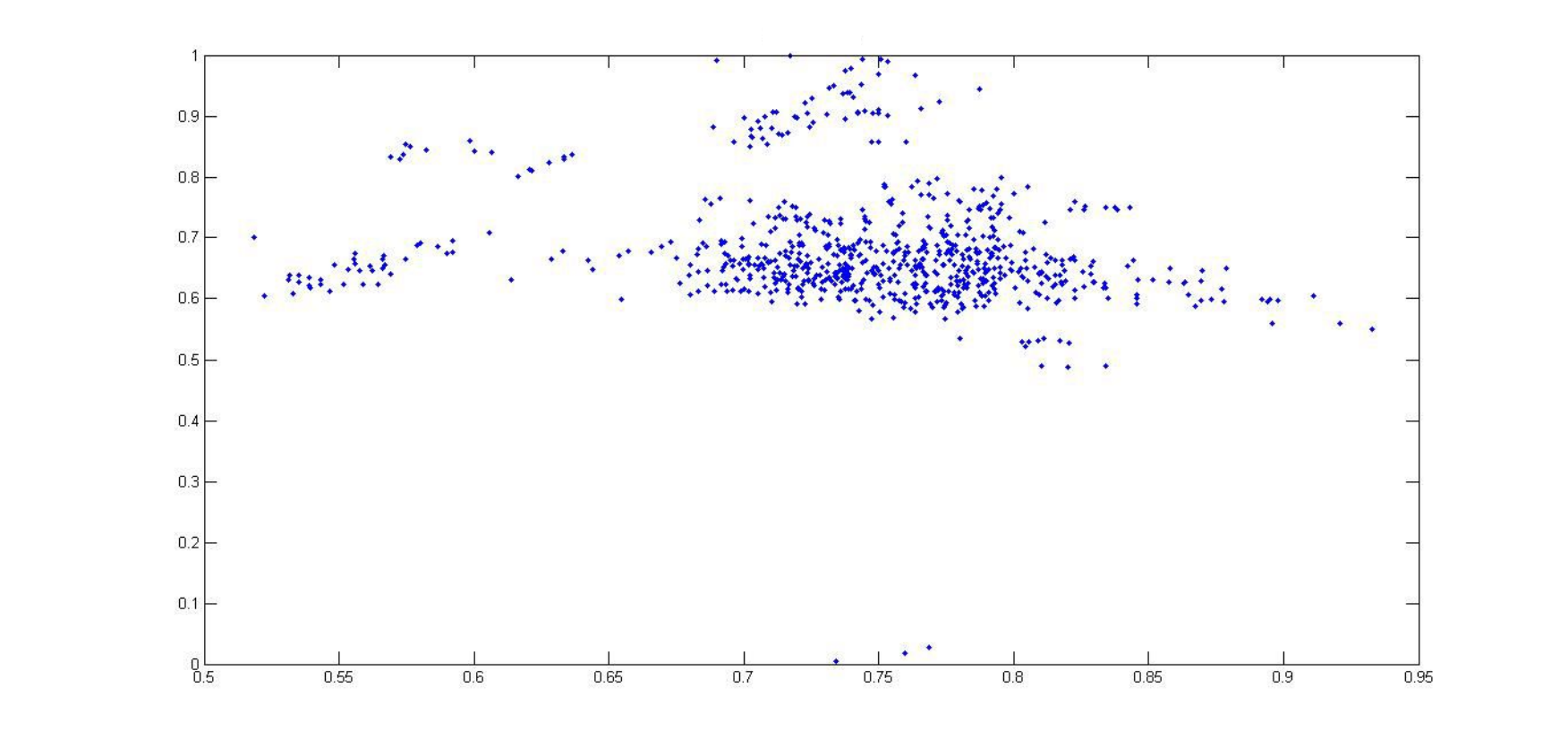}}
\subfigure[\it Size in function of time.]
    {\includegraphics[width=5.6cm,height=5cm]
                 {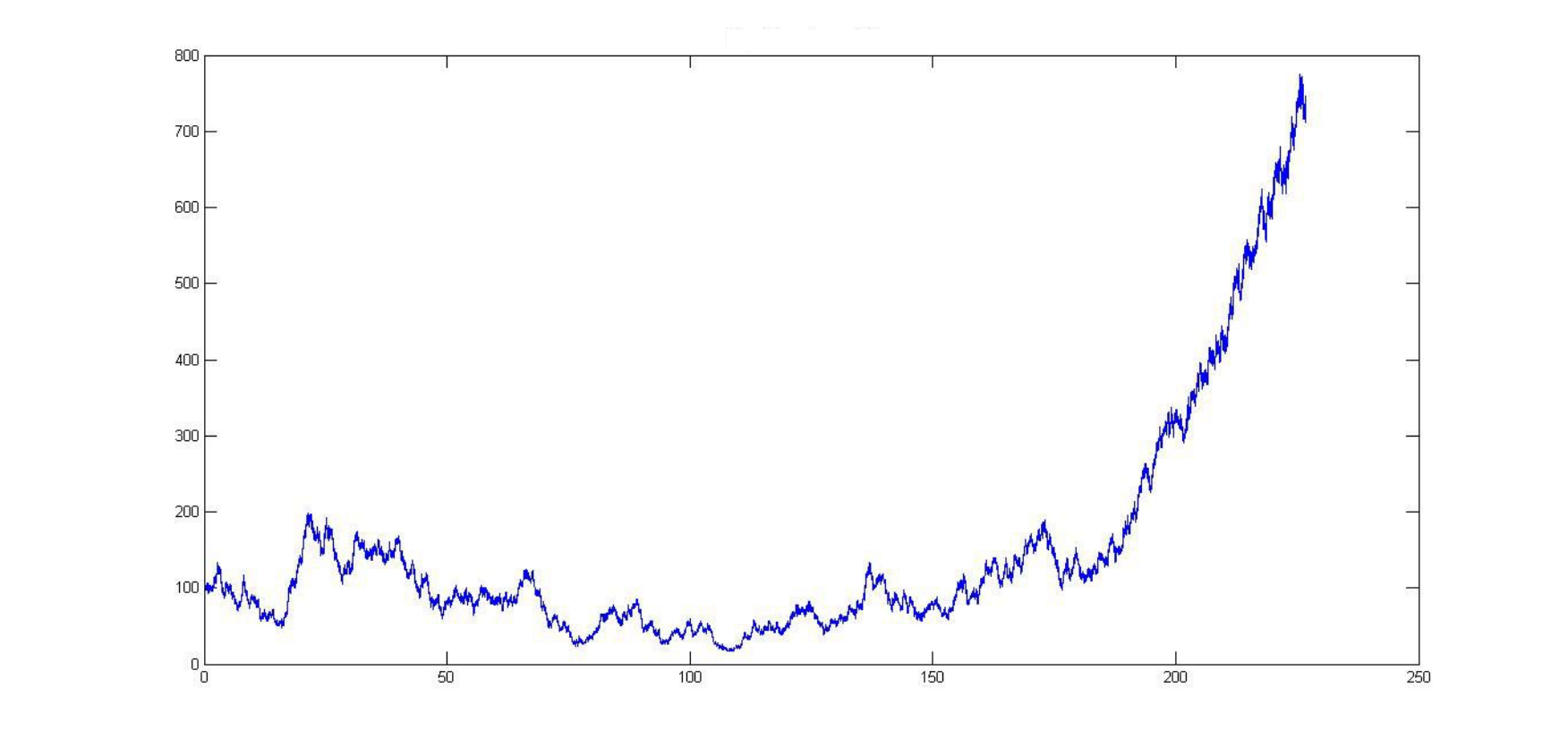}}             
\caption{{\small Figures shown (a) the small initial state $N_0=100$, (b) the final state of the random system $\mathscr{N}_\tau$ and (c) the evolution of the number of vertices in function of time (during $10^5$ updates) using the Monte Carlo scheme with parameters $\alpha=1$, $\beta=1$, $A_{\mathsf{f}}=1$, $a_\mathsf{f}=0.1$ and $\sigma=0.01$.}}
\label{Fig1}
\end{figure}

\begin{figure}
\center
\subfigure[\it $\mathscr{N}_\tau$ at $\tau=0$ with $N_0=10^3$.]
    {\includegraphics[width=5.5cm,height=5cm]
                 {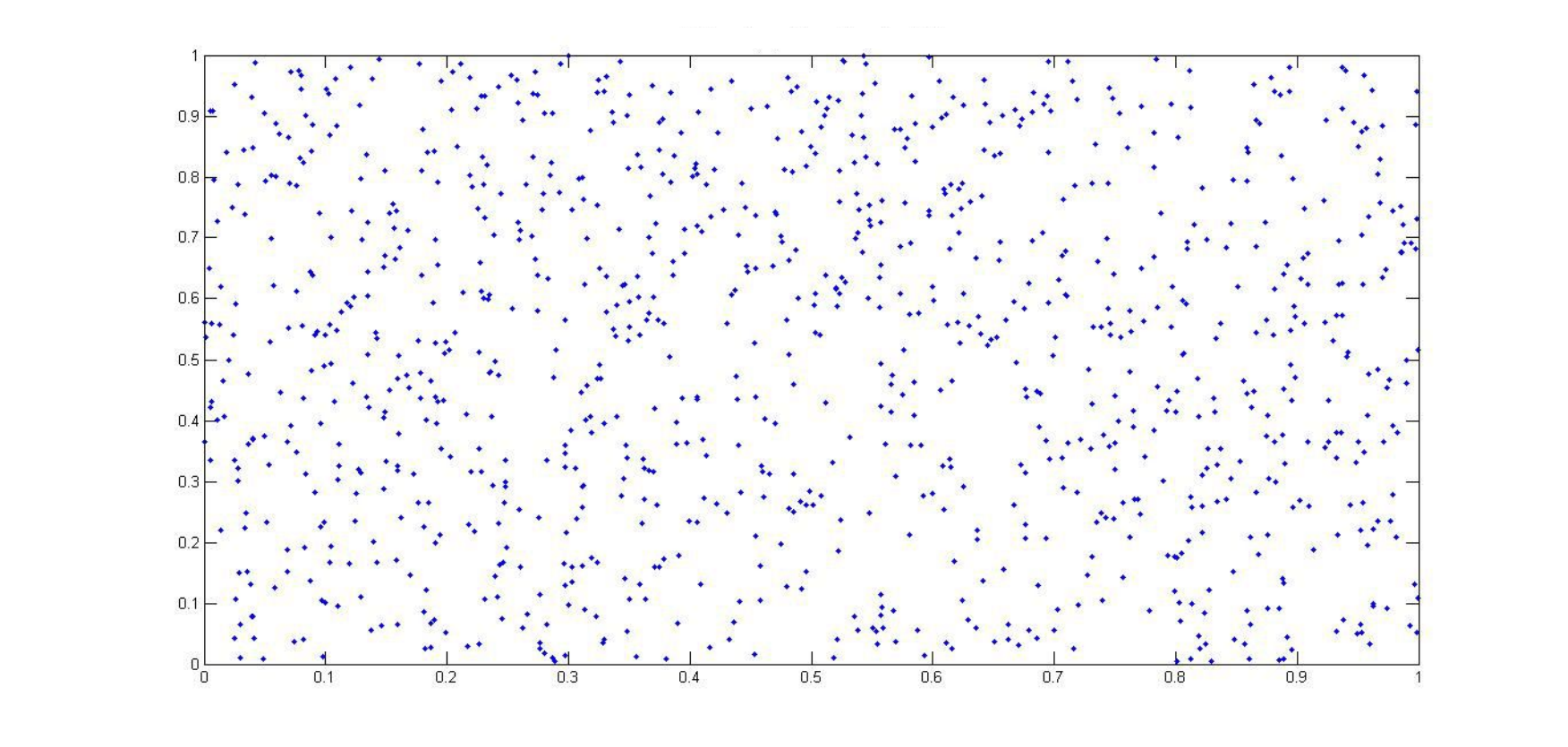}}
\subfigure[\it $\mathscr{N}_\tau$ at $\tau=T_{10^5}$.]
    {\includegraphics[width=5.5cm,height=5cm]
                 {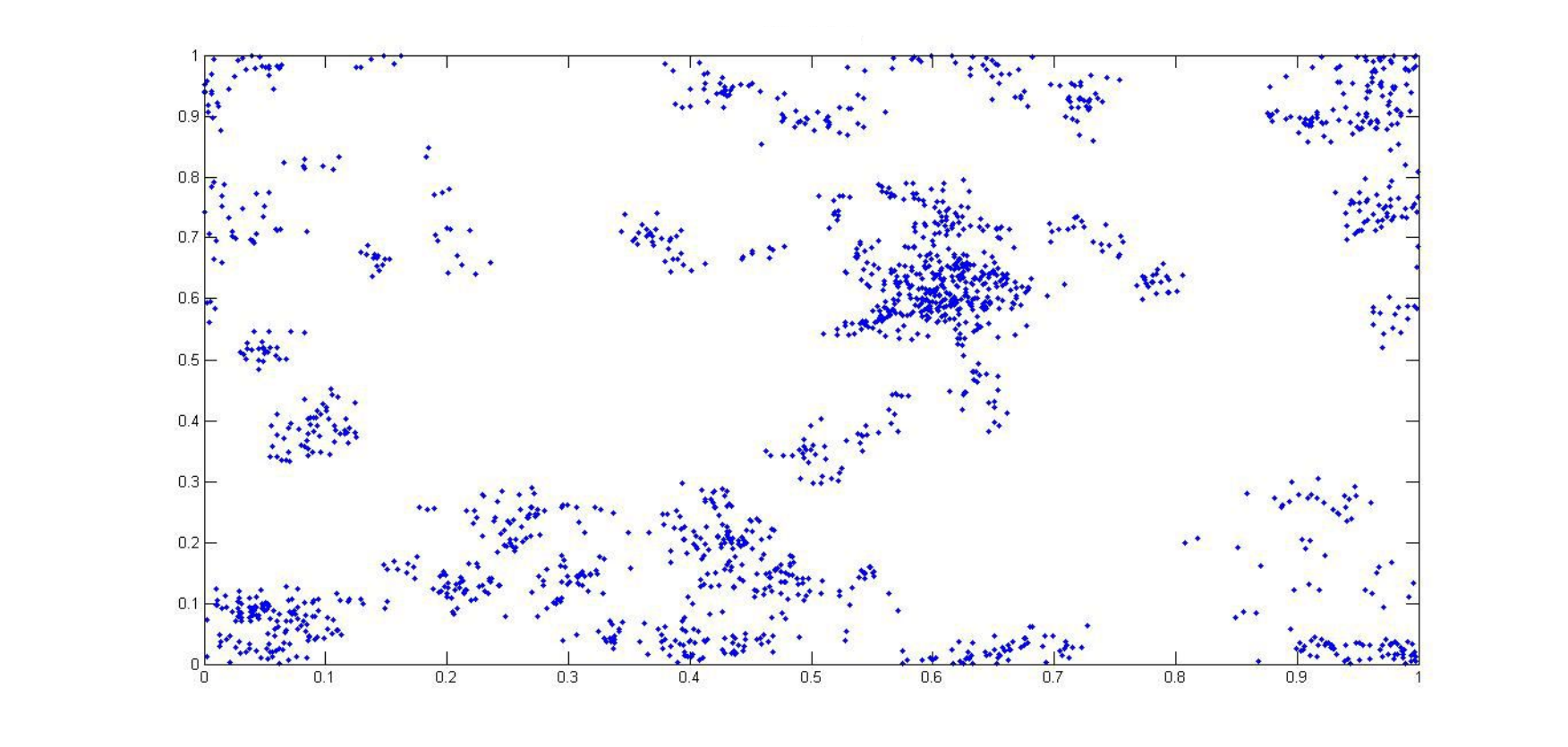}}
\subfigure[\it Size in function of time.]
    {\includegraphics[width=5.5cm,height=5cm]
                 {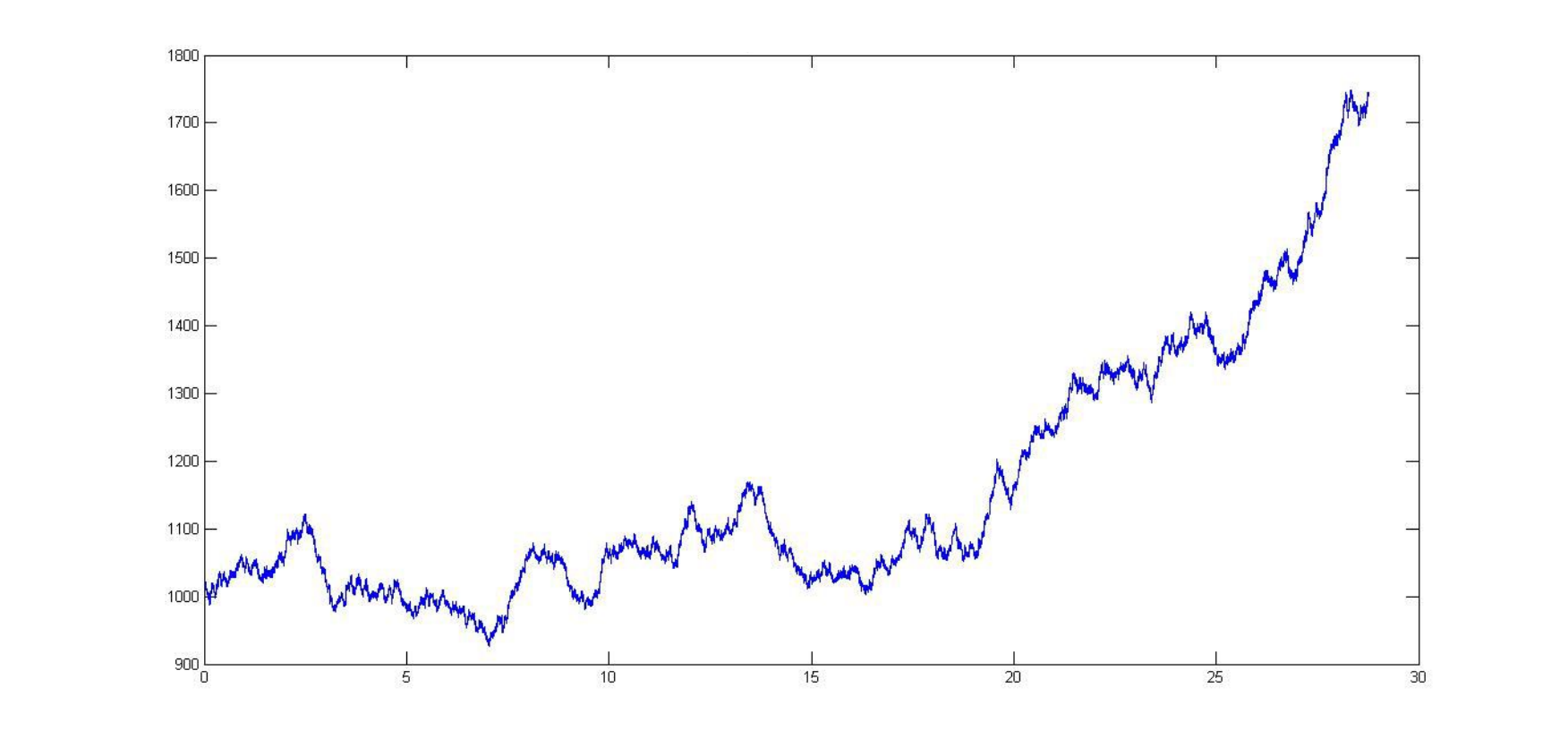}}             
\caption{{\small Figures shown (a) the medium initial state $N_0=1000$, (b) the final state of the random system $\mathscr{N}_\tau$ and (c) the evolution of the number of vertices in function of time (during $10^5$ updates) using the Monte Carlo scheme with parameters $\alpha=1$, $\beta=1$, $A_{\mathsf{f}}=1$, $a_\mathsf{f}=0.1$ and $\sigma=0.01$.}}
\label{Fig2}
\end{figure}

\begin{figure}
\center
\subfigure[\it $\mathscr{N}_\tau$ at $\tau=0$ with $N_0=10^4$.]
    {\includegraphics[width=5.5cm,height=5cm]
                 {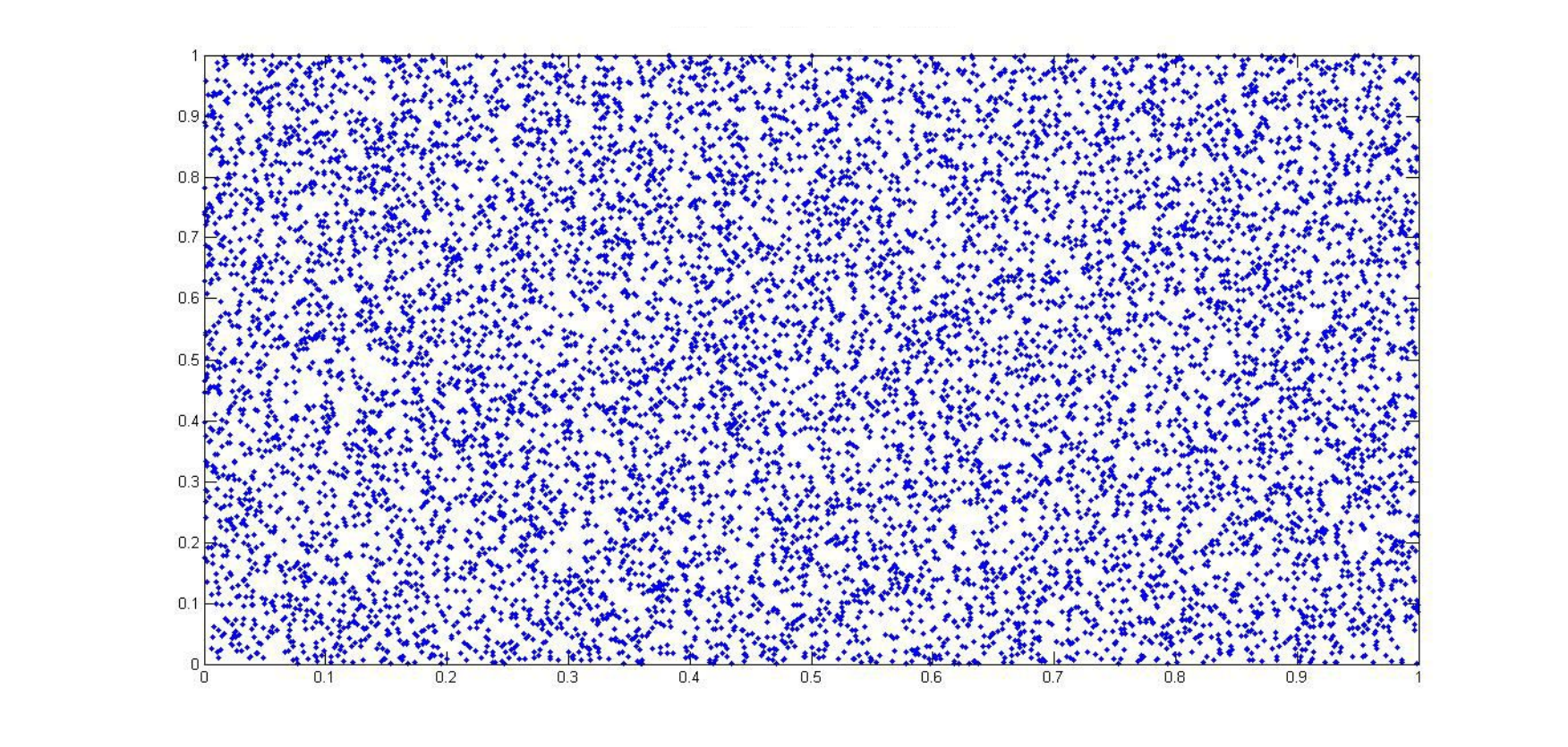}}
\subfigure[\it $\mathscr{N}_\tau$ at $\tau=T_{10^5}$.]
    {\includegraphics[width=5.5cm,height=5cm]
                 {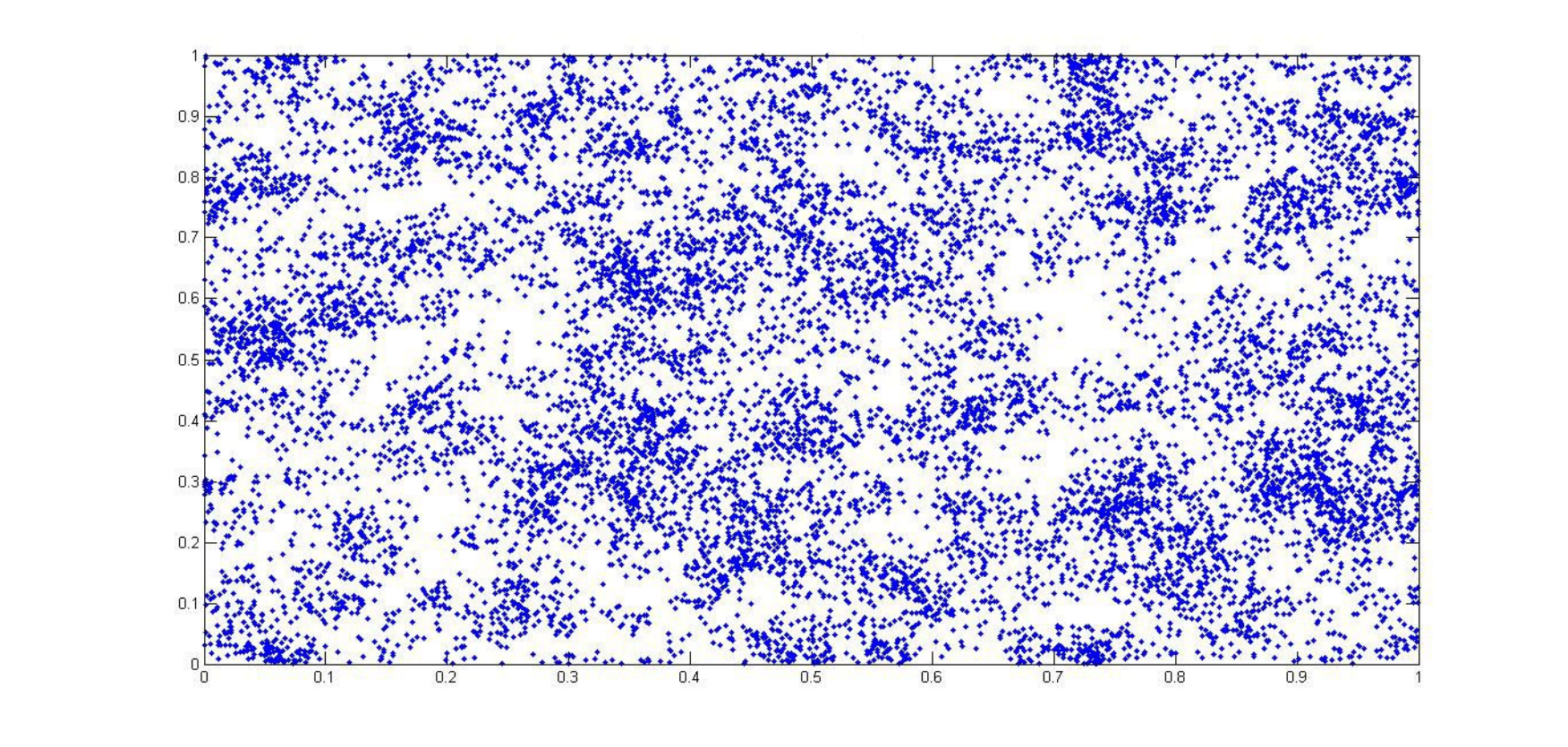}}
\subfigure[\it Size in function of time.]
    {\includegraphics[width=5.5cm,height=5cm]
                 {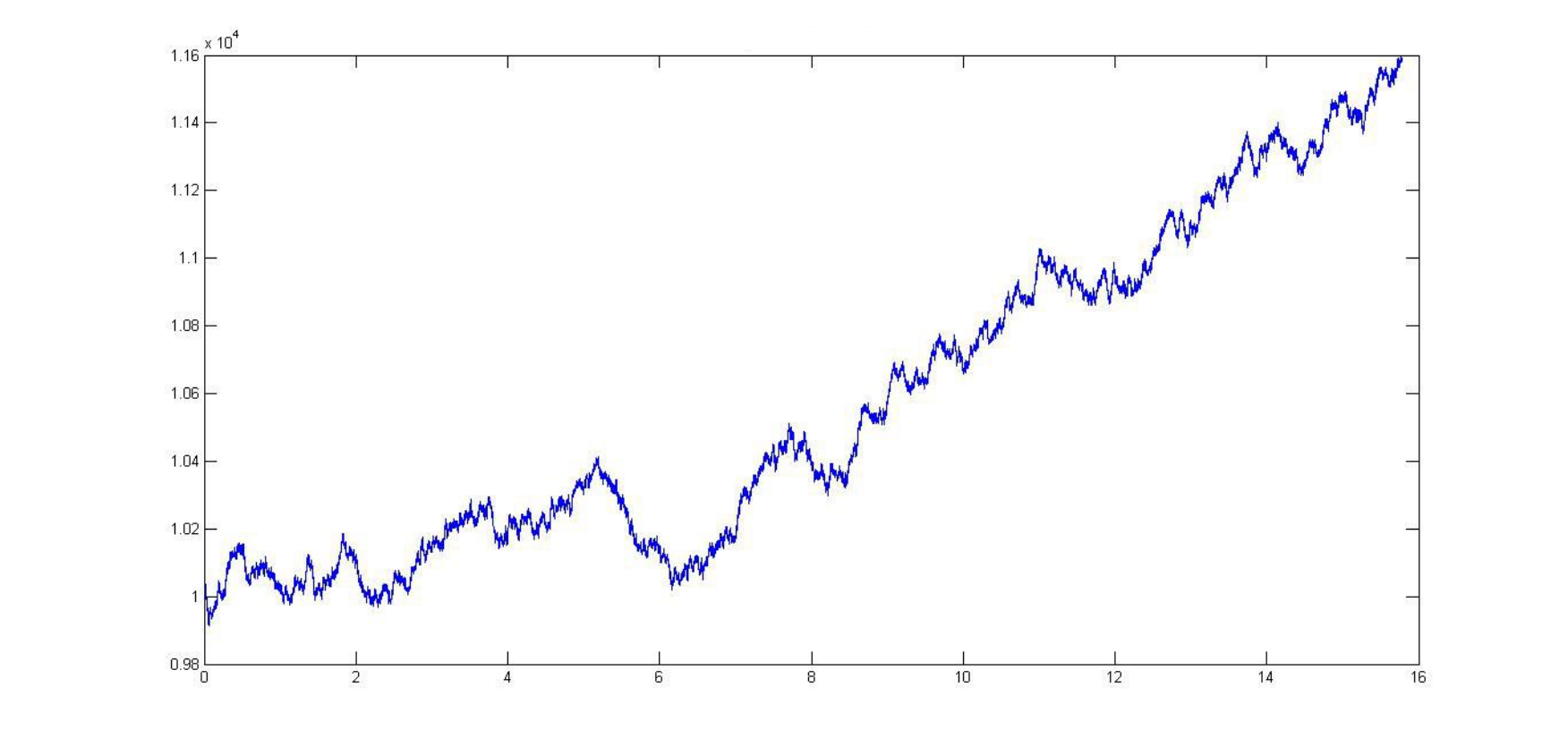}}             
\caption{{\small Figures shown (a) the high initial state $N_0=10000$, (b) the final state of the random system $\mathscr{N}_\tau$ and (c) the evolution of the number of vertices in function of time (during $10^5$ updates) using the Monte Carlo scheme with parameters $\alpha=1$, $\beta=1$, $A_{\mathsf{f}}=1$, $a_\mathsf{f}=0.1$ and $\sigma=0.01$.}}
\label{Fig3}
\end{figure}

From these starting states, we run the numerical scheme for $10^5$ updates with $\alpha=1$, $\beta=1$, $A_{\mathsf{f}}=1$, $a_\mathsf{f}=0.1$ and $\sigma=0.01$. For instance, according to the values taken by the initial size  we aim to test  its effect  on the recruitment potential of the system and to follow the dynamics of the first particles. The systems of small initial size (Figure \ref{Fig1}), of medium initial size (Figure \ref{Fig2}) and of relatively high initial size (Figure \ref{Fig3}) have been simulated during the time interval $\tau\in [0,T_{10^5}]$. The three scenarios with different levels of initial size showed that the dynamics of the simulated systems are dependent on the initial configuration of particles. We report now the results of this numerical experiment and remark that edges are not given here for the clarity of figures.

In the first scenario, the analysis of Figures \ref{Fig1}(a-b) shows that new vertices tend to occupy one particular region and their dispersion is concentrated in the vicinity of some attractive vertices. After $10^5$ updates, we observe clearly the formation of one big community surrounded by three small dispersed communities but the whole system occupies approximately only one half of the unit square. The other half of the unit square appears completely empty where all their first particles present at $\tau =0$ have been disappeared by withdrawal events. We can conclude that only regions with high attractive potential can retain their few particles and recruit new other ones when we start with small size at $\tau =0$. Whereas, if a big region loses their particles at a given time it will be so difficult to recompense after and it remains empty (the two recruitment methods by invitation and affinity appear without any effect in completely empty big region).
As well, we observe in Figure \ref{Fig1}(c) that the network takes certain long time before observing a rapidly growing of its size. This period seems to be due to the period needed for the establishment of  communities. In fact, the communities represent the set of attractive particles that make the network becomes more attractive and then grows up rapidly in size. This behavior is well known by the "preferential attachment mechanism" where vertices prefer linking to the more connected vertices and this phenomenon is responsible of the scale free property observed in several real networks \citep{Bar+Rek99}.

The analysis of Figures \ref{Fig2}(a-b) shows that, in the second scenario, the particles of the system tend to occupy different regions of the space and overall the whole system appears more dispersed than the dynamics seen previously in the first scenario. One remarkable aspect in this numerical test is the formation of several small communities. One approximation may explain the behaviour of the system in Figure \ref{Fig2}(b). This approximation relates to the states of new recruited particles where states with a low local attractiveness of particles have a smaller ability to recruit by invitation or affinity which draws the spatial pattern towards randomness.  This favours the formation of communities. Moreover, as shown in Figure \ref{Fig2}(c), the evolution of the size in function of time shows a more rapidly increasing tendency than the first scenario which is a consequence of the short phase of communities formation.

The analysis of Figures \ref{Fig3}(a-b) shows that, in the third scenario, particles tend  to invade the whole space state $[0,1]^2$. The inter-particle distance is low and the size of particles is very high since the recruitment mechanisms are not negligible now. So, the importance of the initial size appears significant in leading to strong growth of the system density. In Figure  \ref{Fig3}(c), the size of the system was computed at each time for the
count of particles. In fact, as the spatial pattern is clustered and condensate, the count of particles is not over-dispersed with
a greater occurrence of recruitment events expected in simulation of events. These recruitments explain the high number of particles. Further, we can say in this sense that the particle system displays degree assortativity where popular vertices, i.e. those with many edges, are particularly likely to be linked to other popular vertices.

\begin{figure}
\center
\subfigure[\it Histogram of $\mathscr{N}_\tau$ at $\tau=0$.]
    {\includegraphics[width=8cm,height=8cm]
                 {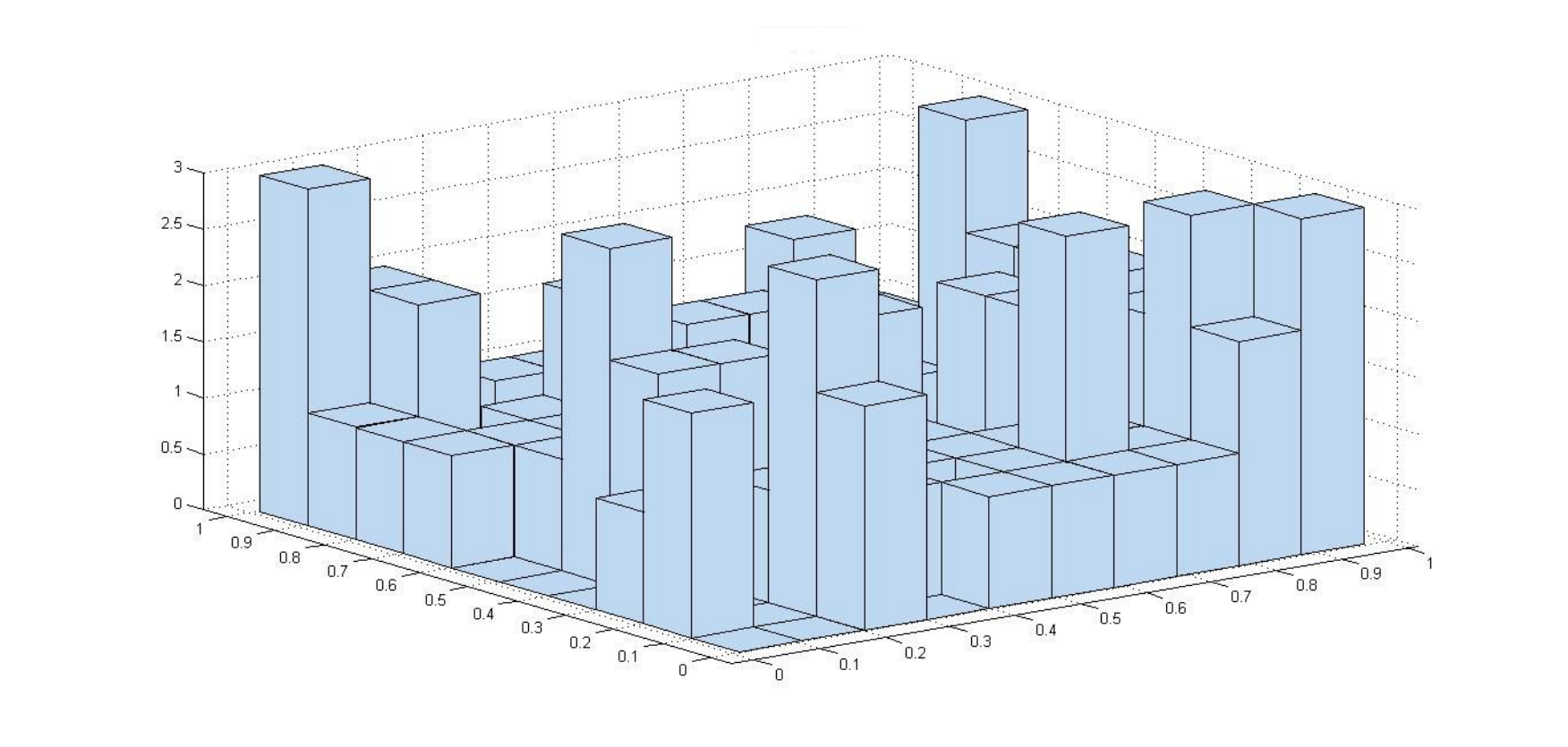}}
\subfigure[\it Histogram of $\mathscr{N}_\tau$ at $\tau=T_{10^5}$.]
    {\includegraphics[width=8cm,height=8cm]
                 {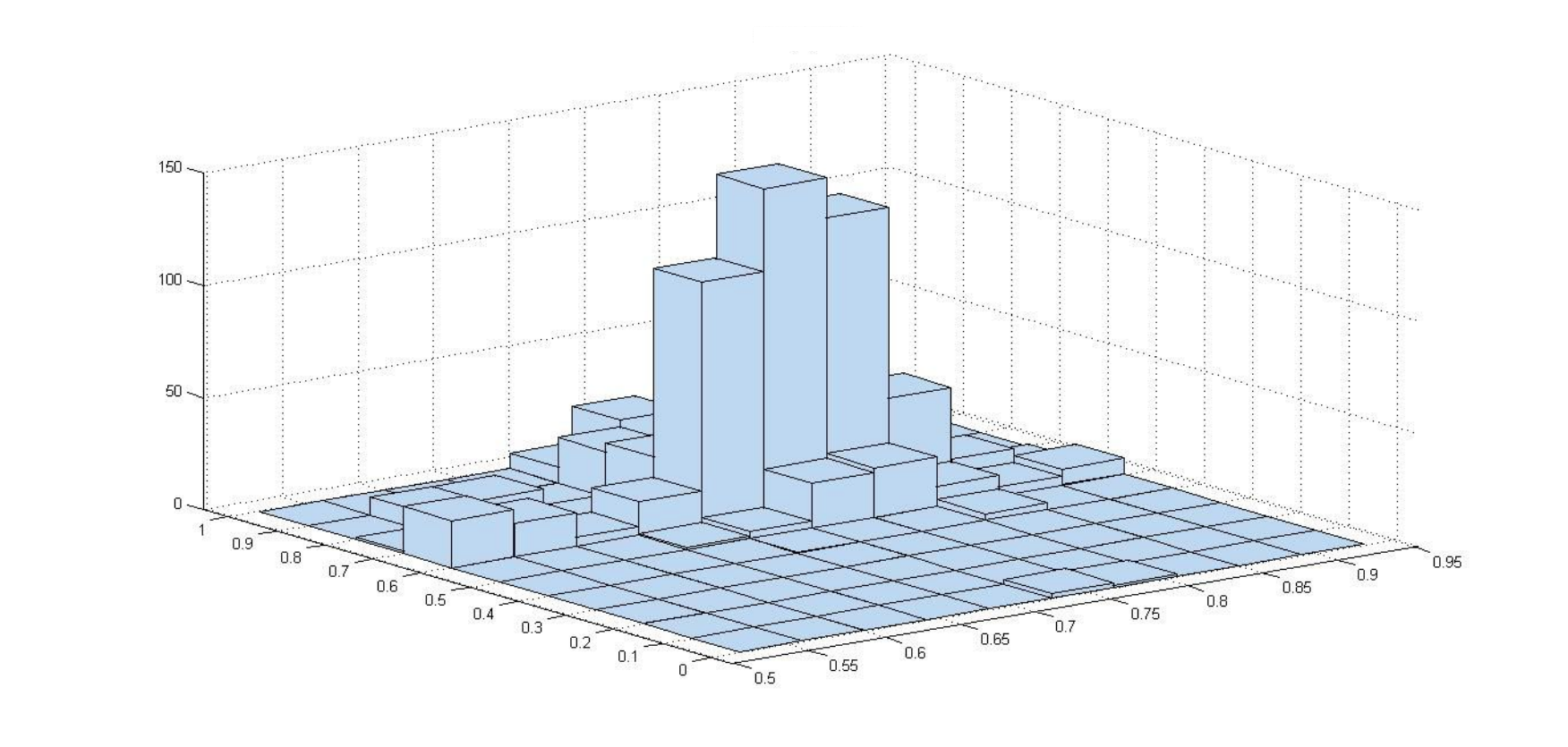}}
             
\caption{{\small Figures shown the histograms of states of vertices distribution  with $N_0=10^2$.}}
\label{Fig4}
\end{figure}

\begin{figure}
\center
\subfigure[\it Histogram of $\mathscr{N}_\tau$ at $\tau=0$.]
    {\includegraphics[width=8cm,height=8cm]
                 {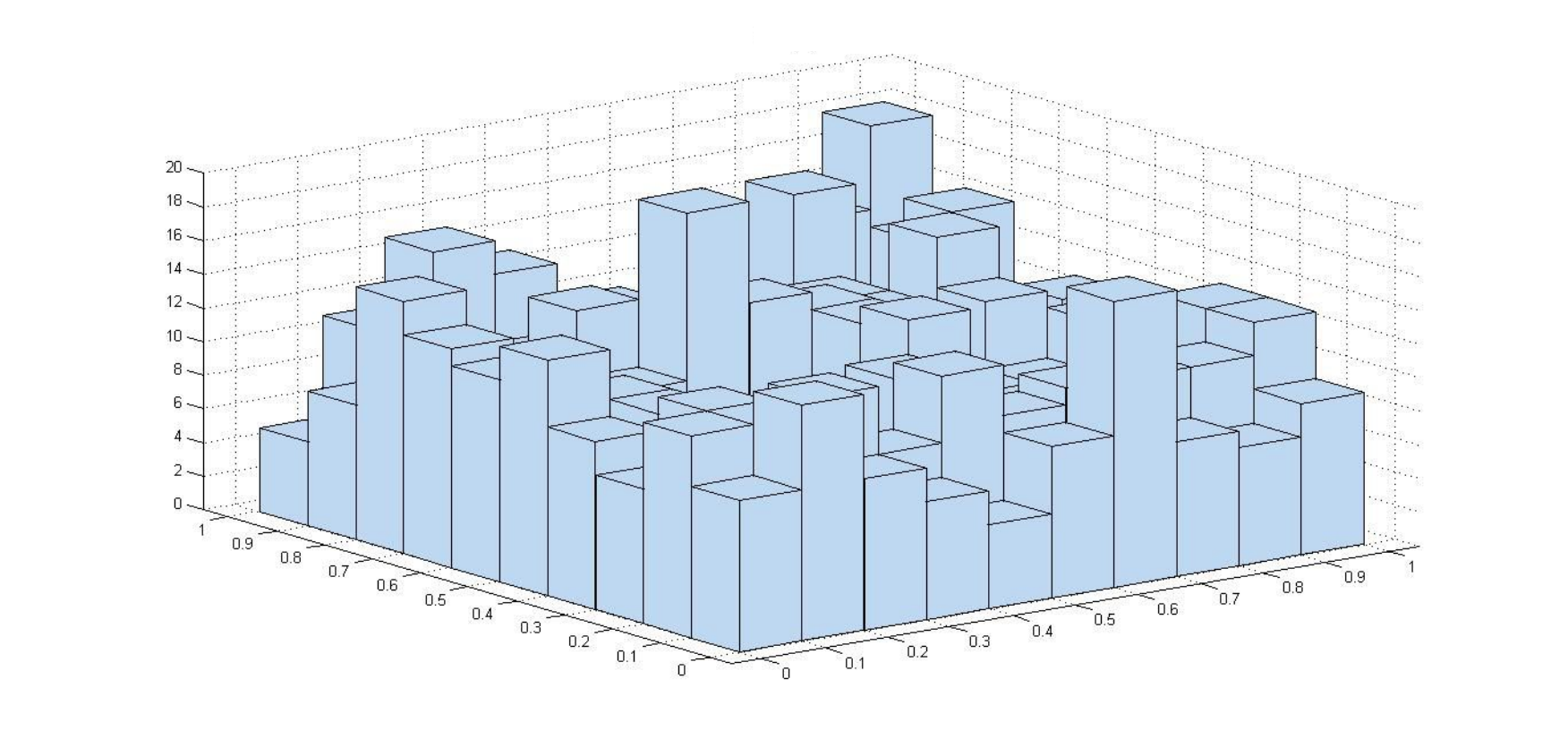}}
\subfigure[\it Histogram of $\mathscr{N}_\tau$ at $\tau=T_{10^5}$.]
    {\includegraphics[width=8cm,height=8cm]
                 {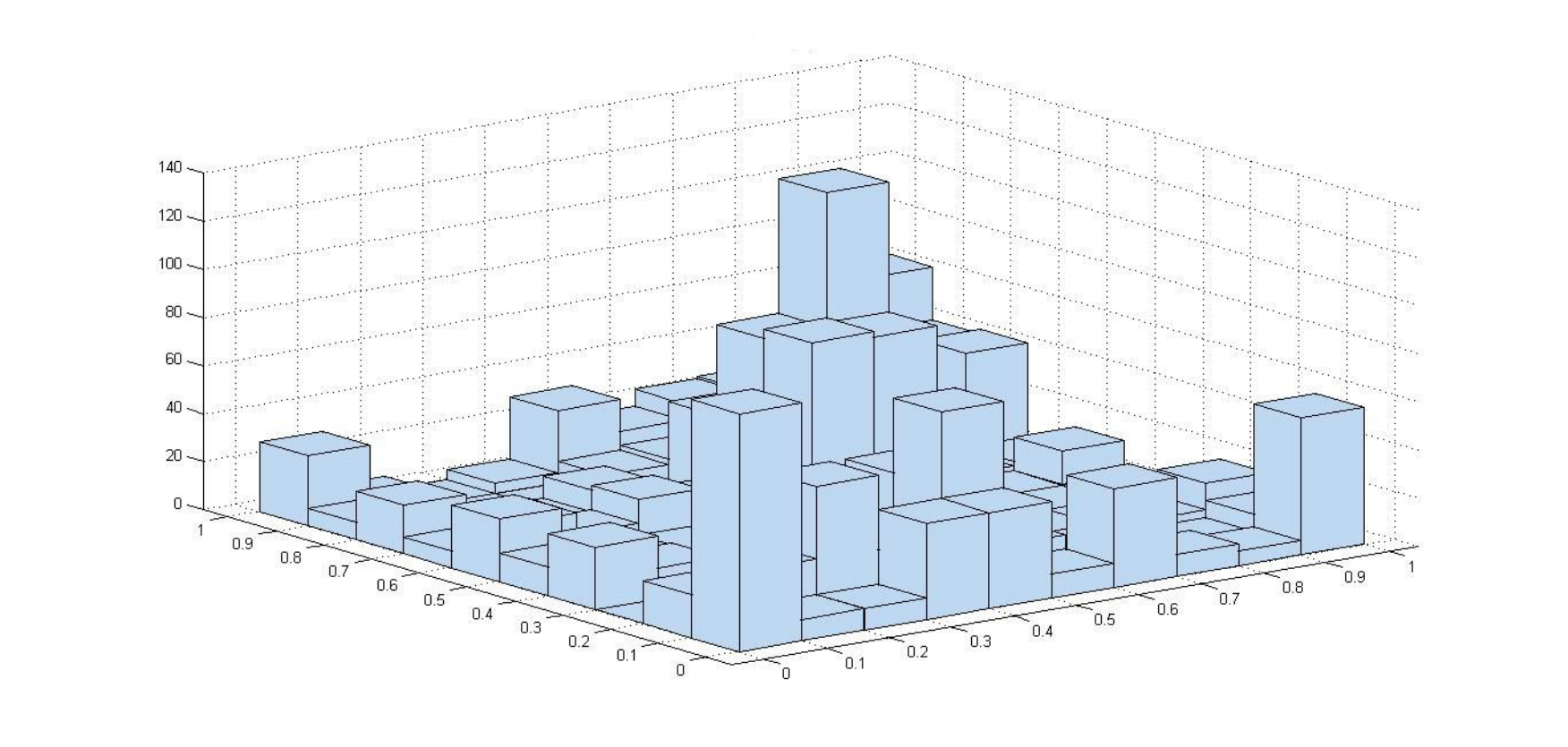}}
             
\caption{{\small Figures shown the histograms of states of vertices distribution  with $N_0=10^3$.}}
\label{Fig5}
\end{figure}

\begin{figure}
\center
\subfigure[\it Histogram of $\mathscr{N}_\tau$ at $\tau=0$.]
    {\includegraphics[width=8cm,height=8cm]
                 {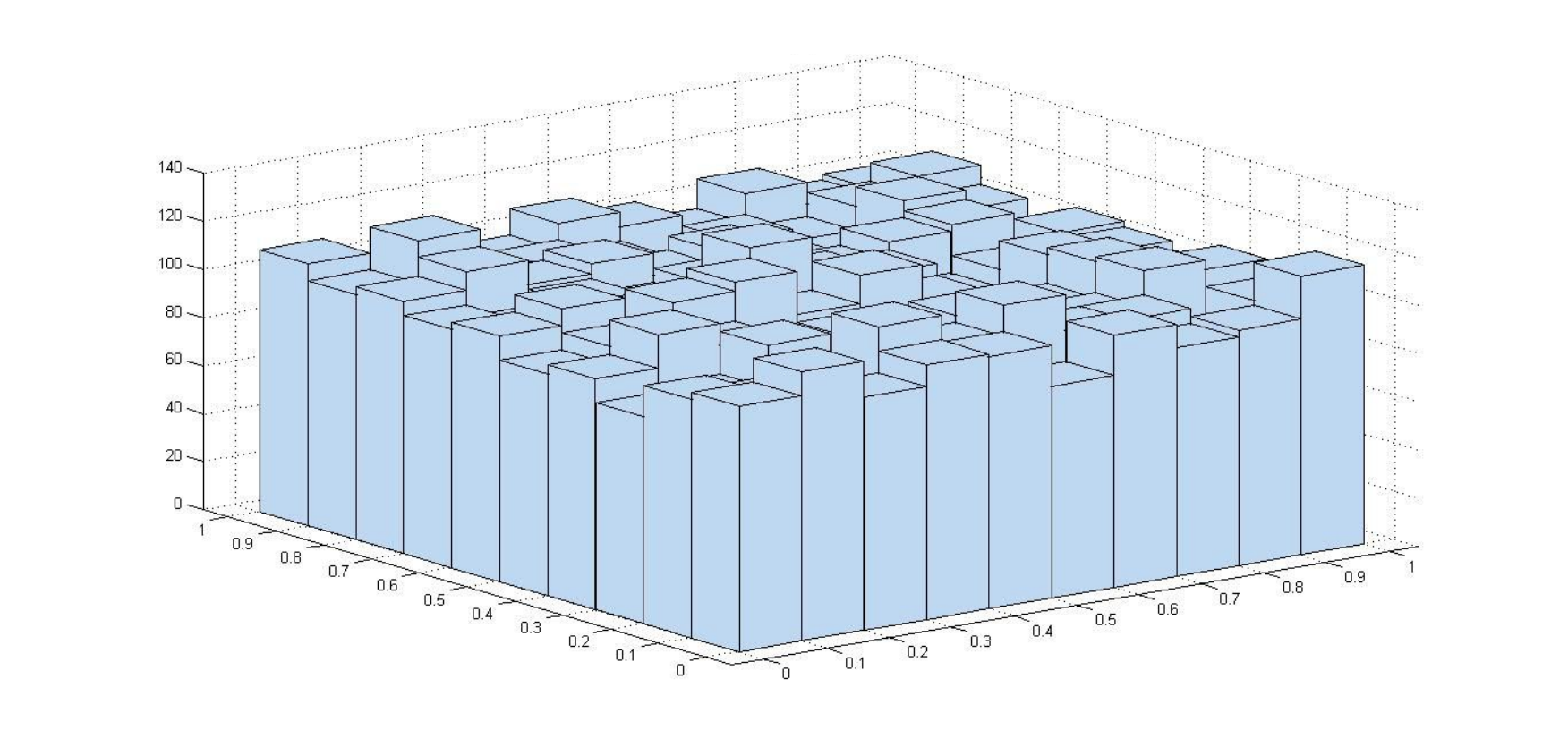}}
\subfigure[\it Histogram of $\mathscr{N}_\tau$ at $\tau=T_{10^5}$.]
    {\includegraphics[width=8cm,height=8cm]
                 {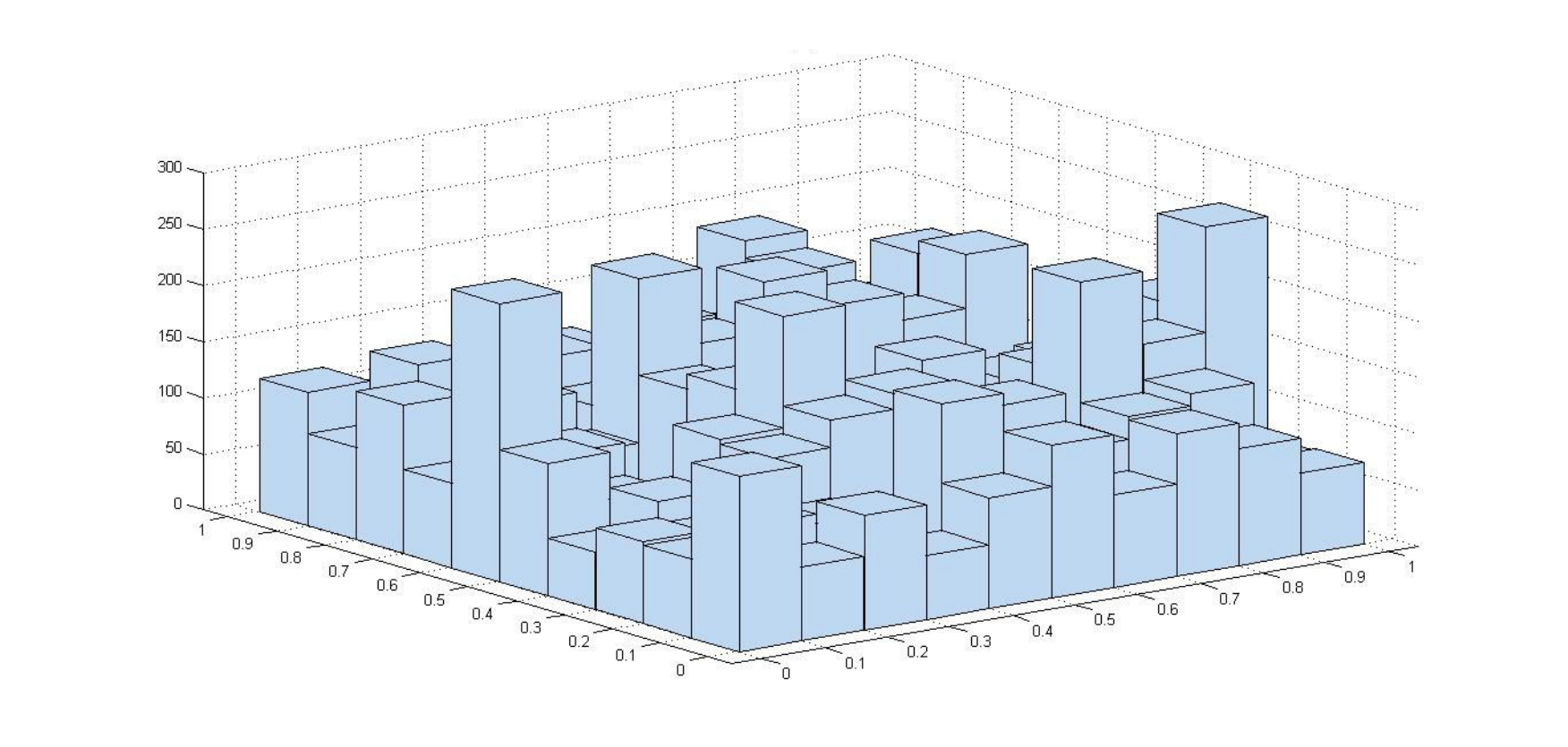}}
             
\caption{{\small Figures shown the histograms of states of vertices distribution  with $N_0=10^4$.}}
\label{Fig6}
\end{figure}

In order to highlight the formation of communities during the dynamics, we plot in Figures \ref{Fig4}, \ref{Fig5} and \ref{Fig6} 
the histograms of the spatial distribution of particles for the three scenarios. As shown, at the initial state $\tau=0$, the particle states are distributed uniformly 
 and after some time communities emerge and the system dynamics evolve through a process of clustering (due to the mechanisms of affinity and dispersion after each recruitment). At this glance, we can conclude that the network based on a measure-valued process designed here resembles real human social networks in a number of fashions. More importantly, it spontaneously demonstrates community structure or clusters of vertices with high modularity.   

\begin{figure}
\center
\subfigure[\it Size in function of time with $N_0=10^2$.]
    {\includegraphics[width=5.5cm,height=5cm]
                 {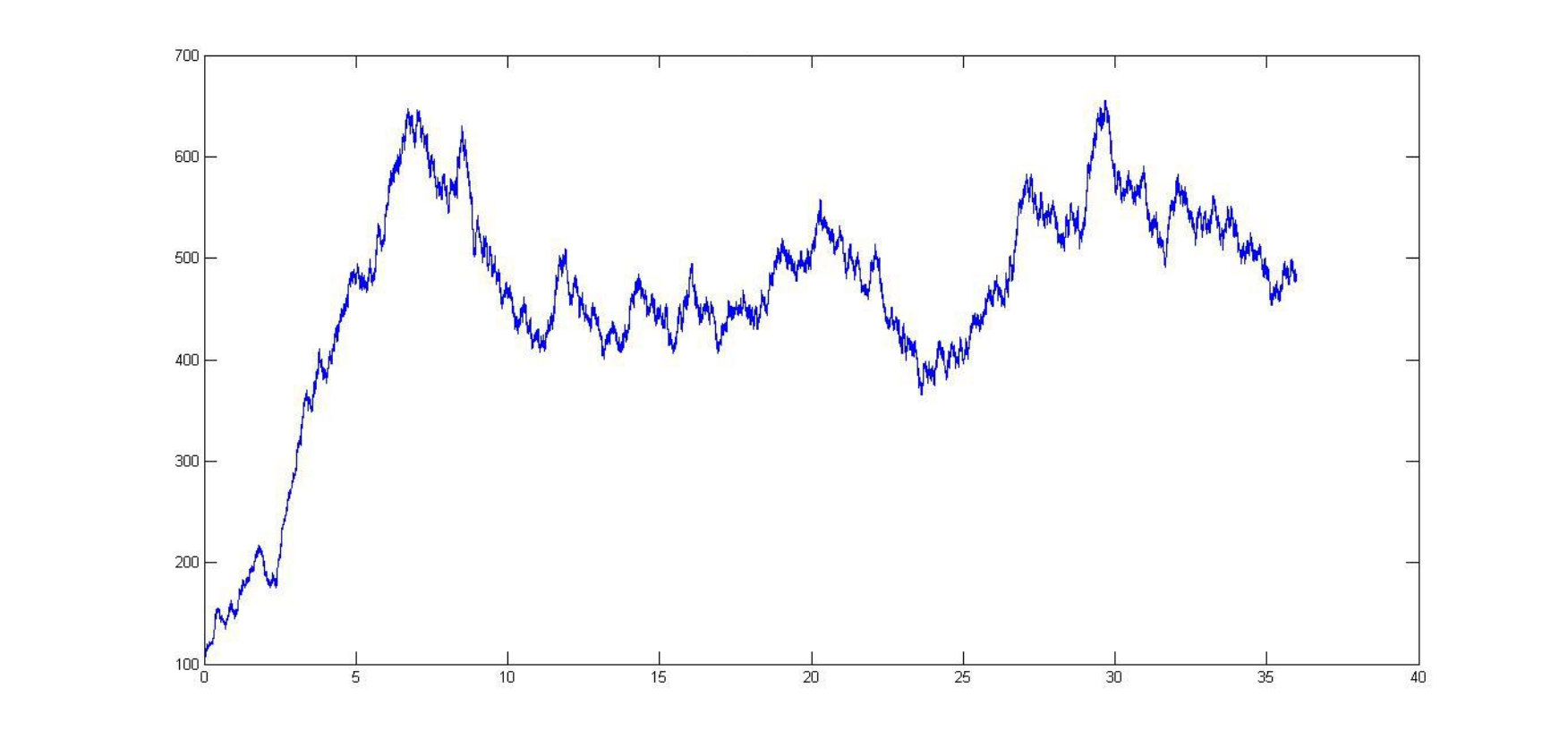}}
\subfigure[\it Size in function of time with $N_0=10^3$.]
    {\includegraphics[width=5.5cm,height=5cm]
                 {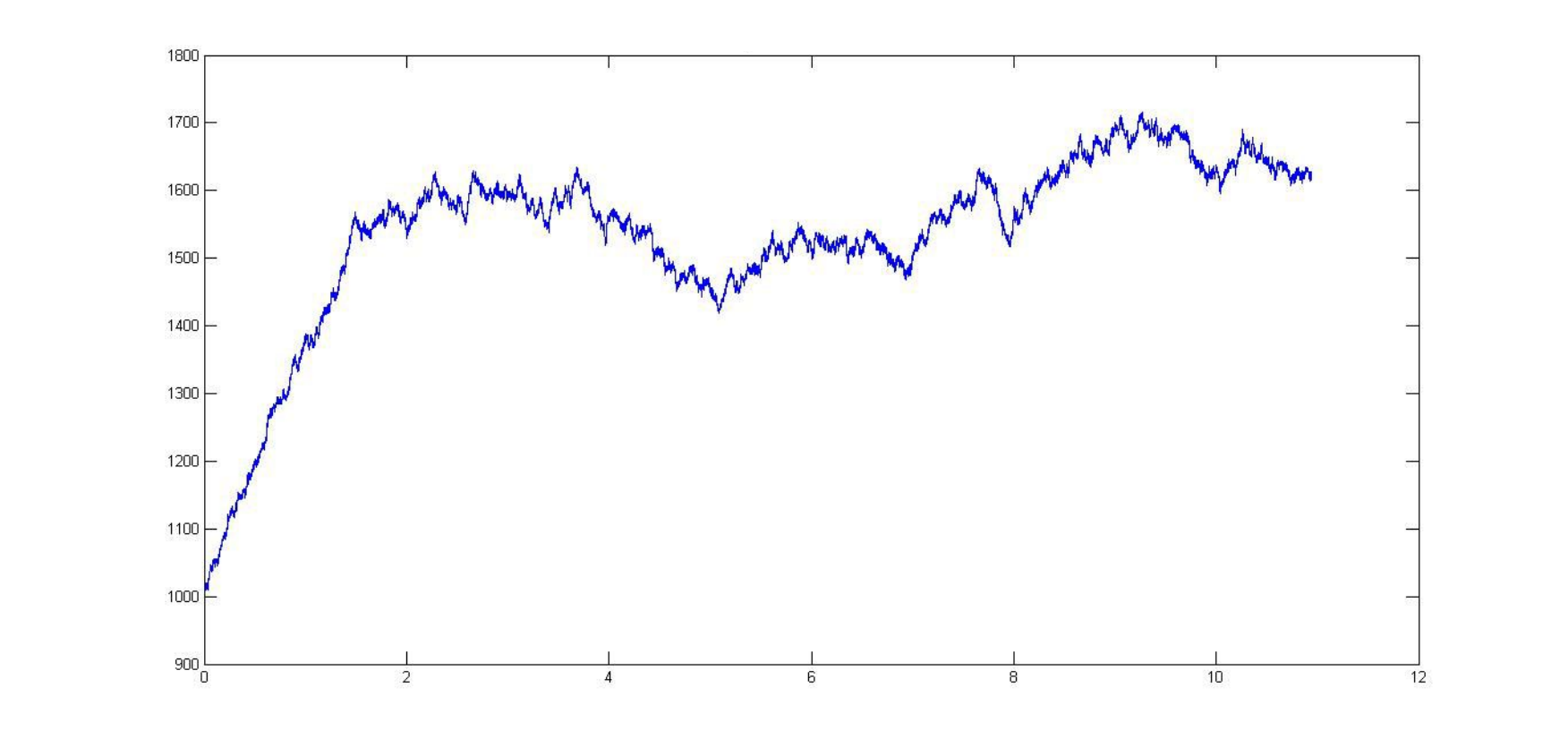}}
\subfigure[\it Size in function of time with $N_0=10^4$.]
    {\includegraphics[width=5.5cm,height=5cm]
                 {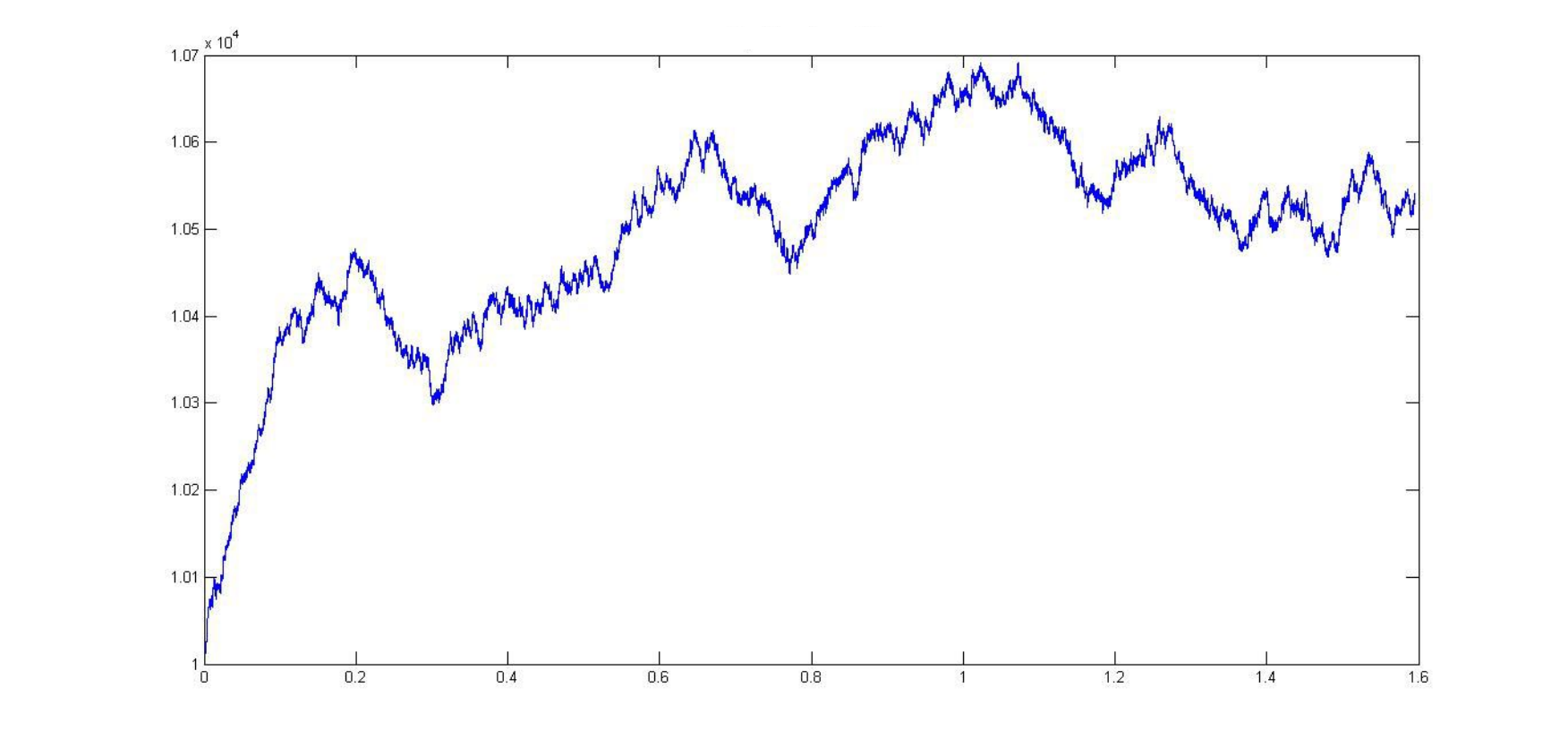}}             
\caption{{\small Figures shown  the time evolution of the number of vertices (during $10^5$ updates) using the model parameter values: $\alpha=2$, $\beta(0)=1.6$, $\epsilon^\star=0.01$, $A_{\mathsf{f}}=2$, $a_\mathsf{f}=0.1$ and $\sigma=0.01$.}}
\label{Fig7}
\end{figure}

Finally, it is interesting to remark that modeling is tuned here by, among others, only rates $\alpha$ and $\beta$ that are space independent and the numerical tests overall showed that after the formation of communities the system size continues to increase because of the increases of attractiveness without reaches some stationary phase as expected. Consequently, the size of the network increases indefinitely during some period before reaches a quasi-stationary state and after that observes an a.s. extinction. Indeed, the system size  grows until a certain equilibrium state and fluctuate temporarily around it. Such  equilibrium state (or full potential state) is no more than the attractive equilibrium of the integro-differential deterministic equation established in Theorem \ref{theo5}. To briefly illustrate this point, we plot in Figure  \ref{Fig7} three trajectories of the system size in function the time where now the numerical test consists of slightly accelerating the  withdrawal rate in a convenient way as follows:  $\beta(T_k)=\beta(T_{k-1})+\epsilon^\star$ for some small positive constant $\epsilon^\star$. Thus, the time increasing withdrawal rate enables us to accelerate the observation of the quasi-stationary states before the ultimate extinction.
Overall, the rates considered here clearly require a thorough analysis and an estimation method for tuning automatically the model as well as how the system is scalable to the massive amounts of data generated today are requested. This will be detailed in a separate paper under preparation.

\section{Discussion}
\label{sec7}

In this work, we investigate a measure-valued process for social networks on a geometric random graph with fixed connectivity threshold $a_\mathsf{f}$ and evolving according to an interacting particle model as follows. Every particle not yet withdrawn is assumed to be capable to invite new particles with positive constant rate $\alpha$. The withdrawal rate $\beta$ is also assumed constant and present particles stay in the system during random exponential times with mean $1/\beta$ during which they interact with their neighbors following a local affinity function of the distance. At the end of the presence period, the particle becomes removed and is no longer considered in the system. The whole system (the set of all particles) can attract and consequently recruits new particles by affinity with a state variable rate $w^{\mathsf{af}}(\cdot,\mathscr{N})$. As we understood, the size $N$ of the system is variable along the continuous time. The particles are related through a random network and are represented by the vertices of an undirected geometric graph. Between two neighbors, we place an edge if their distance is small than some threshold. The graph is nonoriented and an edge between $x$ and $y$ can be seen as two directed edges, one from $x$ to $y$ and the other from $y$ to $x$.

Mean field approximation or large system approximation provides a deterministic equation to describe our system in the asymptotic regime (when the number of particles goes to infinity). Particularly, starting from a random model on finite graph, we derive deterministic equation by increasing the size of the network. The convergence of the continuous-time measure-valued process model to its deterministic limit for large graphs was well established under mild assumptions.

 Nevertheless, this approach creates several open problems with substantial interest in the engineering and computer science communities and more investigations are requested. For instance, an important question addressed partially in this paper is the community detection inside the network. It is of great interest to study this issue.

\bibliographystyle{livre} 
\bibliography{biblio}

\end{document}